\documentclass[a4paper,12pt]{article}
\usepackage{amsfonts,amssymb,amsbsy,amsmath,amsthm,mathrsfs,enumerate,verbatim}
\usepackage{pdfsync}
\usepackage[caption=false,font=footnotesize]{subfig}
\usepackage{grffile}
\usepackage{cite}
\usepackage{color}
 \usepackage{graphicx,color}
\usepackage{algorithmic,algorithm}
\usepackage{breqn}
\usepackage{colortbl}
\usepackage{hhline}

\usepackage{graphics} 
\usepackage{epsfig} 
\usepackage{graphicx}
\usepackage{algorithmic,algorithm}
\usepackage{epstopdf}

\usepackage{amsfonts}
\usepackage{amsmath}
\usepackage{amssymb}
\usepackage{graphicx}
\usepackage{epsf}

\usepackage{mathptmx}       
\usepackage{helvet}         
\usepackage{courier}        
\usepackage{type1cm}        
\usepackage{makeidx}         
\usepackage{graphicx}        
\usepackage{multicol}        
\usepackage[bottom]{footmisc}
\usepackage{amsfonts}
\usepackage{amsmath}
\usepackage{amssymb}
\usepackage{amsfonts}
\usepackage{graphicx}
\usepackage{epsf}
\usepackage{epsfig}

\usepackage{color}
\usepackage{afterpage}
\usepackage{verbatim}
\usepackage{algorithmic,algorithm}

\usepackage{xcolor}
\usepackage{hyperref}

\makeindex             

\begin{document}

\title{An adaptive  finite element/finite difference   domain decomposition method for  applications in microwave imaging}


\author{L.~Beilina   \thanks{Department of Mathematical Sciences, Chalmers University of Technology and University of Gothenburg, SE-42196 Gothenburg, Sweden, e-mail: \texttt{\
  larisa@chalmers.se}}
  \and
 Eric Lindström  \thanks{Department of Mathematical Sciences, Chalmers University of Technology and University of Gothenburg, SE-42196 Gothenburg, Sweden, e-mail: \texttt{\
     erilinds@chalmers.se}}
 \thanks{Journal version of this paper is published in \emph{Electronics} 2022, 11, 1359. https://doi.org/10.3390/electronics/11091359}
  }


\date{}

\maketitle

\begin{abstract}
A new domain decomposition method for Maxwell's equations in
  conductive media is presented. Using this method reconstruction
  algorithms are developed for determination of dielectric
  permittivity function using time-dependent scattered data of electric
  field.  All reconstruction algorithms are based on optimization
  approach to find stationary point of the Lagrangian. Adaptive
  reconstruction algorithms and space-mesh refinement indicators are
  also presented.  Our computational tests show qualitative
  reconstruction of dielectric permittivity function using
  anatomically realistic breast phantom.
\end{abstract}

Keywords: Maxwell's equations; conductive media; microwave imaging; coefficient inverse problem; adaptive finite element method; finite difference method; domain decomposition

 MSC: 65M06; 65M32; 65M55; 65M60


 \graphicspath{
   {./FIGURES/}
   {./FIGURES/TEST2v2/}
   {./FIGURES/gidmesh/}
   {./FIGURES/Eanimation/}
   {./FIGURES/DataTest2/}
   {./FIGURES/DataTest2/PARAVIEW/}
   {./FIGURES/DataTest1/}
   {./FIGURES/MEDIA/}
   {./FIGURES/ForwardP/}
   {./FIGURES/ForwardP/omega40ref3noise10bary/}
   {./FIGURES/ForwardP/omega40ref3noise10rand/}
   {./FIGURES/ForwardP/omega40ref3noise3rand/}
   {./FIGURES/ForwardP/omega40ref3noise3bary/}
   {./FIGURES/TEST2/}
   {./FIGURES/TEST1/}
   {./FIGURES/TEST1/NEW/}
   {./FIGURES/exactepssigma/}
     }

\section{Introduction}

In this work are presented reconstructions algorithms for the problem
of determination of the spatially distributed dielectric permittivity
function in conductive media
using scattered time-dependent data of the electric field at the
boundary of investigated domain.
  Such problems are called Coefficient
Inverse Problems (CIPs). A CIP for a system of time-dependent
Maxwell's equations for electric field is a problem about the
reconstruction of unknown spatially distributed coefficients of this
system  from boundary measurements.



One of the most important application of algorithms of this paper
is microwave imaging including microwave medical imaging and imaging
of improvised explosive devices (IEDs).
 Potential application of algorithms developed in this work are in
 breast cancer detection.  In numerical examples of current paper we
 will focus on microwave medical imaging of realistic breast phantom
 provided by online repository \cite{wisconsin}.
 In this work we
 develop simplified version of reconstruction algorithms which allow
 determine the dielectric permittivity function under the condition
 that the effective conductivity function is known.  Currently we are
 working on the development of similar algorithms for determination of both
 spatially distributed functions, dielectric permittivity and conductivity,
 and we are planning report about obtained results in a near future.


Microwave medical imaging
 is non-invasive imaging.
 Thus, it is very attractive
 addition to the existing imaging technologies like X-ray
 mammography, ultrasound and MRI imaging.
It makes use of the capability of microwaves to
differentiate among tissues based on the contrast in their dielectric
properties.

In  \cite{ieee1} were reported different
 malign-to-normal tissues contrasts, revealing that malign
 tumors have a higher water/liquid content, and thus, higher relative
 permittivity and conductivity values, than normal tissues.
 The challenge is to accurately estimate the relative
permittivity of the internal structures using the information from the
backscattered electromagnetic waves of frequencies around 1 GHz
collected at several detectors.

 Since the 90-s quantitative reconstruction algorithms based on the
 solution of CIPs for Maxwell's system have been developed to provide
 images of the complex permittivity function, see \cite{ieee10} for 2D
 techniques, \cite{ieee15, ruvio1, ieee12,convex1} for 3D techniques in the frequency
 domain and \cite{ieee19, ieee} for time domain (TD)
 techniques.

 In all these works
 microwave medical imaging remained the research field and had 
 little clinical acceptance \cite{meaney1} since the computations are
 inefficient, take too long time,  and produce low contrast values
 for the inside inclusions.
 In all the above cited works
local gradient-based mathematical algorithms
  use
frequency-dependent measurements
which often produce low contrast values of
inclusions and miss small cancerous inclusions.
Moreover, computations in these algorithms are done  often in MATLAB,
sometimes requiring around 40 hours  for solution of inverse problem.

It is well known that CIPs are ill-posed problems \cite{BakKok, T, KSS, itojin}.
Development of non-local  numerical methods is a main
challenge in solution of a such problems.
In works \cite{BK,  BTKM2, TBKF1, TBKF2} was
developed and numerically verified  new non-local approximately
globally convergent  method
for  reconstruction of dielectric permittivity
function.
The two-stage global adaptive optimization method was developed in
\cite{BK} for reconstruction of the dielectric permittivity function.
The
two-stage numerical procedure of \cite{BK} was verified
in several works \cite{BTKM2, TBKF1, TBKF2}
on
experimental data collected by the microwave scattering facility.

The experimental and numerical tests
of above cited works
show that developed methods provide accurate
imaging of all three components of interest in imaging of targets:
shapes, locations and refractive indices of non-conductive media.
 In \cite{convex1}, see also references therein, authors
show reconstruction of complex dielectric permittivity function using
convexification method and frequency-dependent data. 
Potential applications of all above cited works are in the detection and
characterization of improvised explosive devices (IEDs).

The algorithms of the current work can efficiently and accurately
reconstruct the dielectric permittivity function
for one concrete frequency using single measurement data generated by
a plane wave.
 
A such plane wave can be generated by a horn antenna as
it was done in experimental works \cite{BTKM2, TBKF1, TBKF2}.  We are
aware that conventional measurement configuration for detection of
breast cancer consists of antennas placed on the breast skin
\cite{hyperthermi, ruvio1, ruvio2, meaney1, ieee19}. In this work we
use another measurement set-up: we assume that the breast is placed in
a coupling media and then the one component of a time-dependent
electric plane wave is initialized at the boundary of this media. Then
scattered data is collected at the transmitted boundary. This data is
used in reconstruction algorithms developed in this work.  Such
experimental set-up allows avoid multiply measurements and
overdetermination since  we are working with
data resulted   from a single measurement.
An
additional advantage is that in the case of single measurement data
one can use the method of Carleman estimates \cite{klibanov} to prove
the uniqueness of reconstruction of dielectric permittivity
function.

For numerical solution of Maxwell's equations we  have developed finite element/finite difference  domain
decomposition method ( FE/FD DDM).

This approach combines the flexibility of the finite elements and the
efficiency of the finite differences in terms of speed and memory
usage as well as fits the best for reconstruction algorithms of this
paper.
We are unaware of other
works which use similar set-up
  for solution of CIP for
 time-dependent Maxwell's equations
 in conductive media   solved via  FE/FD DDM,
 and this is the first work on this topic.

An outline of the work is as follows: in section \ref{sec:model} we
present the mathematical model and in section \ref{sec:hyb}  we describe the
structure of domain decomposition.  Section
\ref{sec:recalg}  presents reconstruction algorithms including formulation of inverse problem, derivation of finite element and finite difference schemes together with  optimization approach  for solution of inverse problem.
Section \ref{sec:numex}  shows numerical examples of reconstruction of dielectric permittivity function of anatomically realistic breast phantom
 at frequency 6 GHz
of online repository \cite{wisconsin}.  Finally, section \ref{sec:concl}
 discusses obtained results and future research.

\section{The mathematical model}

\label{sec:model}

Our basic model is given in terms of the electric field $E\left( x,t\right)
=\left( E_{1},E_{2},E_{3}\right) \left( x,t\right), x \in
\mathbb{R}^{3}$ changing in the time interval $t \in (0, T)$ under the 
assumption  that the dimensionless relative magnetic permeability of 
the medium is $\mu_r \equiv
  1$.  We consider the Cauchy problem for the Maxwell equations for
  electric field $E\left( x,t\right)$, further assuming that that 
the  electric  volume
charges are equal zero,  to get  the model equation
for $x \in \mathbb{R}^{3},t\in (0,T]$.
\begin{equation}\label{E_gauge}
\begin{split}
 \frac{1}{c^2} \varepsilon_r \frac{\partial^2 E}{\partial t^2} +  \nabla \times \nabla \times E  &= - \mu_0 \sigma \frac{\partial E}{\partial t}, \\
  \nabla \cdot(\varepsilon E) &= 0,  \\
  E(\cdot,0) &= f_0, ~~~\frac{\partial E}{\partial t}(\cdot,0) = f_1.
\end{split}
\end{equation}
Here, $\varepsilon_r(x) = \varepsilon(x)/\varepsilon_0$  is
the dimensionless relative dielectric permittivity 
and
$\sigma(x)$ is the 
effective conductivity function, $\varepsilon_0,
\mu_0$ are the permittivity and permeability of the free space,
respectively, and $c=1/\sqrt{\varepsilon_0 \mu_0}$ is the speed of
light in free space.

We are not able numerically solve the problem (\ref{E_gauge}) in the
unbounded domain, and thus we introduce a convex bounded domain
$\Omega\subset \mathbb{R}^{3}$ with boundary $\partial \Omega$.  For numerical
solution of the problem (\ref{E_gauge}), a domain decomposition finite
element/finite difference method is developed
and summarized in Algorithm 1 of section \ref{sec:hyb}.

A domain decomposition means that we
divide the computational domain $\Omega$ into two subregions,
$\Omega_{\rm FEM}$ and $\Omega_{\rm FDM}$ such that $\Omega = \Omega_{\rm
  FEM} \cup \Omega_{\rm FDM}$ with $ \Omega_{\rm FEM} \subset \Omega$,
see Figure \ref{fig:F3}.
Moreover, we will additionally decompose the domain $
\Omega_{\rm FEM} = \Omega_{\rm IN} \cup \Omega_{\rm OUT}$ with $
\Omega_{\rm IN} \subset \Omega_{\rm FEM}$ such that functions
$\varepsilon_r(x)$ and $\sigma(x)$ of equation (\ref{E_gauge}) should
be determined only in $\Omega_{\rm IN}$, see Figure \ref{fig:F3}.
When solving the inverse problem \textbf{IP} this assumption allows
stable computation of the unknown functions $\varepsilon_r(x)$ and
$\sigma(x)$ even if they have large discontinuities in
$\Omega_{\rm FEM}$.

The communication between $\Omega_{\rm FEM}$ and
$\Omega_{\rm FDM}$ is arranged using a mesh overlapping through a two-element
thick layer around $\Omega_{\rm FEM}$, see  elements in blue color in
Figure \ref{fig:F2}-a),b).  This layer consists of triangles  in $\mathbb{R}^2$
or tetrahedrons in $\mathbb{R}^3$  for
$\Omega_{\rm FEM}$,  and of squares in $\mathbb{R}^2$ or cubes in $\mathbb{R}^3$
for $\Omega_{\rm FDM}$.

The key idea with such a domain decomposition is to apply different
numerical methods in different computational domains.  For the
numerical solution of \eqref{E_gauge} in $\Omega_{\rm FDM}$ we use the
finite difference method on a structured mesh. In $\Omega_{\rm FEM}$, we
use finite elements on a sequence of unstructured meshes $K_h =
\{K\}$, with elements $K$ consisting of
tetrahedron's in $\mathbb{R}^3$ satisfying minimal angle condition \cite{KN}.


We   assume in this paper that for some known constants $d_1>1, d_2>0$, the functions $\varepsilon_r(x)$ and $\sigma(x)$ of equation
 (\ref{E_gauge}) satisfy
\begin{equation} \label{2.3}
\begin{split}
  \varepsilon_r(x) &\in \left[ 1,d_1\right],\quad
  \sigma(x) \in \left[0,d_2\right], ~\text{ for }x\in  \Omega _{\rm IN}, \\
  ~  ~ \varepsilon_r(x) &=1, \quad \sigma(x) = 0\quad  \text{ for }x\in  \Omega _{\rm FDM} , ~~
\varepsilon_r(x), \sigma(x) \in C^{2}\left( \mathbb{R}^{3}\right). 
\end{split}
\end{equation}

Turning to the boundary conditions at $\partial \Omega$, we use the fact that
(\ref{2.3}) and (\ref{E_gauge}) imply that since $\varepsilon_r(x) =
1, \sigma(x)=0$ for $x\in \Omega_{\rm FDM} \cup \Omega_{\rm OUT},$ then a well
known transformation 
\begin{equation}\label{divfree}
\nabla \times \nabla \times E = \nabla (\nabla
\cdot E) - \nabla \cdot ( \nabla E)
\end{equation}
 makes the equations
(\ref{E_gauge}) independent on each other in $\Omega_{\rm FDM} $, and thus, in
$\Omega_{\rm FDM}$ we need to solve the equation
\begin{equation} \label{6.11}
  \frac{\partial^2 E}{\partial t^2}-\Delta E= 0, ~~~(x,t) \in  \Omega_{\rm FDM} \times ( 0,T]. 
\end{equation}
We write $\partial \Omega =\partial \Omega_1 \cup \partial \Omega_2 \cup \partial \Omega_{3}$, meaning
that $\partial \Omega_1$ and $\partial \Omega_2$ are the top and bottom sides of the
domain $\Omega$, while $\partial \Omega_{3}$ is the rest of the
boundary. Because of (\ref{6.11}), it seems natural to impose first
order absorbing boundary condition for the wave equation \cite{EM},
\begin{equation}
\frac{\partial E}{\partial n} + \frac{\partial E}{\partial t} =0,\left( x,t\right) \in \partial \Omega \times ( 0,T].  \label{6.12a}
\end{equation}
Here, we denote the outer normal derivative of
electrical field on $\partial \Omega$ by $\frac{\partial \, \cdot}{\partial n}$, where $n$ denotes the unit
outer normal vector on $\partial \Omega$.

It is well known that for stable implementation of the finite element
solution of Maxwell's equation divergence-free edge elements are the
most satisfactory from a theoretical point of view ~\cite{Nedelec,
  Monk}.  However, the edge elements are less attractive for solution
of time-dependent problems since a linear system of equations should
be solved at every time iteration.  In contrast, P1 elements can be
efficiently used in a fully explicit finite element scheme with lumped
mass matrix \cite{delta, joly}.  It is also well known that numerical
solution of Maxwell equations using nodal finite elements can be
resulted in unstable spurious solutions \cite{MP, PL}.  There are a
number of techniques which are available to remove them, see, for
example, \cite{Jiang1, Jiang2, Jin, div_cor, PL}. 

In the domain decomposition method of this work
we use stabilized
P1 FE method for the numerical solution of (\ref{E_gauge}) in
$\Omega_{\rm FEM}$.
Efficiency of usage an explicit P1 finite element scheme is evident
for solution of CIPs.  In many algorithms which solve electromagnetic
CIPs a qualitative collection of experimental measurements is
necessary on the boundary of the computational domain to determine the
dielectric permittivity function inside it.  In this case the
numerical solution of time-dependent Maxwell's equations are required
in the entire space $\mathbb{R}^{3}$, see for example \cite{BK, 
  BTKM2, BondestaB, TBKF1, TBKF2}, and it is efficient to consider
Maxwell's equations with constant dielectric permittivity function in
a neighborhood of the boundary of the computational domain.
An
explicit P1 finite element scheme with $\sigma=0$ in \eqref{E_gauge}
is numerically tested for solution of time-dependent Maxwell's system
in 2D and 3D in \cite{BMaxwell}.
Convergence analysis of this scheme is presented in \cite{BR1}
and CFL condition is derived in \cite{BR2}.
The scheme of \cite{BMaxwell} is used
for solution of different CIPs for determination of dielectric
permittivity function in non-conductive media in time-dependent
Maxwell's equations using simulated and experimentally generated data,
see \cite{BTKM2, BondestaB, TBKF1,
  TBKF2}.


 The stabilized model problem
 considered in this paper  is:
\begin{equation}\label{eq1domdec}
  \begin{array}{ll}
 \frac{1}{c^2}   \varepsilon_r \frac{\partial^2 E}{\partial t^2} +
 \nabla ( \nabla \cdot E) - \triangle E  -
 \varepsilon_0 \nabla  (\nabla \cdot ( \varepsilon_r  E))
  = - \mu_0\sigma \frac{\partial E}{\partial t} & \mbox{ in } \Omega \times (0, T), \\
    E(\cdot,0) = f_0, \mbox{ and } \frac{\partial E}{\partial t} (\cdot,0) = f_1& \mbox{ in } \Omega, \\
    \frac{\partial E}{\partial n} = - \frac{\partial E}{\partial t} & \mbox{ on } \partial \Omega \times (0,T),
  \end{array}
\end{equation}
with functions $\varepsilon_r, \sigma$
 satisfying conditions \eqref{2.3}.

\section{The domain decomposition algorithm}
\label{sec:hyb}

\begin{figure}[tbp]
\begin{center}
\begin{tabular}{ccc}
 {\includegraphics[scale=0.25, trim={0cm 0cm 0cm 0cm}, clip=]{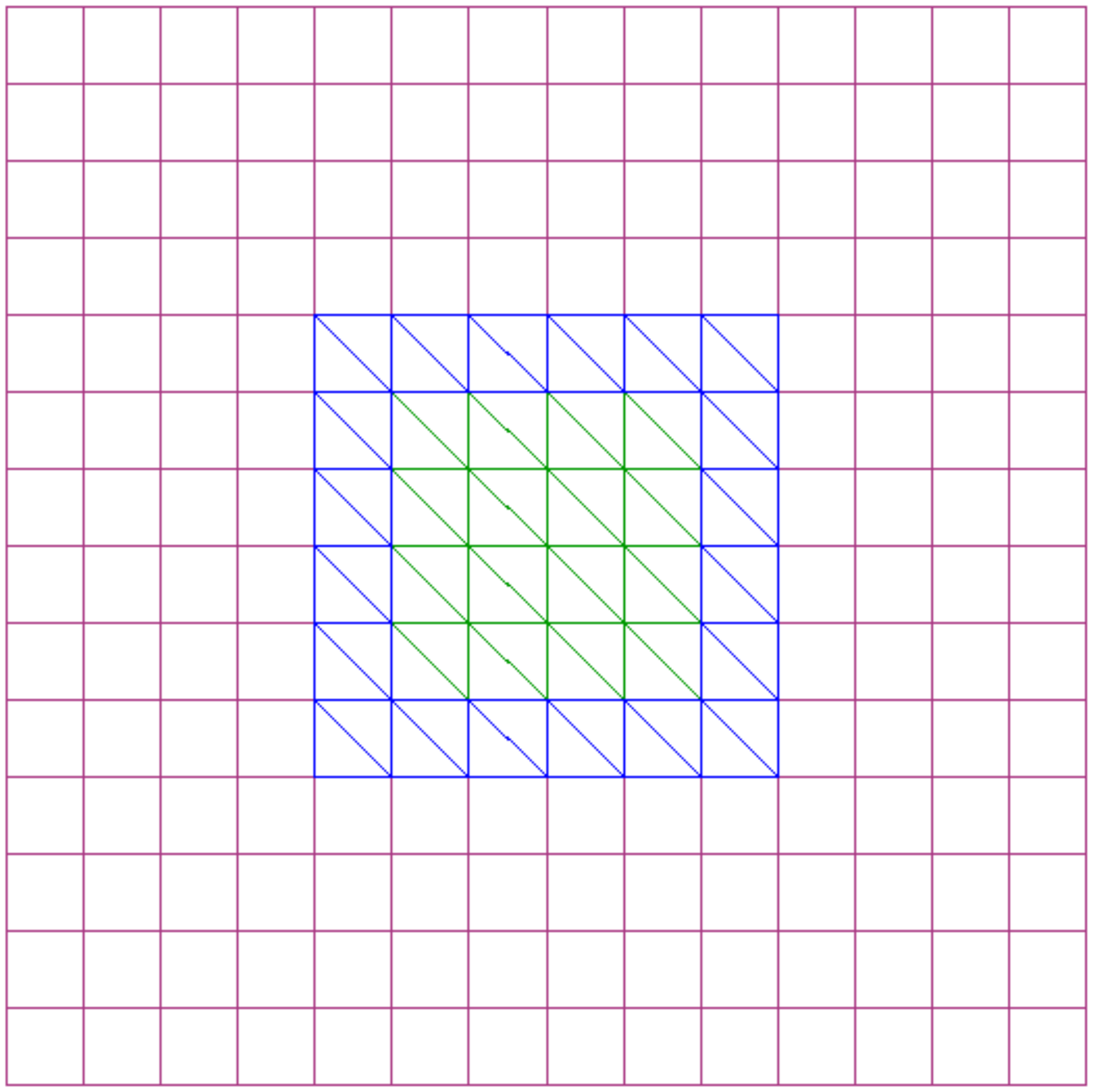}}
  &
 \put(-25,25){\includegraphics[scale=0.25, trim={6cm 5cm 6cm 2cm}, clip=]{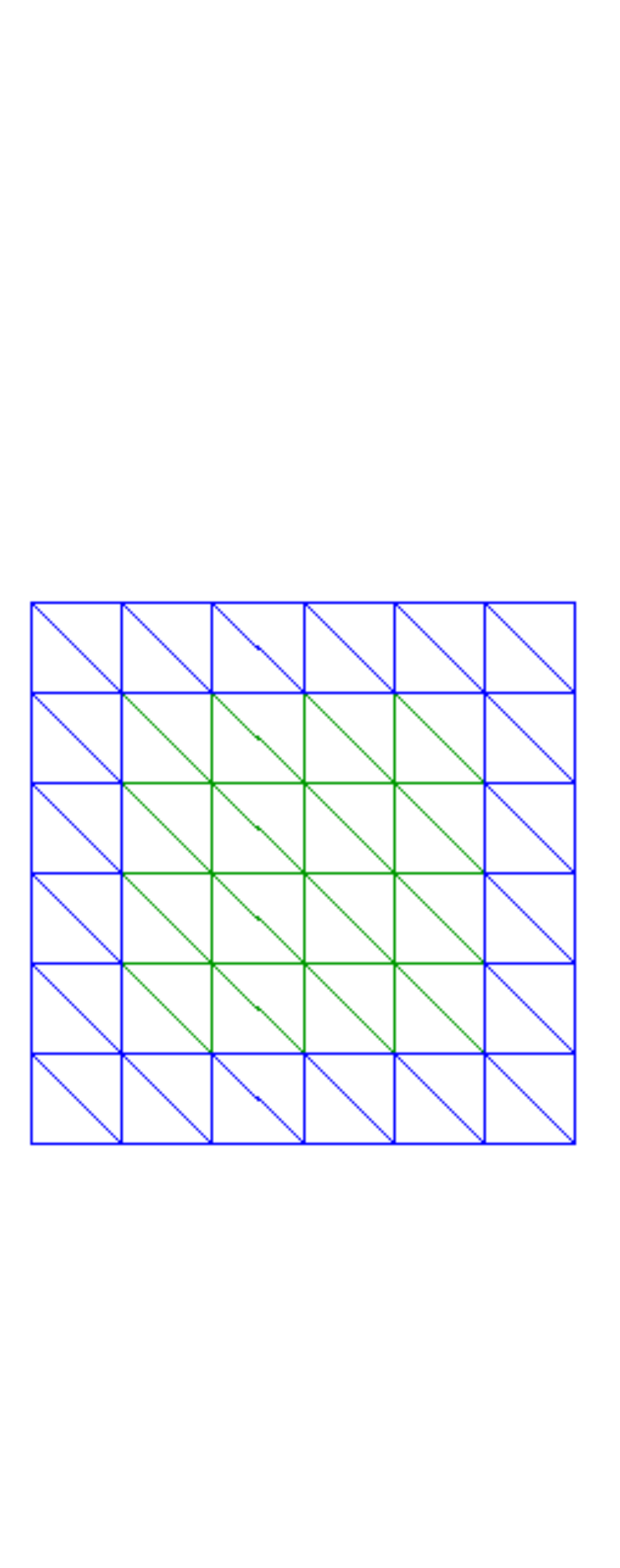}}
\put(0,90){$\Omega_{\rm IN}$}
 \put(0,70){$\Omega_{\rm OUT}$}
 &
    {\includegraphics[scale=0.25,clip=]{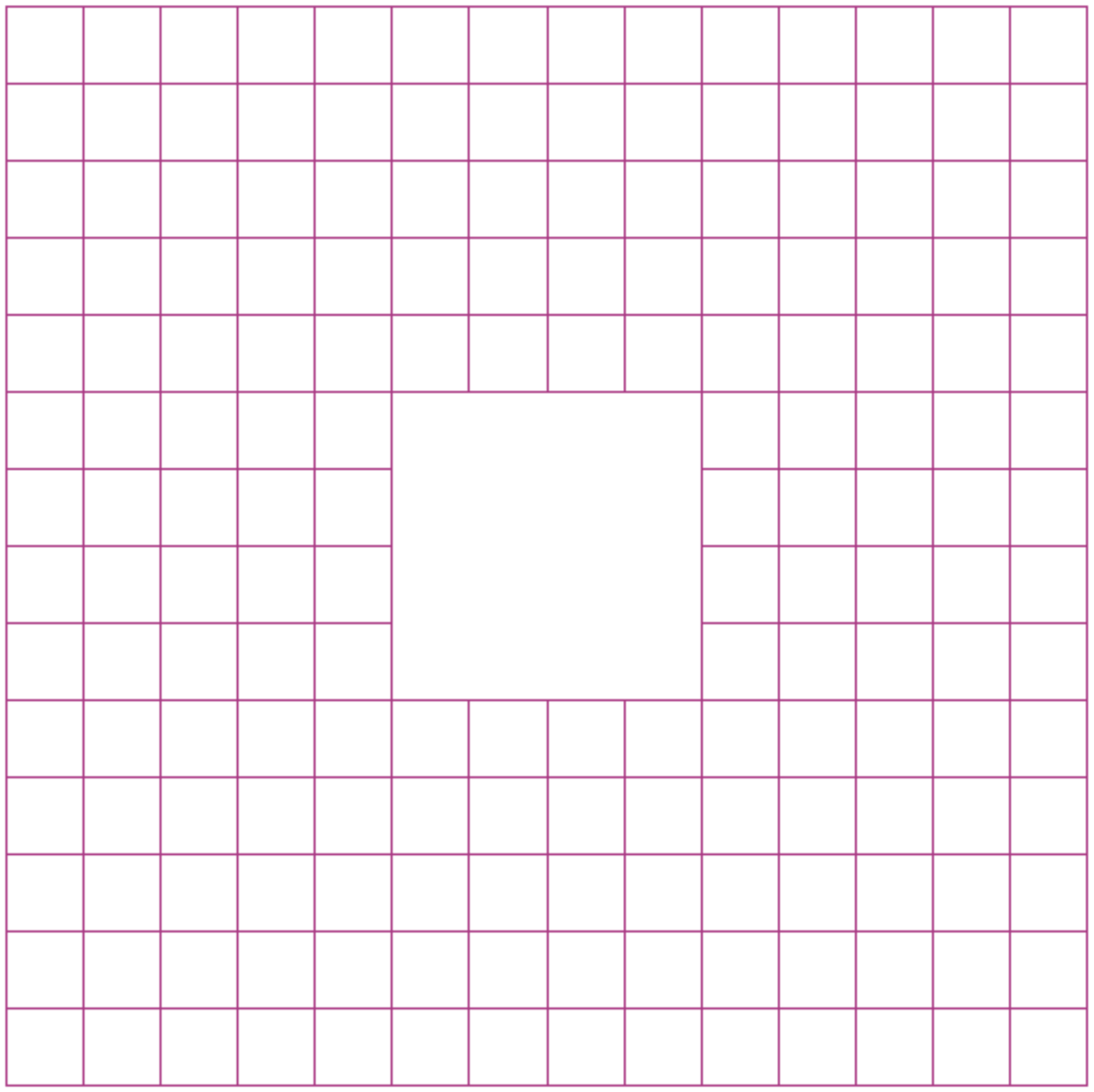}} \\  
 a) $ \Omega = \Omega_{\rm FEM} \cup \Omega_{\rm FDM}$ &
 b) $\Omega_{\rm FEM} = \Omega_{\rm IN} \cup \Omega_{\rm OUT}$ & c) $ \Omega_{\rm FDM}$  
\end{tabular}
\end{center}
\caption{
  \small \emph{\   Domain
    decomposition and mesh discretization in $\Omega$.
    The domain $\Omega$ presented on a) is a combination of
    the quadrilateral finite difference mesh $\Omega_{\rm FDM }$ presented
    on c), and the finite element mesh $\Omega_{\rm FEM}$ presented on
    b).
}}
\label{fig:F2}
\end{figure}


\begin{figure}[tbp]
\begin{center}
\begin{tabular}{cc}
 {\includegraphics[scale=0.24,clip=]{hybriddomainmesh.pdf}}
 &
 {\includegraphics[scale=0.4,  trim={2cm 5cm 2cm 2cm}, clip=]{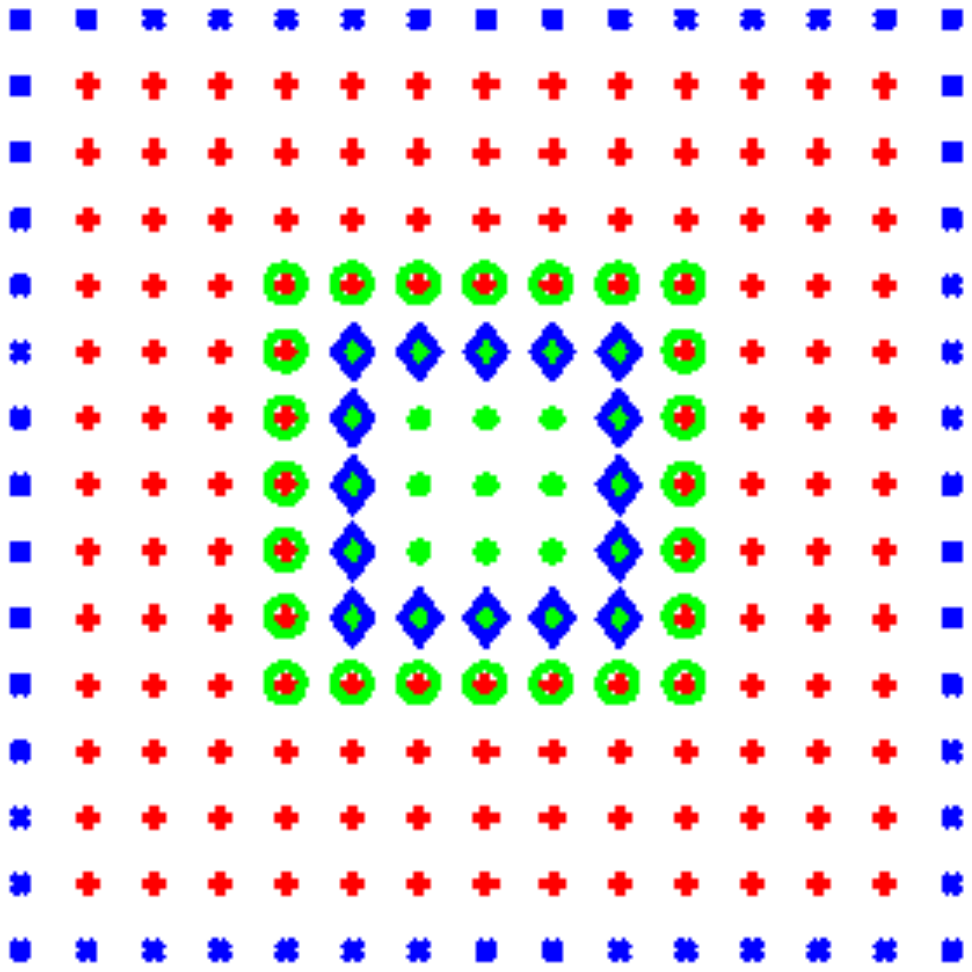}}\\
 a) $ \Omega = \Omega_{\rm FEM} \cup \Omega_{\rm FDM}$  & b) $ \Omega = \Omega_{\rm FEM} \cup \Omega_{\rm FDM}$ \\
\end{tabular}
\end{center}
\caption{\small \emph{\ Coupling between $\Omega_{\rm FEM}$  and $\Omega_{\rm FDM}$. 
  The nodes of the  FE/FD mesh of a) are presented also on b)  as   sets of
  following nodes: $\omega_{\rm o}$ (green circles), $\omega_{\diamond}$ (blue diamonds),
    $\omega_{*}$ (green stars), $\omega_+$ (red pluses), $\omega_{\rm x}$ (blue crosses). These sets  are described in the domain decomposition algorithm. 
}}
\label{fig:F3}
\end{figure}

We now describe the domain decomposition method between two domains
$\Omega_{\rm FEM}$ and $\Omega_{\rm FDM}$ where  FEM is used for computation of the
solution in $\Omega_{\rm FEM}$, and FDM is used in $\Omega_{\rm FDM}$, see Figures
\ref{fig:F2}, \ref{fig:F3}.  Overlapping nodes between $\Omega_{\rm FDM}$ and $\Omega_{\rm FEM}$ are outlined  in Figure  \ref{fig:F3} by green circles (boundary nodes of $\Omega_{\rm FEM}$)  and blue diamonds (inner boundary nodes of $\Omega_{\rm FDM}$).

 The communication between two domains $\Omega_{\rm FEM}$ and
 $\Omega_{\rm FDM}$ is achieved by overlapping of both meshes across a
 two-element thick layer around $\Omega_{\rm FEM}$ - see Figure
 \ref{fig:F3}.  The nodes of the computational domain
 $\Omega$ belong to either of the following sets (see Figure  \ref{fig:F3}-b)):
\begin{itemize}
\item[$\omega_{\rm o}$:] Nodes '$\rm o$' - lie on the boundary $\partial \Omega_{\rm FEM}$ of $\Omega_{\rm FEM}$ and are interior to $\Omega_{\rm FDM}$,
\item[$\omega_{\diamond}$:] Nodes '$\diamond$' -  lie on the inner boundary $\partial \Omega_{\rm FDM}$ of
$\Omega_{\rm FDM}$ and are 
 interior to $\Omega_{\rm FEM}$,
\item[$\omega_{*}$:] Nodes '$*$' are interior to $\Omega_{\rm FEM}$,
\item[$\omega_{+}$:] Nodes '$+$' are interior to $\Omega_{\rm FDM}$, 
  \item[$\omega_{\rm x}$:] Nodes '$\rm x$' lie on the outer boundary $\partial \Omega$ of $\Omega_{\rm FDM}$.
\end{itemize}

Then the main loop in time for the explicit schemes which solves the  problem (\ref{eq1domdec})  with appropriate boundary conditions is shown in Algorithm \ref{alg:alg1}.


\begin{algorithm}[hbt!]
  \centering
    \caption{The domain decomposition algorithm  \label{alg:alg1}}
  \begin{algorithmic}[1]

 \STATE On the structured part of the mesh $\Omega_{\rm FDM}$, where FDM is
  used, update the Finite Difference (FD) solution at nodes $\omega_{+}$ and
  $\omega_{\diamond}$.

  \STATE On the unstructured part of the mesh $\Omega_{\rm FEM}$, where FEM is used,  update the
 Finite Element (FE) solution at nodes $\omega_{*}$ and $\omega_{\rm o}$.

  \STATE Copy FE solution obtained at nodes $\omega_{\diamond}$ as a
  boundary condition for the FD solution in $\Omega_{\rm FDM}$.

  \STATE Copy FD solution obtained at nodes $\omega_{\rm o}$ as a
  boundary condition for the FE solution in $\Omega_{\rm FEM}$.

 \end{algorithmic}
\end{algorithm}

By conditions \eqref{2.3}
functions $\varepsilon_r =1$ and $\sigma=0$ at the overlapping nodes between $\Omega_{\rm FEM}$ and
$\Omega_{\rm FDM}$, and thus,  the Maxwell's equations will transform to the system of uncoupled acoustic wave equations \eqref{6.11}
which leads to the fact
that
the
FEM and FDM discretization schemes coincide on the
common structured overlapping layer.
In this way we avoid
instabilities at interfaces in the domain decomposition algorithm.

\section{Reconstruction algorithms}

\label{sec:recalg}

In this section we develop different optimization
 algorithms which allow
 determination of the relative dielectric permittivity
 function using scattered data of the electric field at
 the boundary of the investigated domain.
 In all algorithms we use assumption that the effective conductivity function
  is known in the investigated domain.

In summary, the main algorithms presented in this section are:

\begin{itemize}

\item Algorithm 2: The domain decomposition algorithm
  for  efficient  solution of forward and   adjoint problems used in
   algorithms 3, 4, 5.

\item Algorithm 3: Optimization algorithm for determination of the
  relative dielectric permittivity function under condition that the
  effective conductivity function is known.

\item Algorithm 4, 5: Adaptive optimization algorithms for
  determination of the relative dielectric permittivity function.
  These algorithms use local adaptive mesh refinement based on a new
  error indicators for improved determination of location, material
  and sizes of the inclusions to be identified.

\end{itemize}

Let  the domain decomposition of the computational domain $\Omega$ be as it is described in Section \ref{sec:hyb}, see also Figures \ref{fig:F3}.
We denote by $\Omega_T
:= \Omega \times (0,T), \partial \Omega_T := \partial \Omega \times
(0,T), T > 0$. 
Let the boundary
$\partial  \Omega  = \partial \Omega_{\rm FDM}^{\rm out} \cup \partial \Omega_{\rm FDM}^{\rm in} $ be the outer boundary $\partial \Omega_{\rm FDM}^{\rm out}$ of $\Omega$
 together with the inner boundary $\partial \Omega_{\rm FDM}^{\rm in}$
  of $\Omega_{\rm FDM}$, and  $\partial
\Omega_{\rm FEM}$ be the boundary of $\Omega_{\rm FEM}$.
Let at $S_T := \partial \Omega_{\rm FDM}^{\rm out} \times (0,T)$ 
we have time-dependent backscattering observations.

Our coefficient inverse problem will be the following.

 \textbf{Inverse Problem (IP) } \emph{Assume that the functions
 }$\varepsilon_r(x),~ \sigma(x)$ \emph{\ satisfy conditions
 (\ref{2.3}) for known }$d_1 >1,$ $d_2 > 0$.\emph{\ Let the function
 $\varepsilon_r$ be unknown in the domain }$\Omega \backslash
 (\Omega_{\rm FDM} \cup \Omega_{\rm OUT})$\emph{. Determine the function}
 $\varepsilon_r(x)$\emph{\ for }$x\in \Omega \backslash
 (\Omega_{\rm FDM} \cup \Omega_{\rm OUT}),$ \emph{\ assuming that the
 function $\sigma(x)$ is known in $\Omega$ and the following function }$\tilde E\left(x,t\right) $\emph{\ is  measured at  $S_T$: }
\begin{equation}\label{2.5}
  E \left(x,t\right) = \tilde E \left(x,t\right) , \forall \left( x,t\right)
  \in   S_T.
\end{equation}

The function $\tilde E\left(x,t\right) $ in (\ref{2.5}) represents
the time-dependent measurements of all components of the electric wave field $E(x,t)= (E_1, E_2, E_3)(x,t)$ at $S_T$.

To  solve IP
we minimize the corresponding
Tikhonov functional and use a Lagrangian approach to do that. We present details of derivation of optimization algorithms in the next section.

\subsection{Derivation of optimization algorithms}

For solution of the \textbf{IP} for Maxwell's system
 (\ref{eq1domdec}) it is natural to minimize the
following Tikhonov functional
\begin{equation}
\begin{split}
J(E, \varepsilon_r) = 
\frac{1}{2} \int_{\Omega_T}(E - \tilde{E} )^2 z_{\delta } \delta_{\rm obs}~ ds dt 
+\frac{1}{2} \gamma  \int_{\Omega}( \varepsilon_r -  \varepsilon^{0})^2~~ dx,
\label{functional}
\end{split}
\end{equation}
where
$\tilde{E}$ is the observed electric field in \eqref{2.5} at the observation points  located at $\partial \Omega_{\rm FDM}^{\rm out}$, $\delta_{\rm obs}  = \sum \delta(\partial \Omega_{\rm FDM}^{\rm out})$ is a sum  of delta-functions at
 the observations points  located at $\partial \Omega_{\rm FDM}^{\rm out}$,
$E$ satisfies the
equations (\ref{eq1domdec}) and thus depends on $\varepsilon_r, \sigma$. We denote by 
$\varepsilon^{0}$ the initial
guess for $\varepsilon_r$, and by $\gamma$ the
regularization parameter.  Here, $z_{\delta}$ is a cut-off
function ensuring the compatibility
conditions  for data, see details in \cite{BondestaB}.

Let us
introduce the following spaces of real valued  functions
\begin{equation}\label{spaces}
\begin{split}
H_E^1(\Omega_T) &:= \{  w \in H^1(\Omega_T):  w( \cdot , 0) = 0 \}, \\
U^{1} &= ((H_{E}^{1}(\Omega_T)^3 \times (H_{E}^{1}(\Omega_T))^3
\times C\left( \overline{\Omega}\right),\\
U^{0} &= (L_{2}\left(\Omega_{T}\right))^3
 \times (L_{2}\left(\Omega_{T}\right))^3 \times
L_{2}\left( \Omega \right).
\end{split}
\end{equation}

  To solve the minimization problem
\begin{equation}\label{minproblem}
\min_{\varepsilon_r} J(E,\varepsilon_r)
\end{equation}
 we take into
account conditions (\ref{2.3}) on the function $\varepsilon_r$ and introduce the Lagrangian
\begin{equation}\label{lagrangian}
\begin{split}
L(u) &= J(E, \varepsilon_r) \\
&+ \int_{\Omega_T}  \lambda \Big ( \frac{1}{c^2}  \varepsilon_r \frac{\partial^2 E}{\partial t^2}
 - \nabla \cdot (\nabla E)
-\nabla \nabla \cdot ((\varepsilon_r \varepsilon_0 - 1)  E)
 + \mu_0 \sigma \frac{\partial E}{\partial t}  \Big)~ dxdt,
\end{split}
\end{equation}
where $u=(E,\lambda, \varepsilon_r)$.

To solve the minimization problem \eqref{minproblem}
we find a stationary point of the Lagrangian with respect to $u$
satisfying $\forall \bar{u}= ( \bar{E}, \bar{\lambda},
\bar{\varepsilon_r})$
\begin{equation}
 L'(u; \bar{u}) = 0 ,  \label{scalar_lagr2}
\end{equation}
where $ L^\prime (u;\cdot )$ is the Jacobian of $L$ at $u$.  For
solution of the minimization problem \eqref{scalar_lagr2} we  develop
conjugate gradient method for  reconstruction of 
parameter $\varepsilon_r$.

To obtain optimality conditions  from \eqref{scalar_lagr2},  we integrate by parts in space and time
the Lagrangian  \eqref{lagrangian},
 assuming that $\lambda( x,T) =\frac{\partial \lambda}{\partial t} \left(
x,T\right) =0, \frac{\partial \lambda}{\partial t} = 
\frac{\partial \lambda}{\partial n}$, and  impose such conditions on the function
$\lambda $  that $ L(E,\lambda, \varepsilon_r)
:=L(u) = J(E, \varepsilon_r).$ 
Using
the facts that $\lambda (x ,T) = \frac{\partial \lambda}{\partial t}
(x,T) =0$, $\nabla \cdot (\varepsilon \lambda) = 0$  and  $\sigma=0, \varepsilon_r=1$ on $\partial \Omega$, together with initial and boundary conditions  of (\ref{eq1domdec}), we get following optimality conditions
for all $\bar{u} \in U^1$,
\begin{equation}\label{forward1}
\begin{split}
0 &= \frac{\partial L}{\partial \lambda}(u)(\bar{\lambda}) =
- \int_{\Omega_T}\frac{1}{c^2}  \varepsilon_r  \frac{\partial \bar{\lambda}}{\partial t} \frac{\partial E}{\partial t}~ dxdt 
+  \int_{\Omega_T}  (  \nabla E) (\nabla  \bar{\lambda}) ~ dxdt\\
&+ \varepsilon_0 \int_{\Omega_T}  (  \nabla \cdot ( \varepsilon_r E)) (\nabla \cdot \bar{\lambda}) ~ dxdt - \int_{\Omega_T}  (  \nabla \cdot E) (\nabla \cdot \bar{\lambda}) ~ dxdt \\
&+ \int_{\Omega_T}   \mu_0 \sigma \frac{\partial E}{\partial t}  \bar{\lambda} ~ dxdt - \int_{\Omega} \frac{\varepsilon_r}{c^2} \bar{\lambda}(x,0) f_1(x) ~dx  \\
& + \int_{\partial \Omega_T} \bar{\lambda} \frac{\partial E}{\partial t} ~d\sigma dt  
 ,~~\forall \bar{\lambda} \in H_{\lambda}^1(\Omega_T);
\end{split}
\end{equation}
\begin{equation} \label{control1}
\begin{split}
0 &= \frac{\partial L}{\partial E}(u)(\bar{E}) =
\int_{\Omega_T}(E- \widetilde{E})~ \bar{E}~ z_{\delta} \delta_{\rm obs}~ d \sigma dt- \int_{\Omega} 
 \frac{\varepsilon_r}{c^2} \frac{\partial{\lambda}}{\partial t}(x,0) \bar{E}(x,0) ~dx \\
&- \int_{\partial \Omega_T}  
\frac{\partial{\lambda}}{\partial t} \bar{E} ~d\sigma dt 
-  \int_{\Omega_T} \frac{\varepsilon_r}{c^2} \frac{\partial \lambda}{\partial t} \frac{\partial \bar{E}}{\partial t}~ dxdt
 + \int_{\Omega_T} ( \nabla  \lambda) (\nabla  \bar{E})  ~ dxdt \\
&+ \varepsilon_0 \int_{\Omega_T}  (  \nabla \cdot ( \varepsilon_r \bar{E})) (\nabla \cdot \lambda) ~ dxdt -
 \int_{\Omega_T}  (  \nabla \cdot \bar{E}) (\nabla \cdot \lambda) ~ dxdt\\
&- \int_{\Omega_T}    \mu_0 \sigma  \bar{E} \frac{\partial \lambda}{\partial t} ~ dxdt, ~\forall \bar{E} \in H_{E}^1(\Omega_T).
\end{split}
\end{equation}
Finally, we obtain the main equation   for iterative update 
$\varepsilon_r$  in the conjugate gradient algorithm
which express that the gradient
with respect to  $\varepsilon_r$  vanishes:
\begin{equation} \label{grad1new} 
\begin{split}
0 &= \frac{\partial L}{\partial  \varepsilon_r}(u)(\bar{\varepsilon}_r)
 =   - \int_{\Omega}  \frac{ \bar{\varepsilon}_r}{c^2} \lambda(x,0)f_1(x)~dx 
 - \int_{\Omega_T}  \frac{ \bar{\varepsilon}_r}{c^2} 
 \frac{\partial \lambda}{\partial t} \frac{\partial E}{\partial t} ~dxdt \\
&+ \varepsilon_0 \int_{\Omega_T} (\nabla \cdot \lambda)( \nabla \cdot ( \bar{\varepsilon}_r E)) ~dxdt
+\gamma \int_{\Omega} (\varepsilon_r - \varepsilon^0) \bar{\varepsilon}_r ~dx,~ x \in \Omega.
\end{split}
\end{equation}

The equation (\ref{forward1}) is the weak formulation of the forward problem
(\ref{eq1domdec}) and the equation (\ref{control1}) is the weak
formulation of the following adjoint problem
\begin{equation}\label{adjoint1}
     \begin{array}{ll}
 \frac{1}{c^2}  \varepsilon_r \frac{\partial^2 \lambda}{\partial t^2} - 
 \triangle \lambda
- \varepsilon_0 \varepsilon_r   \nabla  ( \nabla \cdot  \lambda) +
 \nabla  ( \nabla \cdot  \lambda) - 
\mu_0 \sigma  \frac{\partial \lambda}{\partial t}
= -  (E - \widetilde{E}) z_{\delta} \delta_{\rm obs}
 & \mbox{ in } \Omega_T, \\
\lambda(\cdot, T) =  \frac{\partial \lambda}{\partial t}(\cdot, T) = 0
& \mbox{ in } \Omega, \\
\frac{\partial \lambda}{\partial n} = \frac{\partial \lambda}{\partial t}
& \mbox{ on } S_T.
 \end{array}
\end{equation}

\subsection{The domain decomposition FE/FD method for solution of forward and adjoint problems}

\label{sec:domdec}

\subsubsection{Finite element discretization}
\label{sec:fem}

We denote by $\Omega_{\rm FEM_T} := \Omega_{\rm FEM} \times (0,T), \partial \Omega_{\rm FEM_T} := \partial \Omega_{\rm FEM} \times
(0,T), T > 0$ where
 $\partial \Omega_{\rm FEM}$ is the boundary of  $\Omega_{\rm FEM}$,
 and discretize $\Omega_{{\rm FEM}_T}$ denoting by $K_h = \{K\}$ a partition of
 the domain $\Omega_{\rm FEM}$ into elements $K$  such that
\begin{equation*}
 K_h = \cup_{K \in K_h} K=K_1 \cup K_2...\cup K_l,
\end{equation*}
 where $l$ is the total
 number of elements $K$ in $\overline{\Omega}_{\rm FEM}$.

 Here, $h=h(x)$ is a piecewise-constant mesh function
 defined as 
\begin{equation}\label{meshfunction}
h |_K = h_K ~~~ \forall K \in K_h,
\end{equation}
 representing the local diameter of the elements.
We also denote by  $\partial K_h = \{\partial K\}$ a partition of
 the boundary $\partial \Omega_{\rm FEM}$ into boundaries  $\partial K$ of the elements $K$ such that vertices of these elements belong to  $\partial \Omega_{\rm FEM}$.
We let $J_{\tau}$ be a partition of the time interval $(0,T)$ into time
 intervals $J=(t_{k-1},t_k]$ of uniform length
 $\tau=T/N$
for a given number of time steps $N$.
 We
 assume also a minimal angle condition on the $K_h$ \cite{Brenner, KN}.

To formulate the finite element method in $\Omega$  for
\eqref{scalar_lagr2}
 we
 define the finite element spaces $C_h$, $W_h^E$. 
First, we introduce the finite element trial space $W_h^E$ for every component of the electric field $E$ defined by
\begin{equation}
W_h^E := \{ w \in H_E^1: w|_{K} \in  P_1(K) ,  \forall K \in K_h \}, \nonumber
\end{equation}
where $P_1(K)$  denote the set of piecewise-linear functions on $K$.

To approximate function $\varepsilon_r$ 
   we define the space of piecewise constant functions $C_{h} \subset L_2(\Omega)$, 
\begin{equation}\label{p0}
C_{h}:=\{u\in L_{2}(\Omega ):u|_{K}\in P_{0}(K),\forall K\in  K_h\}, 
\end{equation}
where $P_{0}(K)$ is the piecewise constant function on $K$.
Setting ${\textbf W_h^E(\Omega)} := [W_h^E(\Omega)]^3$
we define $U_h =
{\textbf W_h^E(\Omega)}  \times {\textbf W_h^E(\Omega)} \times C_h$.
 The finite element method  for \eqref{scalar_lagr2} now reads: find $u_h: U_h \times [0,T] \to \mathbb{R}$, such
 that
\begin{equation}
L'(u_h)(\bar{u})=0, ~\forall
\bar{u} \in U_h .  \label{varlagr}
\end{equation}
The equation (\ref{varlagr}) expresses discretized versions of
optimality conditions given by \eqref{forward1}-\eqref{grad1new}.  To
get function $\varepsilon_r$ via optimality condition
\eqref{grad1new} we need solutions first of the
forward problem \eqref{eq1domdec}, and then of the adjoint  problem
 \eqref{adjoint1}.
To solve these problems via the   domain decomposition method, we decompose the computational domain $\Omega = \Omega_{\rm FEM} \cup \Omega_{\rm FDM}$  as it is described in section \ref{sec:hyb}.
Thus, in $\Omega_{\rm FEM}$ we have to solve the following forward problem:
\begin{equation}\label{forwfem1}
  \begin{array}{ll}
\frac{1}{c^2}   \varepsilon_r \frac{\partial^2 E}{\partial t^2} +
 \nabla ( \nabla \cdot E) - \triangle E -
 \varepsilon_0 \nabla  (\nabla \cdot ( \varepsilon_r  E))
  = - \mu_0\sigma \frac{\partial E}{\partial t}
 & \mbox{ in } \Omega_{\rm FEM} \times (0, T), \\
    E(\cdot,0) = f_0, \mbox{ and } \frac{\partial E}{\partial t} (\cdot,0) = f_1 & \mbox{ in } \Omega_{\rm FEM}, \\
    \frac{\partial E}{\partial n} = g & \mbox{ on } \partial \Omega_{\rm FEM} \times (0,T).
  \end{array}
\end{equation}
Here, $g$ is the solution obtained by the finite difference method in $\Omega_{\rm FDM}$ which is saved at $\partial \Omega_{\rm FEM}$.

The equation  (\ref{varlagr}) 
expresses that the finite element method  in $\Omega_{\rm FEM}$ for the solution of the forward problem \eqref{forwfem1} will be:
\emph{Find }$E_h: {\textbf  W_h^E(\Omega_{\rm FEM})}   \times [0,T] \to \mathbb{R}$, 
\emph{ such  that $\forall
 \bar{\lambda}  \in {\textbf  W_h^E(\Omega_{\rm FEM})}$ }
\begin{equation}\label{forwfem2}
  \begin{array}{l}
   \frac{1}{c^2}  \left ( {\varepsilon_r}_h \frac{\partial^2 E_h}{\partial t^2},  \bar{\lambda}  \right ) + (\nabla E_h,\nabla  \bar{\lambda} )+ \varepsilon_0 (\nabla \cdot ({\varepsilon_r}_h E_h), \nabla \cdot  \bar{\lambda}  ) 
    -(\nabla \cdot E_h, \nabla \cdot  \bar{\lambda} ) \\
+ (g_h,  \bar{\lambda} )_{\partial \Omega_{\rm FEM}} + \mu_0 (\sigma_h \frac{\partial E_{h}}{\partial t},  \bar{\lambda}  ) = 0,  \\
   E_h(\cdot,0) = {f_0}_h \mbox{ and } \frac{\partial E_h}{\partial t}(\cdot,0) = {f_1}_h  \mbox{ in } \Omega_{\rm FEM}.  
  \end{array}
\end{equation}
Here, we define ${f_0}_h, {f_1}_h, g_h, {\varepsilon_r}_h, \sigma_h$
to be the usual ${\textbf  W_h^E}$-interpolate of 
$f_0,f_1, g, \varepsilon_r, \sigma$  in \eqref{eq1domdec} in $\Omega_{\rm FEM}$.

To get the  discrete scheme for
 \eqref{forwfem2}   we approximate $E_h(k\tau)$
by $E_h^k$ for $k=1,2,...,N$ using the following  scheme for $k=1,2,\ldots,N-1$
 and  $\forall \bar{\lambda} \in
 {\textbf W_h^E(\Omega_{\rm FEM})}$

\begin{equation}\label{forwfem3}
  \begin{array}{l}
    \frac{1}{c^2}  \left ( {\varepsilon_r}_h
    \frac{E_{h}^{k+1} - 2 E_h^k + E_h^{k-1}}{\tau^2},  \bar{\lambda}  \right ) + (\nabla E_h^k,\nabla  \bar{\lambda} )+
    \varepsilon_0 (\nabla \cdot ({\varepsilon_r}_h E_h^k), \nabla \cdot  \bar{\lambda}  ) 
    -(\nabla \cdot E_h^k, \nabla \cdot  \bar{\lambda} ) \\
    + (g_h^k,  \bar{\lambda} )_{\partial \Omega_{\rm FEM}}
    + \mu_0 (\sigma_h \frac{E_{h}^{k+1} - E_h^{k-1}}{ 2\tau},  \bar{\lambda}  )
    = 0,  \\
   {E_h}^0 = {f_0}_h \mbox{ and } {E_h}^1 =  {E_h}^0 + \tau {f_1}_h \mbox{ in } \Omega_{\rm FEM}.  
  \end{array}
\end{equation}
Rearranging terms in \eqref{forwfem3} we get for $k=1,2,\ldots,N-1$
 and  $\forall \bar{\lambda} \in
 {\textbf W_h^E(\Omega_{\rm FEM})}$
\begin{equation}\label{forwfem5}
  \begin{split}
 &   \left ((1 + \tau c^2 \mu_0\frac{\sigma_h}{2 {\varepsilon_r}_h}) E_h^{k+1}, \bar{\lambda} \right)
  =   \left(2 E_h^k, \bar{\lambda} \right) - \left(E_h^{k-1}, \bar{\lambda} \right)   - \tau^2c^2 (1/ {\varepsilon_r}_h \nabla E_h^k, \nabla \bar{\lambda})\\
 & - \tau^2 c^2 \varepsilon_0 (1/ {\varepsilon_r}_h \nabla \cdot ({\varepsilon_r}_h E_h^k), \nabla \cdot \bar{\lambda}) 
+   \tau^2 c^2 
 (1/ {\varepsilon_r}_h \nabla \cdot E_h^k, \nabla \cdot \bar{\lambda}) \\
&+\tau^2 c^2 \left (\frac{g_h^k}{{\varepsilon_r}_h}, \bar{\lambda} \right)_{\partial \Omega_{\rm FEM}} 
+ \tau c^2 \mu_0(\frac{\sigma_h}{2 {\varepsilon_r}_h}  E_h^{k-1}, \bar{\lambda}), \\
& E_h^0 = {f_0}_h \mbox{ and } E_h^1 =
          E_h^0 + \tau {f_1}_h \mbox{ in } \Omega_{\rm FEM}.
  \end{split}
\end{equation}

The adjoint problem in $\Omega_{\rm FEM}$  will be the following:
\begin{equation}\label{adjoint2}
\begin{split} 
 &\frac{1}{c^2}  \varepsilon_r \frac{\partial^2 \lambda}{\partial t^2} - 
  \triangle \lambda
- \varepsilon_0 \varepsilon_r   \nabla  ( \nabla \cdot  \lambda) +
 \nabla  ( \nabla \cdot  \lambda) - 
\mu_0 \sigma  \frac{\partial \lambda}{\partial t}
= -  (E - \widetilde{E}) z_{\delta} \delta_{\rm obs} \mbox{ in } \Omega_{\rm FEM} \times (0, T),   \\
&\lambda(\cdot, T) =  \frac{\partial \lambda}{\partial t}(\cdot, T) = 0 \mbox{~
~~ for~~~~  } x \in \Omega_{\rm FEM}, \\
&\frac{\partial \lambda}{\partial n} = p  ~\mbox{~~~on~~~}~ \partial \Omega_{\rm FEM_T}.
\end{split}
\end{equation}

The finite element method for the solution of adjoint problem
(\ref{adjoint2}) in $\Omega_{\rm FEM}$ reads: 
\emph{Find }$\lambda_{h}\in {\textbf  W_h^E(\Omega_{\rm FEM})}$ \emph{ such  that $\forall
 \bar{E}  \in {\textbf  W_h^E(\Omega_{\rm FEM})}$ }
\begin{equation}\label{adjfem1}
  \begin{array}{l}
   \frac{1}{c^2}  \left ( {\varepsilon_r}_h \frac{\partial^2 \lambda_h}{\partial t^2},  \bar{E}  \right ) + (\nabla \lambda_h,\nabla  \bar{E} )+ 
\varepsilon_0 (\nabla \cdot \lambda_h, \nabla \cdot ({\varepsilon_r}_h \bar{E} ) ) 
    - (\nabla \cdot \lambda_h, \nabla \cdot  \bar{E} ) \\
- (p_h,  \bar{E} )_{\partial \Omega_{\rm FEM}}- \mu_0 (\sigma_h \frac{\lambda_{h}}{t},  \bar{\lambda}  ) =  -  ((E_h - \widetilde{E}_h) z_{\delta} \delta_{\rm obs}, \bar{E} ).
  \end{array}
\end{equation}
Here, we define $E_h, \widetilde{E}_h, p_h$  to be the usual ${\textbf  W_h^E}$-interpolate of 
$E,\widetilde{E}, p$  in \eqref{adjoint2} in $\Omega_{\rm FEM}$.

We note that the adjoint problem should be solved backwards in time, from
 time $t=T$ to $t=0$.
To get the  discrete scheme for
 \eqref{adjfem1} we approximate $\lambda_h(k\tau)$
by $\lambda_h^k$ for $k=N,N-1,...,1$ using the following  scheme for $k= N-1,
 \ldots,1$:
\begin{equation}\label{adjfem2}
  \begin{array}{l}
  \frac{1}{c^2}    \left ({\varepsilon_r}_h \frac{\lambda_h^{k+1} - 2 \lambda_h^k
 + \lambda_h^{k-1}}{\tau^2}, \bar{E} \right) + (\nabla \lambda_h^k, \nabla \bar{E}) 
     + \varepsilon_0 (\nabla \cdot  \lambda_h^k, \nabla \cdot
({\varepsilon_r}_h  \bar{E})) - (\nabla \cdot \lambda_h^k, \nabla \cdot \bar{E}) \\
		 -
                 \left (p_h^k, \bar{E} \right)_{\partial \Omega_{\rm FEM}} - \mu_0 (\sigma_h \frac{\lambda_h^{k+1} - \lambda_h^{k-1}}{2\tau}, \bar{E})  = -  ((E_h^k - \widetilde{E}_h^k) z_{\delta} \delta_{\rm obs}, \bar{E} )\; \forall \bar{\lambda} \in {\textbf  W_h^E(\Omega_{\rm FEM})}.
  \end{array}
\end{equation}

Multiplying both sides of \eqref{adjfem2} by $\tau^2  c^2/{\varepsilon_r}_h$ and rearranging the terms we obtain:
\begin{equation}\label{adjfem4}
  \begin{split}
 &   \left ((1 + \tau c^2 \mu_0\frac{\sigma_h}{2 {\varepsilon_r}_h})
 \lambda_h^{k-1}, \bar{E} \right)
  =   \left(2 \lambda_h^k, \bar{E} \right) - \left(\lambda_h^{k+1}, \bar{E} \right)   - \tau^2c^2 (1/ {\varepsilon_r}_h \nabla \lambda_h^k, \nabla \bar{E})\\
 & - \tau^2 c^2 \varepsilon_0(1/ {\varepsilon_r}_h 
\nabla \cdot \lambda_h^k, \nabla \cdot ({\varepsilon_r}_h \bar{E})) 
+   \tau^2 c^2 
 (1/ {\varepsilon_r}_h \nabla \cdot \lambda_h^k, \nabla \cdot \bar{E}) \\
&+\tau^2 c^2 \left (\frac{p_h^k}{{\varepsilon_r}_h}, \bar{E} \right)_{\partial \Omega_{\rm FEM}} 
+ \tau c^2 \mu_0(\frac{\sigma_h}{2 {\varepsilon_r}_h}  \lambda_h^{k+1}, \bar{E})
-\tau^2  c^2  (1/ {\varepsilon_r}_h   (E_h^k - \widetilde{E}_h^k) z_{\delta}
\delta_{\rm obs}, \bar{E} ),
  \end{split}
\end{equation}
for $k= N-1,
 \ldots,1$,
 $ \forall \bar{E} \in
          {\textbf W_h^E(\Omega_{\rm FEM})}$

We note that usually $\dim U_{h}<\infty $ and $U_{h}\subset U^{1}$ as
a set and we consider $U_{h}$ as a discrete analogue of the space
$U^{1}.$ We introduce the same norm in $U_{h}$ as the one in $U^{0}$,
\begin{equation}\label{equiv}
\left\Vert \bullet \right\Vert _{U_{h}}:=\left\Vert \bullet \right\Vert
_{U^{0}},
\end{equation}
 where $U_0$ is defined in (\ref{spaces}).  From (\ref{equiv}) follows
 that  all norms  in finite dimensional spaces  are equivalent. This
 allows us in numerical simulations of  section \ref{sec:numex}
 compute the discrete function ${\varepsilon_r}_h$, which is approximation of $\varepsilon_r(x)$, in the space
 $C_h$.

\subsection{Fully discrete scheme in $\Omega_{\rm FEM}$}
\label{sec:discrete}

In this section we   present schemes for computations of the solutions 
of forward  \eqref{eq1domdec} and adjoint (\ref{adjoint1}) problems   in $\Omega_{\rm FEM}$.
  After expanding functions $E_h(x)$ and $\lambda_h(x)$ in
terms of the standard continuous piecewise linear functions
$\{\varphi_i(x)\}_{i=1}^M$ in space  as 
\begin{equation*}
\begin{split}
E_h(x) = \sum_{i=1}^M E_{h_{i}}\varphi_i(x),~ 
\lambda_h(x) &= \sum_{i=1}^M \lambda_{h_{i}} \varphi_i(x), 
\end{split}
\end{equation*}
 where $E_{h_{i}}$ and
$\lambda_{h_{i}}$ denote unknown coefficients at the mesh point $x_i \in
K_h, i =1,...,M$, substitute them into
(\ref{forwfem5}) and (\ref{adjfem4}), correspondingly, with
$\bar{\lambda}(x,t) =\bar{E}(x,t) = \sum_{j=1}^M \varphi_j(x)$,
and obtain the system of linear equations  for computation of the forward problem \eqref{eq1domdec}:
\begin{equation} \label{femod1new}
\begin{split}
 M_1 E^{k+1} &=   
2 M E^{k}  - M E^{k-1}  - \tau^2 c^2  G_1 E^{k}  -  \tau^2 c^2 \varepsilon_0
 G_2 E^{k}  \\
&+  
\tau^2 c^2 G_3 E^{k} + \tau^2 c^2 F^k + \tau  c^2 \mu_0 M_2 E^{k-1}.
\end{split}
\end{equation}
  Here, $M, M_1, M_2$ are the assembled block mass matrices in space, $G_1,
  G_2, G_3$ are the assembled block matrices in space, $F^k$ is the assembled load vector at the time iteration $k$,
 $E^k$ denote the nodal values
  of $E_h(\cdot,t_k)$, $\tau$ is the time step.
Now we define the mapping $F_K$ for the reference element $\hat{K}$
such that $F_K(\hat{K})=K$ and let $\hat{\varphi}$ be the piecewise
linear local basis function on the reference element $\hat{K}$ such
that $\varphi \circ F_K = \hat{\varphi}$.  Then, the explicit formulas
for the entries in system of equations (\ref{femod1new}) at each element $K$ can be
given as:
\begin{equation}\label{eq11}
\begin{split}
 M_{{i,j}}^{K}  &=    (\varphi_i(x) \circ F_K, \varphi_j(x)\circ F_K)_K, \\
  {M_1}_{{i,j}}^{K}  &=    ((1 + \tau c^2 \mu_0\frac{\sigma_h}{2 \varepsilon_{r_h}}) 
\varphi_i(x) \circ F_K, \varphi_j(x)\circ F_K)_K, \\
 {M_2}_{{i,j}}^{K}  &=  ( \frac{\sigma_h}{2 {\varepsilon_r}_h} \varphi_i(x) \circ F_K,
  \varphi_j(x)\circ F_K)_K, \\
 {G_1}_{{i,j}}^{ K} & =(\frac{1}{{\varepsilon_r}_h} \nabla  \varphi_i \circ F_K, \nabla  \varphi_j \circ F_K)_K, \\
{G_2}_{{i,j}}^{ K} & =(\frac{1}{{\varepsilon_r}_h} \nabla \cdot ({\varepsilon_r}_h \varphi_i) \circ F_K, \nabla \cdot \varphi_j \circ F_K)_K, \\
{G_3}_{{i,j}}^{ K} & =(\frac{1}{\varepsilon_{r_h}} \nabla \cdot \varphi_i \circ F_K, \nabla \cdot \varphi_j \circ F_K)_K, \\
  F_{j}^k &= (  \frac{g_h^k}{ \varepsilon_{r_h}}, \varphi_j \circ F_K )_{\partial K},
\end{split}
 \end{equation}
where $(\cdot,\cdot)_K$ denotes the $L_2(K)$ scalar product and $\partial K$ is the part of the boundary of element $K$ which lies at $\partial \Omega_{\rm FEM}$.

For the case of adjoint problem (\ref{adjfem4}) we get the 
 system of  linear equations:
\begin{equation} \label{systemadj}
\begin{split}
 M_1 \lambda^{k-1} &=   
2 M \lambda^{k}  - M \lambda^{k+1}  - \tau^2 c^2  G_1 \lambda^{k}  
-  \tau^2 c^2 \varepsilon_0 G_2^{T} \lambda^{k}  \\
&+  \tau^2 c^2 G_3 \lambda^{k} + \tau^2 c^2 P_1^k
+ \tau  c^2 \mu_0 M_2 \lambda^{k+1} - \tau^2 c^2 P_2^k.
\end{split}
\end{equation}
  Here, $M, M_1, M_2, G_1,
  G_2, G_3$ are the assembled block matrices in space with explicit entries
 given in \eqref{eq11}, and $P_1^k, P_2^k$ are assembled load vectors at the time iteration $k$  with explicit entries
\begin{equation}\label{eq12}
\begin{split}
  {P_1}_{j}^k &= (  \frac{p_h^k}{ {\varepsilon_r}_h}, \varphi_j \circ F_K )_{\partial K},\\
  {P_2}_{j}^k &= (1/ {\varepsilon_r}_h   (E_h^k - \widetilde{E}_h^k) z_{\delta}
  \delta_{\rm obs},  \varphi_j \circ F_K )_{K},
\end{split}
 \end{equation}
 $\lambda^k$ denote the nodal values
  of $\lambda_h(\cdot,t_k)$, $\tau$ is the time step.

Finally, for reconstructing  $\varepsilon_r(x)$ in
$\Omega_{IN}$ we can use a gradient-based method with an appropriate
initial guess values $\varepsilon^0$.  The discrete versions
in space of the gradients given in \eqref{grad1new},
after integrating by parts in space of the third term in the right
hand side of \eqref{grad1new}, have the form  $\forall x \in \Omega_{\rm IN}$:
\begin{equation} \label{gradient1} 
\begin{split}
{g_h} &=  -\frac{1}{c^2} \lambda_h(x,0){f_1}_h
 - \frac{1}{c^2} \int_{0}^T 
 \frac{\partial \lambda_h}{\partial t} \frac{\partial E_h}{\partial t} ~dt \\
&+  \varepsilon_0 \int_{0}^T (\nabla \cdot \lambda_h)( \nabla \cdot E_h ) ~dt
+\gamma({\varepsilon_r}_h - \varepsilon_h^0),
\end{split}
\end{equation}
where $\varepsilon_h^0$  is interpolant of $\varepsilon^0$.
We note that because of usage of the domain decomposition method, gradient
\eqref{gradient1}   should be updated only in $\Omega_{IN}$
since in $\Omega_{\rm FDM}$  and in  $\Omega_{\rm OUT}$  by condition \eqref{2.3} we have  $\varepsilon_r=1, \sigma=0$.
In \eqref{gradient1} $E_h$ and $\lambda_h$   are computed values of the forward  and adjoint
  problems  using schemes
(\ref{femod1new}), \eqref{systemadj},  correspondingly, 
 and ${\varepsilon_r}_h$  is  approximate value of the computed
 relative  dielectric permittivity  function $\varepsilon_r$.

  \subsubsection{Finite difference formulation}
\label{sec:fdm}

We recall now that from  conditions (\ref{2.3})  it follows that in
$\Omega_{\rm FDM}$
 the function $\varepsilon_r(x)= 1, \sigma=0$. This means that in
 $\Omega_{\rm FDM}$ the model problem (\ref{eq1domdec})
 transforms to the following forward problem for uncoupled system of acoustic wave equations  for $E=(E_1,E_2,E_3)$:
\begin{align*}\label{FDMmodel2new}
&\frac{\partial^2 E}{\partial t^2}  -  \Delta  E   = 0  &&~~\mbox{in}~~
 \Omega_{\rm FDM} \times (0,T), \\
  &E(\cdot,0) = f_0, ~~~\frac{\partial E}{\partial t}(\cdot,0) = f_1 &&~~ \mbox{in}~~ \Omega_{\rm FDM},  \\
&\frac{\partial E}{\partial n} =-\frac{\partial E}{\partial t} &&~~\mbox{on}~ S_T,\\
&\frac{\partial E}{\partial n} = \frac{\partial E_{\rm FEM}}{\partial n} &&~~\mbox{on}~ \partial \Omega_{\rm FDM}^{\rm in},
\end{align*}
where     $\frac{\partial E_{\rm FEM}}{\partial n}$  are
known  values at $\partial \Omega_{\rm FDM}^{\rm in}$.

Using standard finite difference discretization of the  first
equation in ~(\ref{FDMmodel2new}) in $\Omega_{\rm FDM}$ we obtain the following explicit
scheme for  every component  of the solution $E$ of the forward problem 
\eqref{FDMmodel2new}
\begin{equation}  \label{fdmschemeforward2}  
  E_{l,j,m}^{k+1} =  \tau^2 \Delta E_{l,j,m}^k + 2 E_{l,j,m}^k - E_{l,j,m}^{k-1},
\end{equation}
with correspondingly discretized absorbing boundary conditions.  In equations above,
$E_{l,j,m}^{k}$ is the finite difference solution on the time iteration $k$ at the
discrete point $(l,j,m)$,
 $\tau$ is the time step, and $ \Delta
E_{l,j,m}^k$ is the discrete Laplacian.

The   adjoint problem in $\Omega_{\rm FDM}$ will be:
\begin{equation}
\begin{split} \label{FDMadjoint}
 \frac{\partial^2 \lambda}{\partial t^2} - 
 \Delta   \lambda  = -  (E - \widetilde{E}) z_{\delta} \delta_{\rm obs} & ~~\mbox{in}~~
 \Omega_{\rm FDM} \times (0,T),  \\
 \lambda(\cdot, T) =  \frac{\partial \lambda}{\partial t}(\cdot, T) = 0 &
~~ \mbox{in}~~ \Omega_{\rm FDM},  \\
 \frac{\partial \lambda}{\partial n} = \frac{\partial \lambda}{\partial t} &
 ~~\mbox{on}~~  S_T,\\
\frac{\partial \lambda}{\partial n}  = \frac{\partial \lambda_{\rm FEM}}{\partial n} &  ~~ \mbox{on}~~
\partial \Omega_{\rm FDM}^{\rm in},\\
\end{split}
\end{equation}
where     $\frac{\partial \lambda_{\rm FEM}}{\partial n}$  are
known  values at $\partial \Omega_{\rm FDM}^{\rm in}$.

  Similarly with \eqref{fdmschemeadjoint}  
we get the following explicit scheme for the solution of adjoint problem 
 \eqref{FDMadjoint}  in $\Omega_{FDM}$ which we solve backward in time:
\begin{equation}  \label{fdmschemeadjoint}  
  \lambda_{l,j,m}^{k-1} = -\tau^2 (E - \widetilde{E})_{l,j,m}^k z_{\delta}
  \delta_{\rm obs} 
 + \tau^2 \Delta \lambda_{l,j,m}^k + 2 \lambda_{l,j,m}^k - \lambda_{l,j,m}^{k+1},
\end{equation}
with  corresponding boundary conditions.
  In equations \eqref{fdmschemeforward2}, \eqref{fdmschemeadjoint}  
 $(\cdot)_{l,j,m}^{k}$ is the solution on the
time iteration $k$ at the discrete point $(l,j,m)$, 

 We note that  we use  FDM only
 inside $\Omega_{\rm FDM}$,
 and thus  computed values of  $\frac{\partial E_{\rm FEM}}{\partial n}$ and 
$\frac{\partial \lambda_{\rm FEM}}{\partial n}$ can be approximated and
 will be known at $\partial \Omega_{\rm FDM}^{\rm in}$
 through the finite element solution in $\Omega_{\rm FEM}$, see details in the domain decomposition Algorithm 2.

\subsection{The domain decomposition algorithm to solve forward and adjoint problems}
\label{sec:hybalg}

First we  present domain decomposition algorithm for the
solution of state and adjoint problems.
  We  note that because of
using explicit   finite difference scheme  in $\Omega_{\rm FDM}$  we need to choose
time step $\tau$ accordingly to the CFL stability condition \cite{BR1,BR2,Cohen} such that the whole scheme remains stable. 


\vspace{3mm}

\begin{algorithm}[hbt!]
  \centering
    \caption{ The domain decomposition algorithm  to solve forward and adjoint problems}
    \begin{algorithmic}[1]

\STATE 
Construct the finite element mesh $K_{h}$ in $\Omega_{\rm FEM}$
and the finite difference mesh in $\Omega_{\rm FDM}$ as well as   time partition $J_{\tau}$ of the time interval $\left( 0,T\right)$.
At every time step $k$ we perform the
following operations:

\STATE On  the mesh in $\Omega_{\rm FDM}$ compute
  $E^{k+1}$, $\lambda^{k-1}$ from (\ref{fdmschemeforward2}),
  (\ref{fdmschemeadjoint}), correspondingly, using absorbing boundary conditions  at the outer boundary $\partial \Omega$, with $E^k, E^{k-1}$ and
$\lambda^k, \lambda^{k+1}$ known.

\STATE On the  mesh $K_h$ in  $\Omega_{\rm FEM}$  compute
$E^{k+1}, \lambda^{k-1}$  using the finite element schemes
(\ref{femod1new}),  \eqref{systemadj}, correspondingly,  with $E^k, E^{k-1}$ and
$\lambda^k, \lambda^{k+1}$ known.

\STATE Use the values of the functions $E^{k+1}, \lambda^{k-1}$
  at nodes $\omega_{*}$  overlapping with nodes $\omega_{\diamond}$, which are computed using the finite
  element schemes (\ref{femod1new}),  \eqref{systemadj},  correspondingly,  as a boundary conditions at the inner
  boundary $\partial \Omega_{\rm FDM}^{\rm in}$ for the finite difference method in
  $\Omega_{\rm FDM}$.
  
\STATE Use the values of the functions $E^{k+1}, \lambda^{k-1}$ at
  nodes $\omega_{\rm o}$  overlapping with nodes $\omega_{\rm +}$, which are computed using the finite difference
  schemes (\ref{fdmschemeforward2}), (\ref{fdmschemeadjoint}),
  correspondingly, as a boundary conditions at $\partial \Omega_{\rm FEM}$ for the finite element
  method in $\Omega_{\rm FEM}$.
  
\STATE
 Apply swap of the solutions for the computed functions $E^{k+1},
  \lambda^{k-1}$. Set $k=k+1$ for forward problem and $k = k-1$ for adjoint problem  and go to step 2.

 \end{algorithmic}
\end{algorithm}

\subsection{Reconstruction algorithm for the solution of  inverse problem  IP}
\label{subsec:ad_alg}

We use conjugate gradient method (CGM) for iterative update of
approximation ${\varepsilon_r}_{h}^{m}$ 
of the function ${\varepsilon_r}_{h}$, where $m$ is the
number of iteration in the optimization algorithm. We introduce the following function
\begin{equation} \label{gradient1h} 
\begin{split}
{g}_h^m(x) &=  -\frac{1}{c^2} \lambda_h^m(x,0){f_1}_h(x)
 - \frac{1}{c^2}  \int_{0}^T
 \frac{\partial \lambda_h^m}{\partial t} \frac{\partial E_h^m}{\partial t} ~dt \\
&+ \varepsilon_0  \int_{0}^T (\nabla \cdot \lambda_h^m)( \nabla \cdot E_h^m ) ~dt
+\gamma({\varepsilon_r}_h^m - \varepsilon_h^0),
\end{split}
\end{equation}
where functions $E_{h}^m,\lambda _{h}^m$\ are computed by solving
the state and adjoint problems
with $\varepsilon_r:={\varepsilon_r}_{h}^{m}, \sigma:=\sigma_h^{m}$.\\

\begin{algorithm}[hbt!]
  \centering
  \caption{ Conjugate gradient algorithm for determination of the relative dielectric permittivity  function}
    \begin{algorithmic}[1]

\STATE    Initialize the mesh  in $\Omega$ and time partition $J_{\tau}$ of the time interval $\left( 0,T\right) .$
Start with the initial approximation
  ${\varepsilon_r}_{h}^{0}= \varepsilon_h^{0}$  with known ${\sigma}_{h}$,
 and compute the
sequence of ${\varepsilon_r}_{h}^{m} $ via the following steps:

\STATE Compute solutions $E_{h}( x,t,{\varepsilon_r}_{h}^{m}, \sigma_h)$
  and $\lambda _{h}\left( x,t, {\varepsilon_r}_{h}^{m}, \sigma_h\right) $ of forward and
  adjoint problems on $K_{h}$ and $J_{\tau}$ using the domain
  decomposition algorithm (Algorithm 2).

\STATE  Update the function
 $ {\varepsilon_r}_{h}:= {\varepsilon_r}_{h}^{m+1}$   on $K_{h}$ and $J_{\tau}$ using the CGM  as
\begin{equation*}
\begin{split}
 {\varepsilon_r}_{h}^{m+1} &=  {\varepsilon_r}_{h}^{m}  + \alpha d^m(x),\\
\end{split}
\end{equation*}
where $\alpha$ is the step-size in the gradient update \cite{Peron} and
\begin{equation*}
\begin{split}
 d^m(x)&=  -g_h^m(x)  + \beta^m  d^{m-1}(x),\\
\end{split}
\end{equation*}
with
\begin{equation*}
\begin{split}
 \beta^m &= \frac{\| g_h^m(x)\|^2}{ \| g_h^{m-1}(x)\|^2},\\
\end{split}
\end{equation*}
Here, $d^0(x)= -g_h^0(x)$.

\STATE  Stop computing ${\varepsilon_r}_{h}^m$ at the iteration
  $M:=m$ and obtain the function ${\varepsilon_r}_h^M :={\varepsilon_h}_r^m$ if
  either $\|g_h^{m}\|_{L_{2}( \Omega)}\leq \theta$ or norms
  $\|{\varepsilon_r}_{h}^{m}\|_{L_{2}(\Omega)}$ are stabilized. Here,
  $\theta$ is the tolerance  chosen by the user.
  Otherwise set
  $m:=m+1$ and go to step 2.

 \end{algorithmic}
\end{algorithm}

\subsection{Adaptive algorithms for solution of the inverse problem  IP}

\label{sec:ref}

Adaptive algorithm allows improvement of already computed relative dielectric
permittivity function ${\varepsilon_r}_h^M$ obtained on the initially
non-refined mesh in the previous optimization algorithm (Algorithm 3).
The idea of the local mesh refinement (note that we need it only in
$\Omega_{IN}$) is that it should be refined in all neighborhoods of
all points in the mesh $K_h$ where the function $|h {\varepsilon_r}_h|
$ achieves its maximum value, or where
$|J_{\varepsilon_r}'({\varepsilon_r}_h)|$
achieves its maximal
values.  These local mesh refinements recommendations are based on a
posteriori error estimates for the error $| \varepsilon_r -
{\varepsilon_r}_h|$ in the reconstructed function $\varepsilon_r$ (
see the first mesh refinement indicator), and for the error
$|J({\varepsilon_r}) - J({\varepsilon_r}_h)|$ in the Tikhonov's functional
(see the second mesh refinement indicator), respectively. The proofs of
these a posteriori error estimates for arbitrary Tikhonov's functional is given
in \cite{hyperthermi}.
 A posteriori error for the Tikhonov's functional \eqref{functional} can
be derived  using technique of  \cite{BondestaB}, and it is a topic of
ongoing research. Assuming that we have proof of these a posteriori error indicators, let us show how to   compute them.

We  define by $E(\varepsilon_r, \sigma), \lambda(\varepsilon_r,\sigma)$
 the exact solutions of the
forward and adjoint problems for exact $\varepsilon_r, \sigma$, respectively.
 Then by defining  
\begin{equation*}
u(\varepsilon_r,\sigma) = (E(\varepsilon_r, \sigma),
 \lambda(\varepsilon_r, \sigma), \varepsilon_r) \in U^1,
\end{equation*}
and  using the fact that  for exact solutions  $E(\varepsilon_r, \sigma),
 \lambda(\varepsilon_r, \sigma)$
 we have
\begin{equation}
J(E(\varepsilon_r,\sigma),\varepsilon_r) = L(u(\varepsilon_r, \sigma)).
\end{equation}
~ Assuming  now that solutions $E(\varepsilon_r, \sigma), \lambda(\varepsilon_r, \sigma) $  are sufficiently stable
   we can write that the  Frech\'{e}t derivative of the Tikhonov functional is the following function
\begin{equation}\label{derfunc}
\begin{split}
J'_{\varepsilon_r} (\varepsilon_r, \sigma)  =  
\frac{\partial J}{\partial \varepsilon_r}(E(\varepsilon_r, \sigma), \varepsilon_r) =  \frac{\partial L}{\partial \varepsilon_r}(u(\varepsilon_r, \sigma)).
\end{split}
\end{equation}

Inserting (\ref{grad1new}) into  (\ref{derfunc}),
we get
\begin{equation} \label{derfunc2}
  \begin{split}
   J'_{\varepsilon_r} (\varepsilon_r, \sigma)  &= - \frac{1}{c^2} \lambda(x,0)f_1(x) -
\frac{1}{c^2}  \int_0^T
 \frac{\partial \lambda}{\partial t} \frac{\partial E}{\partial t} ~dt \\
 &- \varepsilon_0 \int_0^T  E \nabla (\nabla \cdot \lambda) ~dt
+ \gamma (\varepsilon_r - \varepsilon^0)(x).
\end{split}
\end{equation}

In the second mesh refinement indicator is used discretized version of
\eqref{derfunc2}  computed for approximations $ ({\varepsilon_r}_h, \sigma_h)$.

\begin{itemize}

\item  \textbf{The First Mesh Refinement Indicator}
 \emph{Refine the mesh in neighborhoods of those points of
 }$K_h$\emph{\ where the function }$|h {\varepsilon_r}_h| $\emph{\ attains its maximal values. In other words, refine
   the mesh in such subdomains of }$K_h $\emph{\ where}
\begin{equation*}
|h {\varepsilon_r}_h| 
\geq \widetilde{\beta} \max \limits_{K_h} |h {\varepsilon_r}_h|.
\end{equation*}
\emph{Here, $ \widetilde{\beta} \in (0,1)$ is a number which should
be chosen computationally and $h$ is the mesh function (\ref{meshfunction}) of the finite element mesh $K_h$}.


  \item \textbf{The Second Mesh Refinement Indicator} \emph{Refine the mesh in
  neighborhoods of those points of } $K_h$ \emph{\ where the
  function } $|J_{\varepsilon_r}'(E,{\varepsilon_r}_h)|$ 
\emph{\ attains its maximal
  values. More precisely, let }$\beta \in (0,1)
$\emph{\ be the tolerance number which should be chosen in
  computational experiments. Refine the mesh $K_h$ in such subdomains  where}
\begin{equation}\label{ind2}
  |J_{\varepsilon_r}'(E, {\varepsilon_r}_h)|
  \geq \beta \max_{K_h} |J_{\varepsilon_r}'(E, {\varepsilon_r}_h)|.
\end{equation}

\end{itemize}

\textbf{Remarks}

\begin{itemize}

\item
1. We note that in (\ref{ind2}) exact values of $E(x,t),
\lambda(x,t)$  are used obtained with the already computed functions
$({\varepsilon_r}_h, \sigma_h)$, see \eqref{derfunc2}. However, in our algorithms and in
computations we approximate exact values of $E(x,t), \lambda(x,t)$ by
the computed ones $E_h(x,t), \lambda_h(x,t)$.

\item 2. In both mesh refinement indicators we used the fact that
  functions $\varepsilon_r, \sigma$ are unknown only in $\Omega_{IN}$.

\end{itemize}

We define the minimizer of the Tikhonov functional (\ref{functional})
and its approximated finite element solution on $k$ times adaptively
refined mesh $K_{h_k}$ by $\varepsilon_r$ and
${\varepsilon_r}_k$, correspondingly.  In our both mesh refinement recommendations
we need compute the functions ${\varepsilon_r}_k$
on the mesh $K_{h_k}$. To do that we apply Algorithm 3 (conjugate gradient
algorithm). We will define by
${\varepsilon_r}_k :={\varepsilon_r}_h^M$ values obtained at steps 3 of
the conjugate gradient algorithm.

\begin{algorithm}[hbt!]
  \centering
  \caption{ Adaptive Algorithm, first version}
    \begin{algorithmic}[1]

      \STATE Construct  the finite difference mesh in $\Omega_{\rm FDM}$.
      Choose an initial space-time mesh ${K_h}_{0} \times
  J_{\tau_0}$ in $\Omega_{\rm FEM} \times [0,T]$.
 Compute the sequence of
  ${\varepsilon_r}_k, k >0$, via following steps:

\STATE  Obtain numerical solution  ${\varepsilon_r}_k$ with known function $\sigma_k$  on $K_{h_k}$ using the Algorithm 3 (Conjugate Gradient Method).

\STATE 
Refine such elements in the mesh $K_{h_k}$  where  the first mesh refinement indicator 
\begin{equation} \label{alg2_2}
  |h {\varepsilon_r}_k| \geq \widetilde{\beta}_k \max_{{K_h}_k}
  |h {\varepsilon_r}_k|
\end{equation}
is satisfied. Here, the tolerance numbers $ \widetilde{\beta}_k  \in \left(
0,1\right) $ are chosen by the user.

\STATE  Define a new refined mesh as $K_{h_{k+1}}$ and
  construct a new time partition $J_{\tau_{k+1}}$ such that the CFL
  condition 
  is satisfied.  Interpolate ${\varepsilon_r}_k, \sigma_k$ on a new mesh
  $K_{h_{k+1}}$ and perform steps 2-4 on the space-time mesh
    $K_{h_{k+1}} \times J_{\tau_{k+1}}$. Stop mesh refinements when
  $||{\varepsilon_r}_k - {\varepsilon_r}_{k-1}|| < tol_1$
  or $|| g_h^k(x)|| < tol_2$, where $tol_i, i=1,2$ are
  tolerances chosen by the user.

 \end{algorithmic}
\end{algorithm}

\vspace{0.5cm}

\begin{algorithm}[hbt!]
  \centering
  \caption{ Adaptive Algorithm, second version}
    \begin{algorithmic}[1]

      \STATE   Choose an initial space-time mesh $K_{h_{0}}\times
  J_{\tau_0}$ in $\Omega_{\rm FEM}$.
Compute the sequence ${\varepsilon_r}_k, k >0$  with known $\sigma_k$, on a refined meshes $K_{h_k}$
  via following steps:

\STATE Obtain numerical solutions  ${\varepsilon_r}_k$   on $K_{h_k} \times
  J_{\tau_k}$ using
the Algorithm 3 (Conjugate Gradient Method).

\STATE
Refine the mesh $K_{h_k}$ at all points where the second mesh refinement indicator
\begin{equation}
| g_h^k(x)|  \geq \beta_k  \max_{{K_h}_k} | g_h^k(x)|,  \label{62}
\end{equation}
is satisfied. Here,
 indicator  $g_h^k$  is defined in
(\ref{gradient1h}).  Tolerance
number $\beta_k \in \left( 0,1\right) $ should be chosen in numerical examples.

\STATE Define a new refined mesh as $K_{h_{k+1}}$ and
  construct a new time partition $J_{\tau_{k+1}}$ such that the CFL
  condition 
  is satisfied. Interpolate ${\varepsilon_r}_k, \sigma_k$ on a new mesh
  $K_{h_{k+1}}$ and perform steps 1-3 on the space-time mesh $K_{h_{k+1}}
  \times J_{\tau_{k+1}}$. Stop mesh refinements when
  $||{\varepsilon_r}_k - {\varepsilon_r}_{k-1}|| < tol_1$,
  or $|| g_h^k(x)|| < tol_2$, where $tol_i, i=1,2$ are tolerances chosen by the user.

 \end{algorithmic}
\end{algorithm}

\vspace{0.5cm}

\textbf{Remarks}

\begin{itemize}
  
\item 1. First we make comments how to choose the tolerance numbers
  $\widetilde{\beta_k }, \beta_k$ in (\ref{alg2_2}),  (\ref{62}).
  Their
  values depend on the concrete values of
  $\max \limits_{\Omega_{IN}}
  |h {\varepsilon_r}_k|$ and
  $ \max \limits_{\Omega_{IN}} | g_h^k(x)|$, correspondingly.  If we will take
  values of $\beta_k, \widetilde{\beta_k }$ which are very close to $1$ then we
  will refine the mesh in very narrow region of the $\Omega_{IN}$,
  and if we will choose $\beta_k, \widetilde{\beta_k} \approx 0$ then
  almost all elements in the finite element mesh  will be refined, and thus, we
  will get global and not local mesh refinement. 

\item 2.  To compute $L_2$ norms $||{\varepsilon_r}_k - {\varepsilon_r}_{k-1}||$,
  in step 3 of adaptive algorithms the
    reconstruction ${\varepsilon_r}_{k-1}$ is interpolated from the
  mesh $K_{h_{k-1}}$ to the mesh $K_{h_k}$.

\item 3. The computational mesh is refined only in $\Omega_{\rm FEM}$
  such that no new nodes are added in the overlapping elements between two domains,   $\Omega_{\rm FEM}$ and  $\Omega_{\rm FDM}$. Thus,   the mesh    in $\Omega_{\rm FDM}$,  where finite diffirence method is used,  always  remains  unchanged.

\end{itemize}

\section{Numerical examples}

\label{sec:numex}

In this section, we present numerical simulations of the
reconstruction of permittivity function of three-dimensional
anatomically realistic breast phantom taken from online repository
\cite{wisconsin} using an adaptive reconstruction Algorithm 4 of
section \eqref{sec:ref}. We have tested performance of an adaptive
Algorithm 5 and it is slightly more computationally expensive in terms
of time compared to the performance of Algorithm 4. Additionally,
 relative errors in the reconstructions of  dielectric permittivity function
are slightly
smaller for Algorithm 4 and thus, in this section we present results
of reconstruction for Algorithm 4.

\begin{figure}[tbp]
  \begin{center}
    \begin{tabular}{c}
      \hline \\
      $\varepsilon_r$
      \\
    \begin{tabular}{cc}
  {\includegraphics[scale=0.6, trim = 6.2cm  7.0cm 4.7cm 7.0cm, clip=true,]{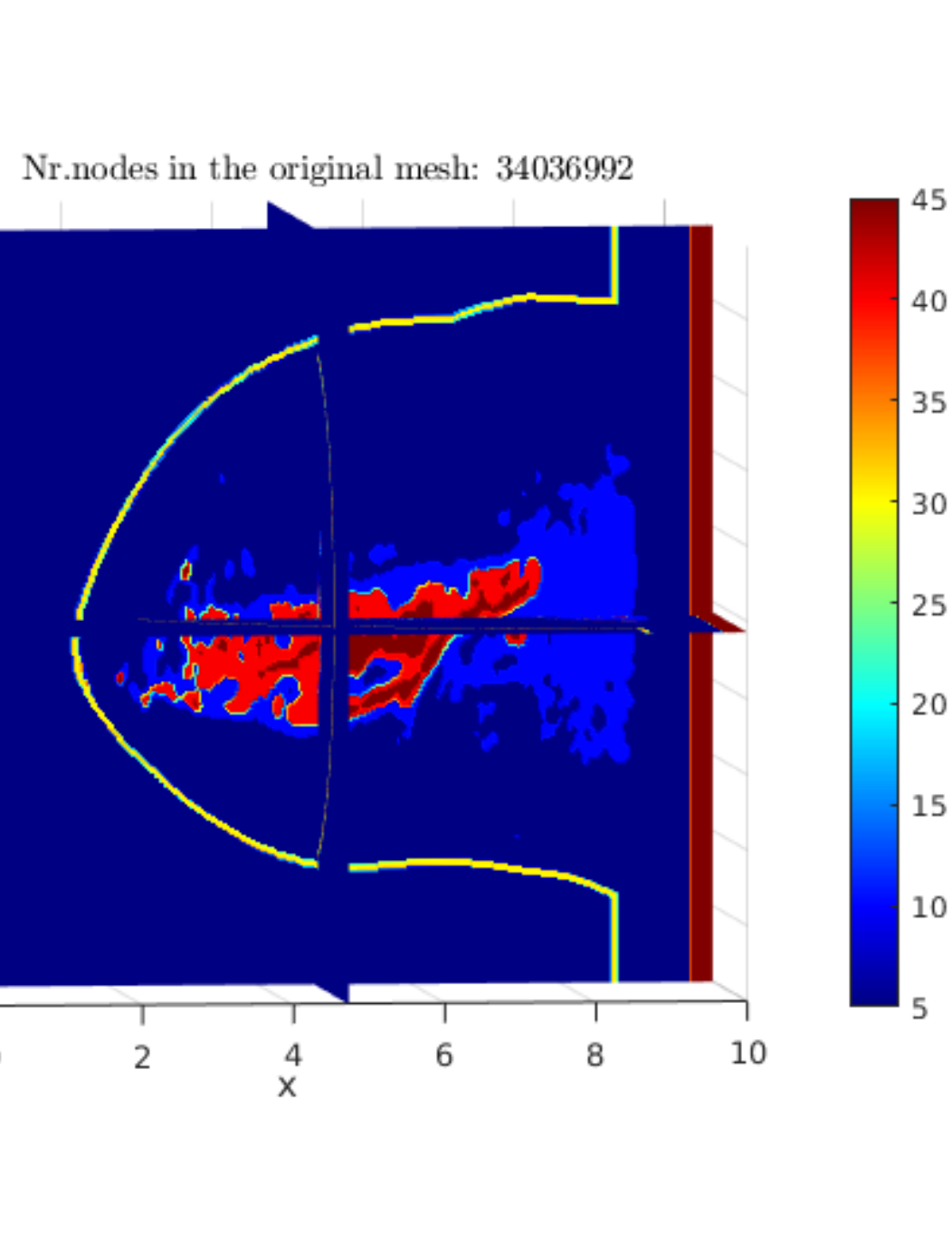}}  &
  {\includegraphics[scale=0.6, trim = 6.4cm 7.0cm 4.7cm 7.0cm, clip=true,]{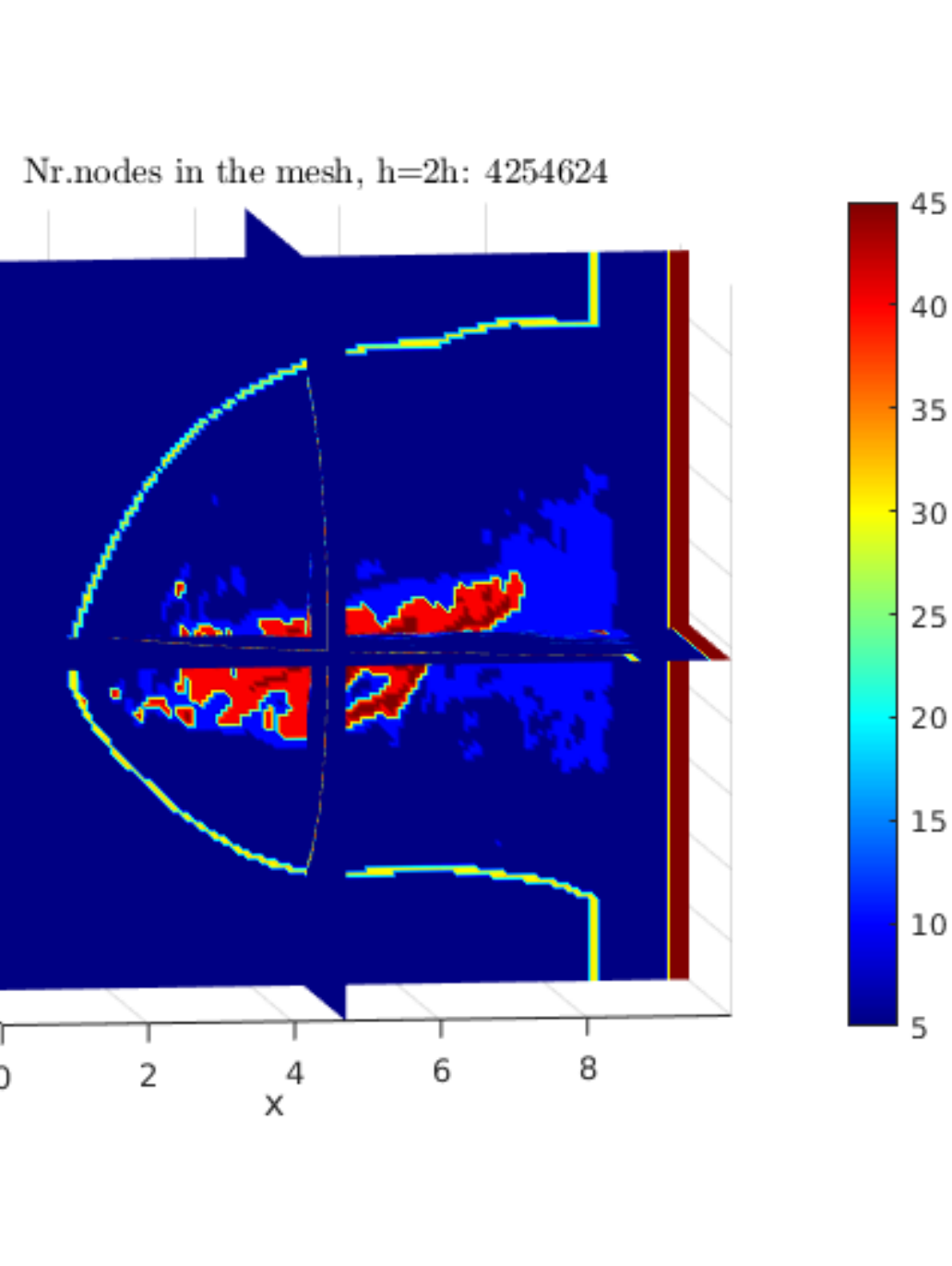}}  \\
  a) & b) \\
  {\includegraphics[scale=0.6, trim = 6.3cm 7.0cm 4.7cm 7.0cm, clip=true,]{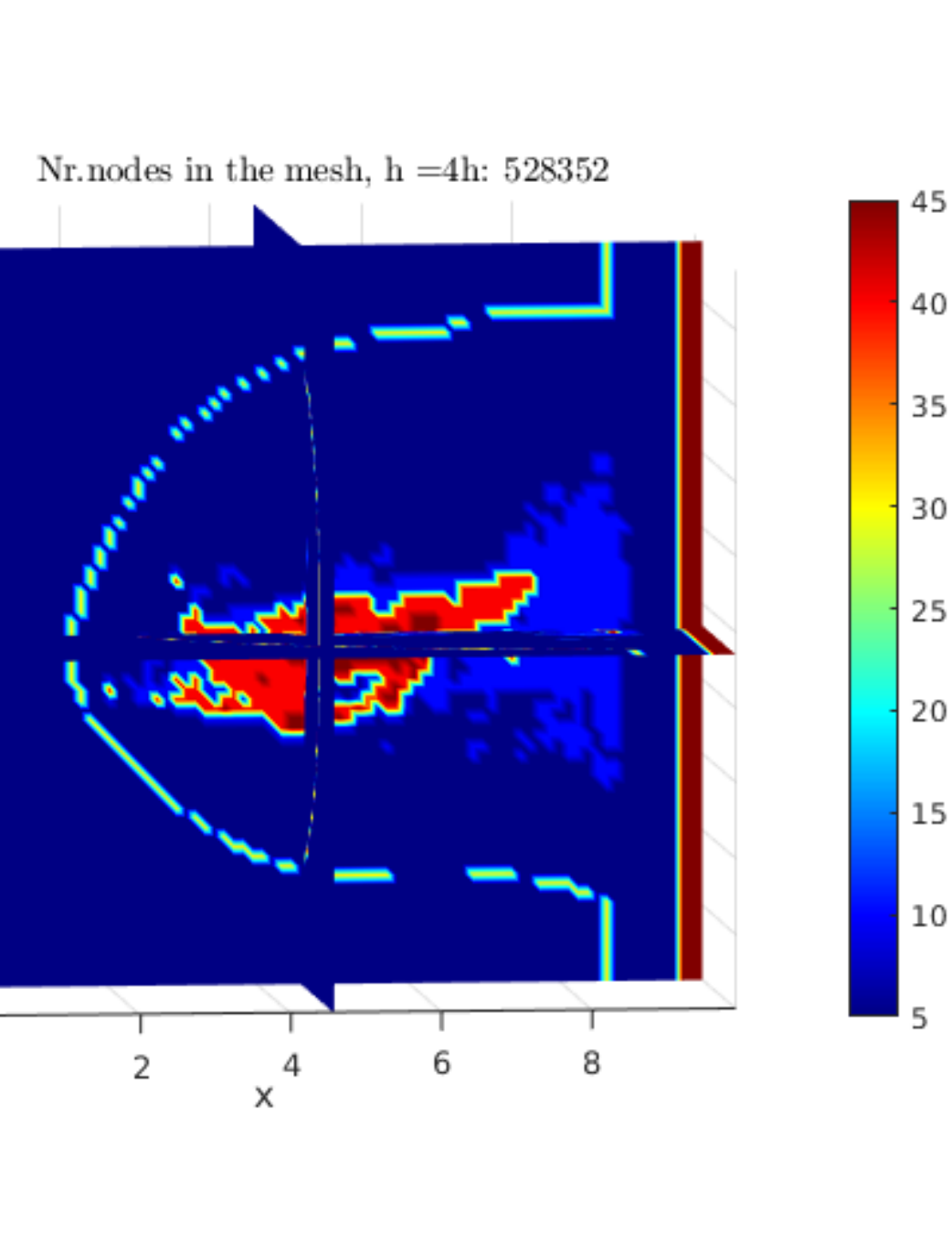}}  &
    {\includegraphics[scale=0.6, trim = 6.2cm 7.0cm 4.7cm 7.0cm, clip=true,]{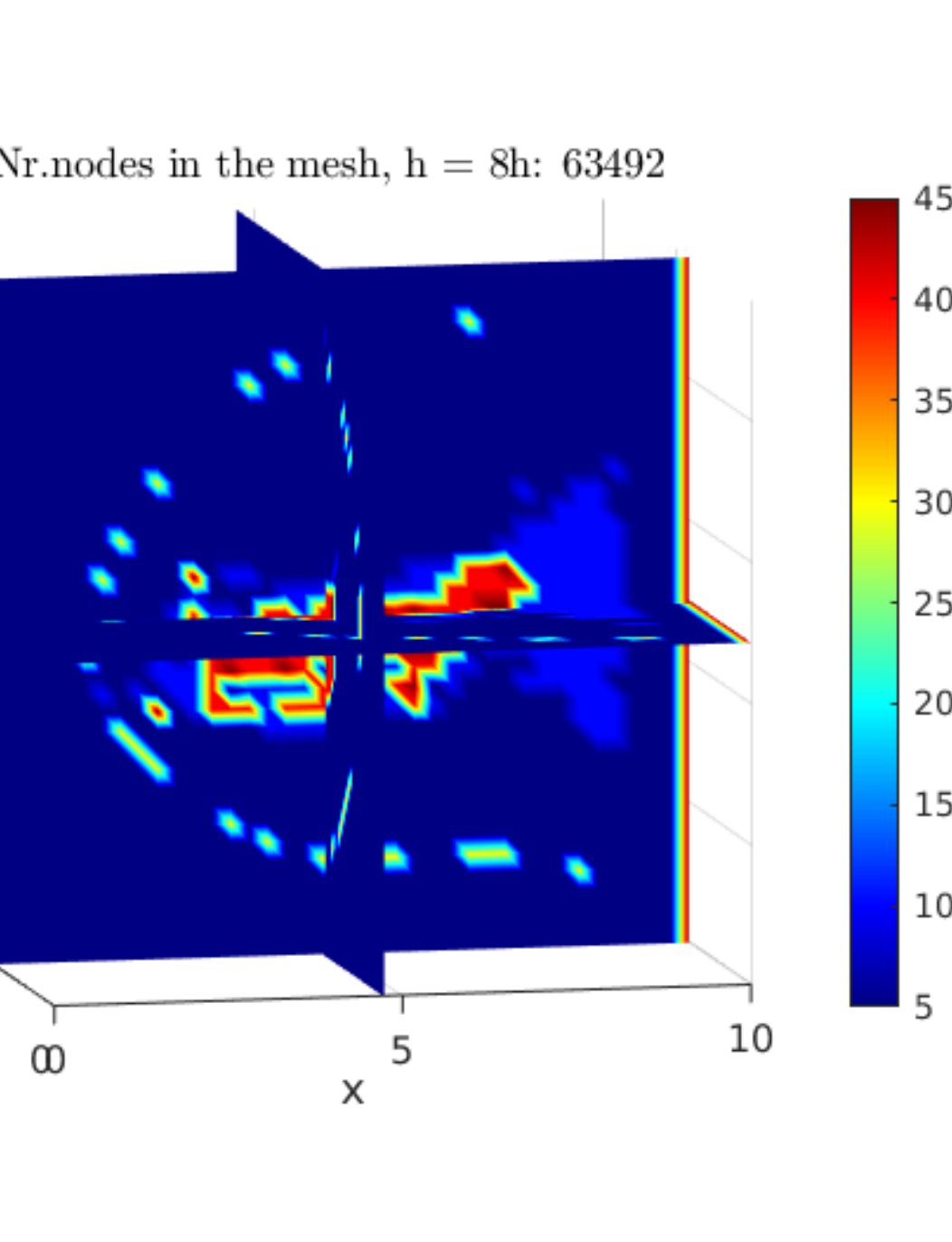}}  \\
    c)  & d) 
    \end{tabular}
    \\
    \hline
    \end{tabular}
 \end{center}
  \caption{\small \emph{\
    Spatial distribution of realistic ultrawideband
   dielectric properties of 3D breast phantom of database
   \cite{wisconsin} developed at the Department of Electrical and
   Computer Engineering at University of Wisconsin-Madison, USA. 
  Figure a) shows  original values of  $\varepsilon_r$
  at  6 GHz for object $ID\_012204$ of database  \cite{wisconsin}.
  Figures b)-d) present sampled version of $\varepsilon_r$.
}}
 \label{fig:eps}
\end{figure}

\begin{figure}
    \centering
    \begin{tabular}{c}
    \hline
    \\ \quad \\
$\sigma$ \\
  \begin{tabular}{cc}
  {\includegraphics[scale=0.6, trim = 6.1cm  7.0cm 4.9cm 7.0cm, clip=true,]{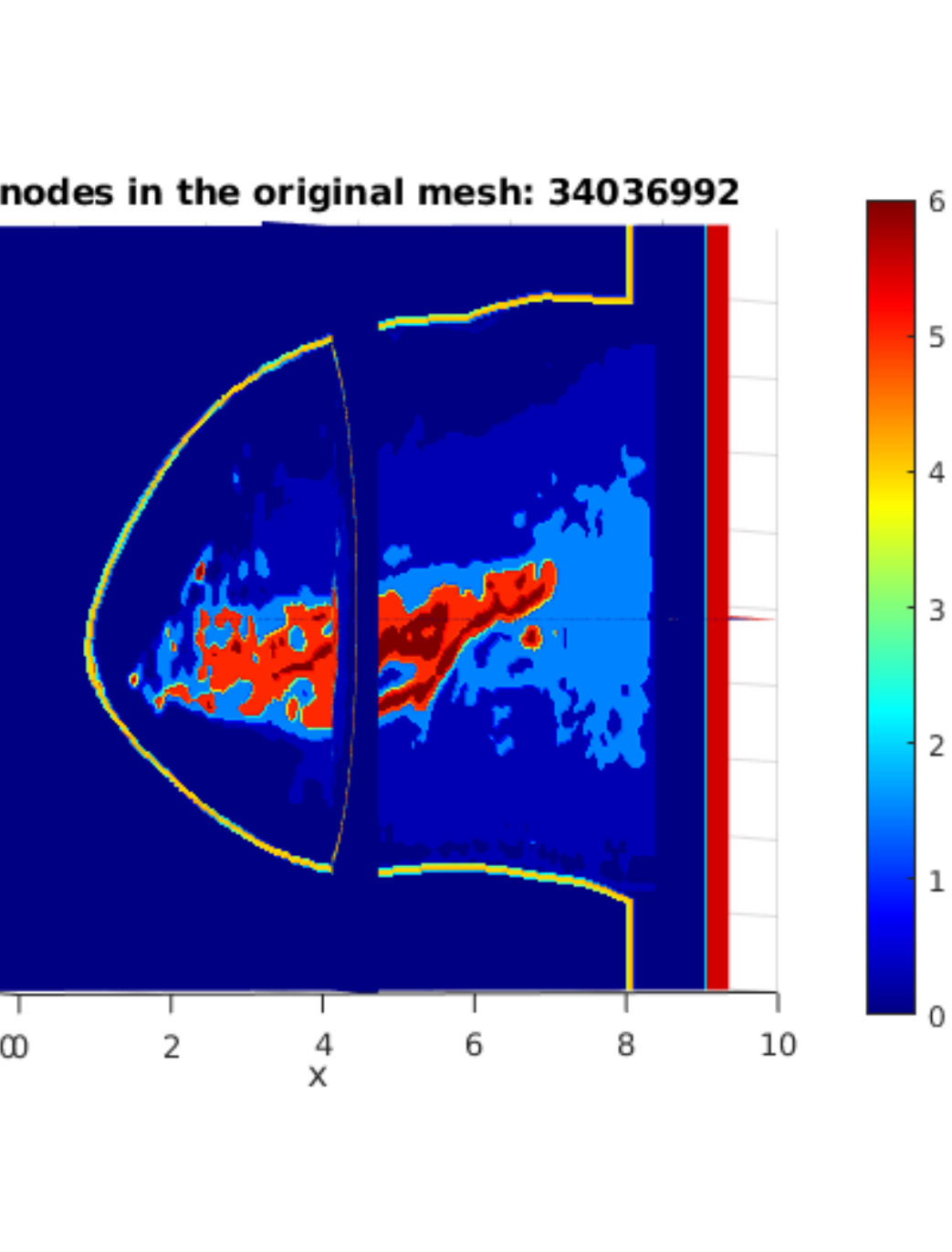}}  &
  {\includegraphics[scale=0.6, trim = 6.35cm 7.0cm 4.9cm 7.0cm, clip=true,]{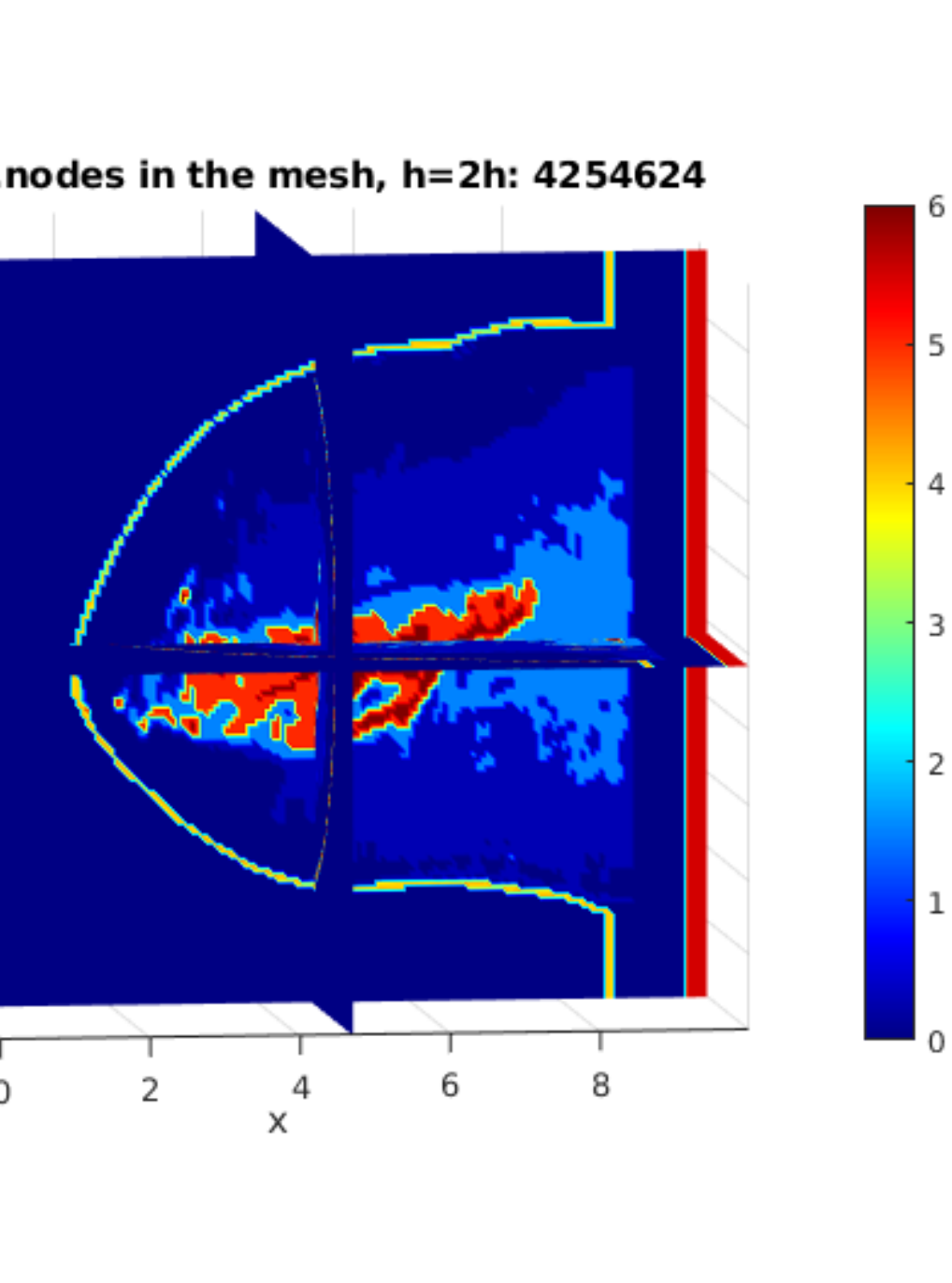}}  \\
  a) & b) \\
  {\includegraphics[scale=0.6, trim = 6.45cm 7.0cm 4.9cm 7.0cm, clip=true,]{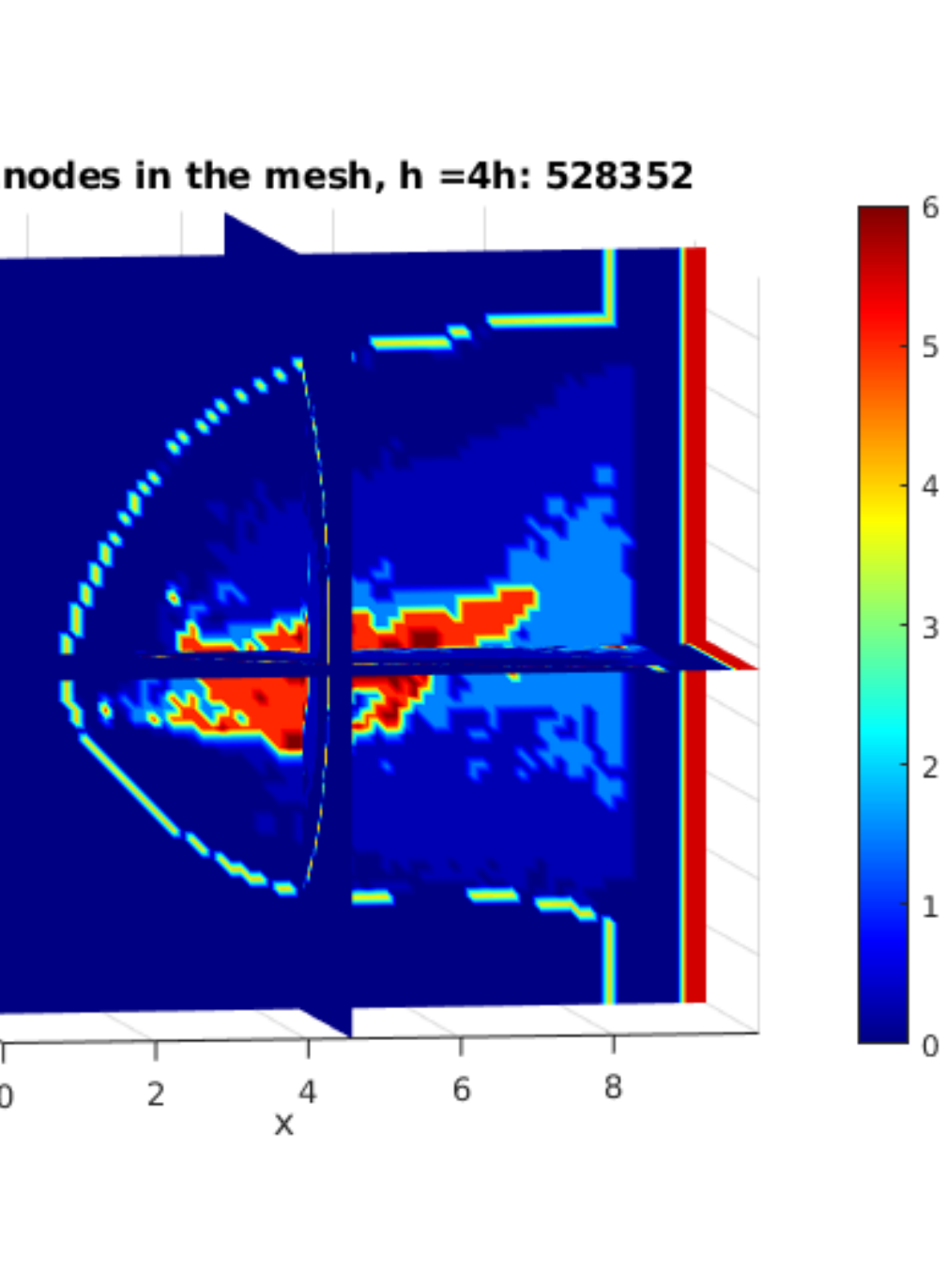}}  &
  {\includegraphics[scale=0.6, trim = 6.5cm 7.0cm 4.9cm 7.0cm, clip=true,]{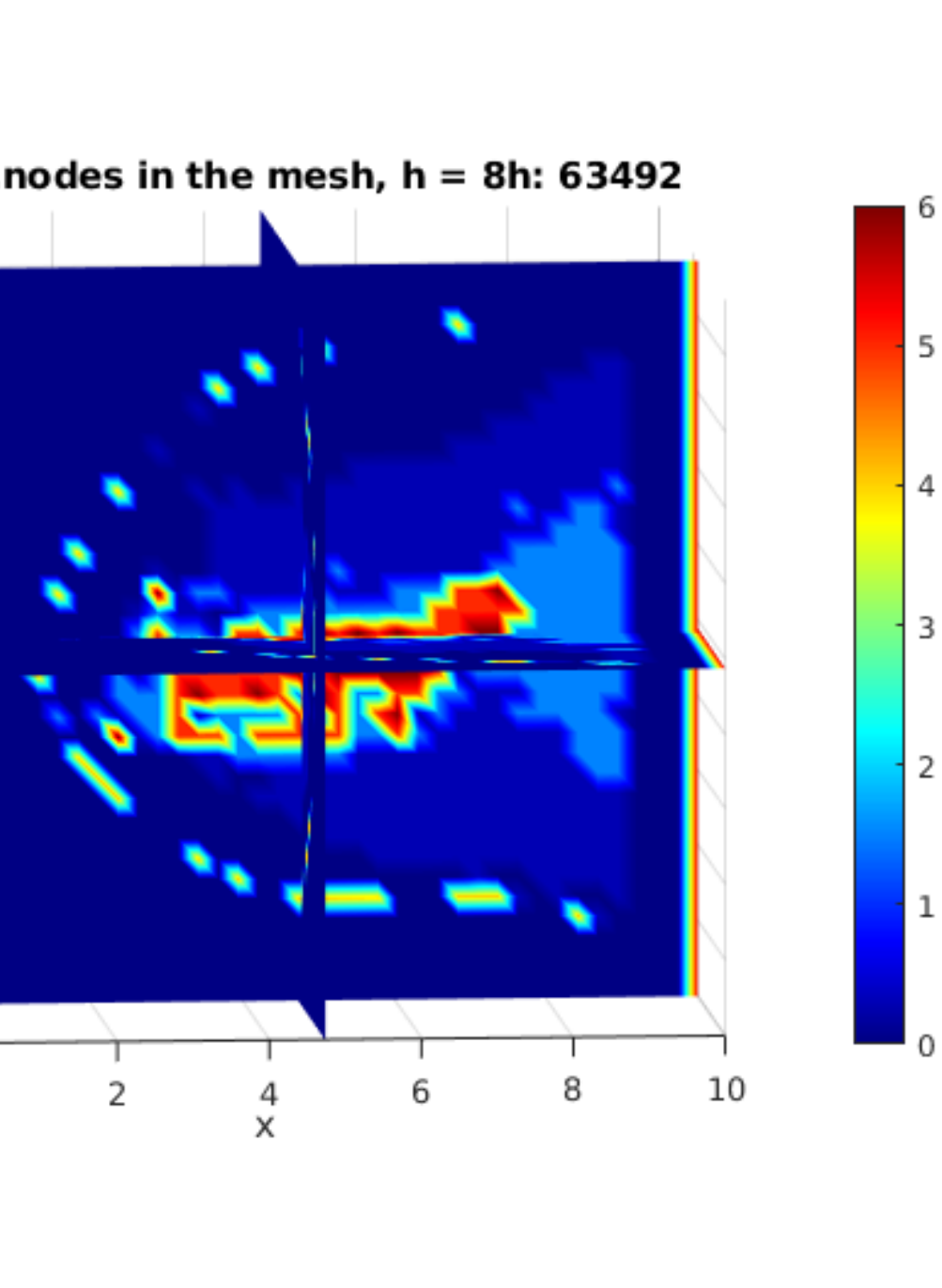}}  \\
  c) & d) 
  \end{tabular}
  \\
  \hline
  \\
  \end{tabular}
    \caption{
      \small \emph{\ Spatial distribution of realistic ultrawideband
   dielectric properties of 3D breast phantom of database
   \cite{wisconsin} developed at the Department of Electrical and
   Computer Engineering at University of Wisconsin-Madison, USA. 
  Figure a) shows  original values of $\sigma$  (S/m)
  at  6 GHz for object $ID\_012204$ of database  \cite{wisconsin}.
  Figures b)-d) present sampled versions of $\sigma$  (S/m).
}}
 \label{fig:sigma}
\end{figure}

\begin{figure}[t]
 \begin{center}
    \begin{tabular}{c}
    \\
   \begin{tabular}{cc}
     {\includegraphics[scale=0.5, trim = 1.0cm 0.0cm 0.0cm 0.0cm,clip=true,]{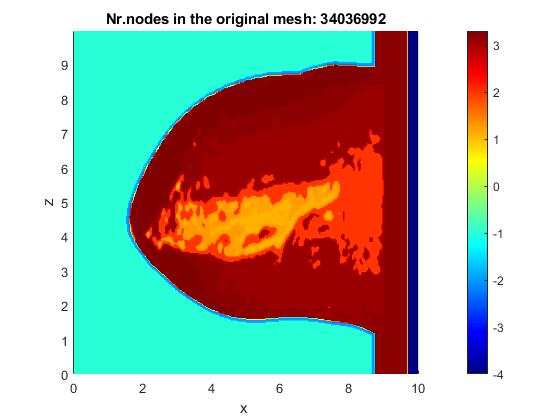}} &
     {\includegraphics[scale=0.5, trim = 1.0cm 0.0cm 0.0cm 0.0cm, clip=true,]{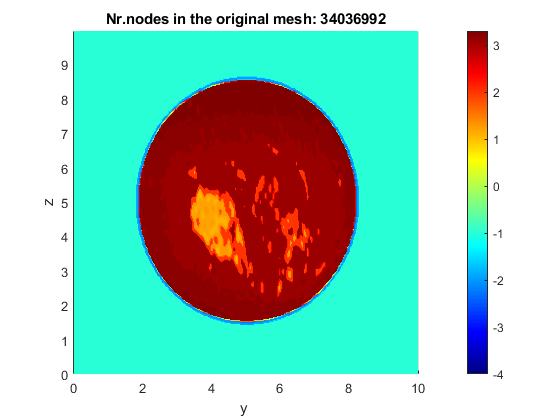}} \\
     a) $x_1 x_3$ view & b) $x_2 x_3$ view  \\
     \hline
     \\
    {\includegraphics[scale=0.55, trim = 2cm 0.0cm 0.0cm 0.0cm, clip=true,]{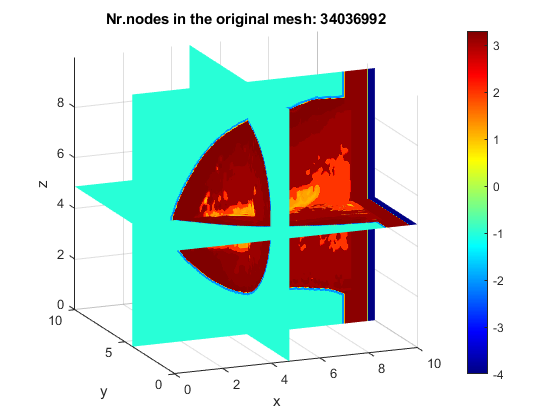}} & {\includegraphics[scale=0.55, trim = 2cm 0.0cm 0.0cm 0.0cm, clip=true,]{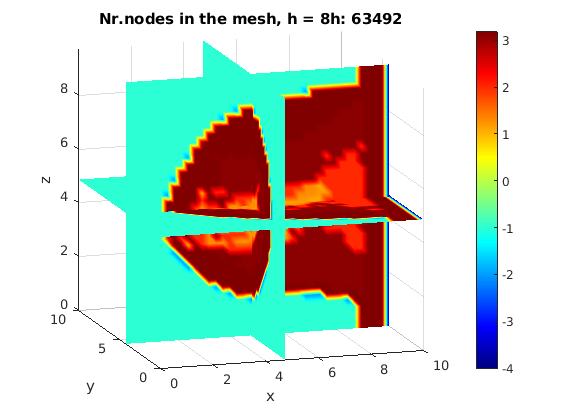}} \\
    c) prospect view   & d) prospect view
\end{tabular}
\end{tabular}
 \end{center}
 \caption{
   \small \emph{\ 
 a)-c) Original values and d) sampled values of the spatial distribution of media numbers  of Table \ref{tab:table1}
 for breast phantom of object $ID\_012204$ of database   \cite{wisconsin}.
 Table   \ref{tab:table1}  clarifies description of  media numbers and corresponding tissue types.
 }}
\label{fig:typemat}
\end{figure}

\begin{table}[ht]
    \centering
    \begin{tabular}{| m{15em} | m{1.4cm}| m{1.4cm}| m{1.4cm}| m{1.4cm} | m{1.4cm}| }
    \hline
    & Media  & Test 1 &  Test 1 &  Test 2 & Test 2 \\
    Tissue type & number & $\varepsilon_r/5$  & $\sigma/5$  & $\varepsilon_r/5$  & $\sigma/5$ \\
    \hline
    Immersion medium  & -1            & 1 & 0    & 1 & 0 \\
    \hline
    Skin & -2                         & 1 & 0    & 1 & 0 \\
    \hline
    Muscle            & -4            & 1  & 0   & 1 & 0 \\
    \hline
    Fibroconnective/glandular-1 & 1.1 & 9 &  1.2 & 9 &  1.2\\
    \hline
    Fibroconnective/glandular-2 & 1.2 &  8 & 1   & 1 &   0 \\
    \hline
    Fibroconnective/glandular-3 & 1.3 &  8 & 1   & 1 &  0  \\
    \hline
    Transitional                 & 2  &  1 & 0   & 1 &  0  \\
    \hline
    Fatty-1 & 3.1                     & 1 & 0    & 1 &  0 \\
    \hline
    Fatty-2 & 3.2                     & 1 & 0    & 1 &  0 \\
    \hline
    Fatty-3 & 3.3                     & 1 & 0    & 1 &  0  \\
    \hline
    \end{tabular}
    \caption{\textit{
        Tissue types and corresponding media numbers  of database \cite{wisconsin}  together with realistic
      weighted  values of $\varepsilon_r$  and $\sigma$  (S/m)  for breast phantom  with  ID=012204
      used in numerical experiments of section \ref{sec:setup}.
      Figure \ref{fig:typemat}  presents media numbers of this table on original and sampled meshes.
    }}
    \label{tab:table1}
\end{table}

\subsection{Description of anatomically realistic data}

We have tested our reconstruction algorithm using three-dimensional
realistic breast phantom with ID = 012204 provided in the online
repository \cite{wisconsin}. The phantom comprises the structural
heterogeneity of normal breast tissue for realistic dispersive
properties of normal breast tissue at 6 GHz reported in \cite{laz1, laz2}.
The breast phantoms of database \cite{wisconsin}  are derived using T1-weighted MRIs
of patients in prone position. Every  phantom presents 3D mesh of cubic voxels of the size
 $0.5 \times 0.5 \times 0.5 $ mm.

Tissue types and corresponding media numbers of breast phantoms 
 are taken from
 \cite{wisconsin} and  are given in Table \ref{tab:table1}.
 Spatial distribution of these media
numbers   for phantom with ID = 012204    is presented in Figure
\ref{fig:typemat}.
Figures  \ref{fig:typemat}-a)-c)

 demonstrate 
distribution of media numbers  on the original coarse mesh consisting of
34 036 992 nodes. Clearly, performing computations on a such big mesh is
computationally demanding task, and thus, we have sampled the original mesh.
In all our computations we have used the   mesh
  consisting of 63492 nodes as a coarse finite element mesh which was obtained by
taking every 8-th node in $x_1, x_2$ and $x_3$ directions of the
original mesh.
 Figures \ref{fig:eps}-\ref{fig:sigma}
 shows spatial distribution of dielectric permittivity
 $\varepsilon_r$ and
  effective conductivity  $\sigma$ (S/m)  on original and sampled meshes.
  
\begin{figure}[tbp]
  \begin{center}
    \begin{tabular}{cc}
      {\includegraphics[scale=0.4, trim = 2.0cm 0.0cm 0.8cm 0.0cm, clip=true,]{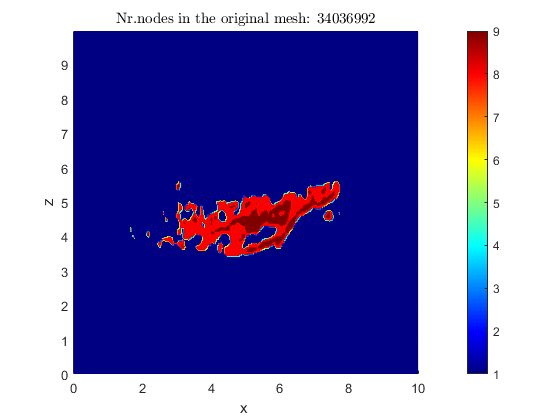}} &
      {\includegraphics[scale=0.4, trim = 2.0cm 0.0cm 0.8cm 0.0cm, clip=true,]{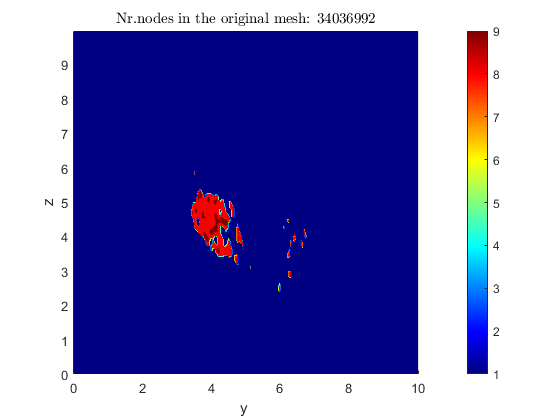}} \\
      a) $x_1 x_3$ view & b)  $x_2 x_3$  view \\
        {\includegraphics[scale=0.4, trim = 2.0cm 0.0cm 0.8cm 0.0cm, clip=true,]{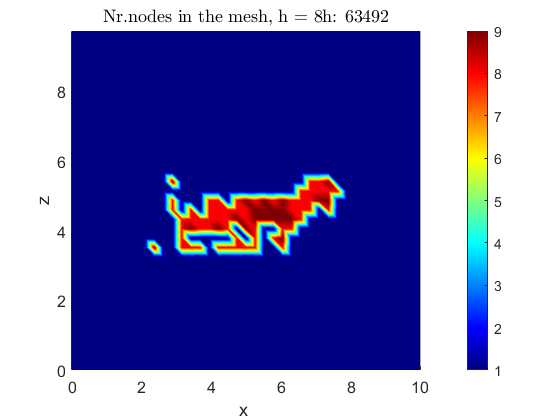}} &
      {\includegraphics[scale=0.4,  trim = 2.0cm 0.0cm 0.8cm 0.0cm, clip=true,]{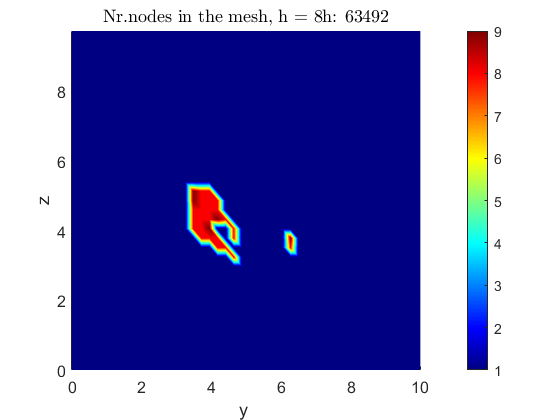}} \\
          c)  $x_1 x_3$  view & d)  $x_2 x_3$  view \\
    \end{tabular}
  \end{center}
  \caption{\small\emph{ Test 1.  Slices of  weighted exact $\varepsilon_r$, see Table   \ref{tab:table1}  for description of different tissue types  and values of  weighted $\varepsilon_r$.  a), b) Slices on original mesh  with mesh size $h$. c), d) Slices on sampled mesh  with mesh size $8h$.}}
  \label{fig:exacteps8h}
\end{figure}

Figure  \ref{fig:typemat}-d) demonstrates distribution of media numbers on finally sampled mesh.
Figure  \ref{fig:exacteps8h}   presents spatial distribution of  weighted values of  $\varepsilon_r$  on original and finally sampled mesh   for Test 1.
Testing of our algorithms on
  other sampled meshes is computationally expensive task,
requiring running of programs in parallel
infrastructure, and can be considered as a topic for future research.

We note that in  all our computations we scaled original values of $\varepsilon_r$ and $\sigma$ of database \cite{wisconsin}  presented in
Figures \ref{fig:eps}-\ref{fig:sigma}
and considered  weighted versions of these parameters,
 in order to satisfy
 conditions \eqref{2.3}
 as well as for efficient implementation of  FE/FD DDM for solution of forward and adjoint problems.
 Table \ref{tab:table1} presents weighted values of $\varepsilon_r$ and $\sigma$
 used in numerical tests of this section.
 Thus, in this way 
we get computational set-up corresponding to the   domain decomposition method which was used in Algorithms 2-5.

\begin{figure}[tbp]
  \begin{center}
    \begin{tabular}{ccc}
     Test 1 &   Test 2  \\
      {\includegraphics[scale=0.3, trim = 8.0cm 0.0cm 8.0cm 0.0cm, clip=true,]{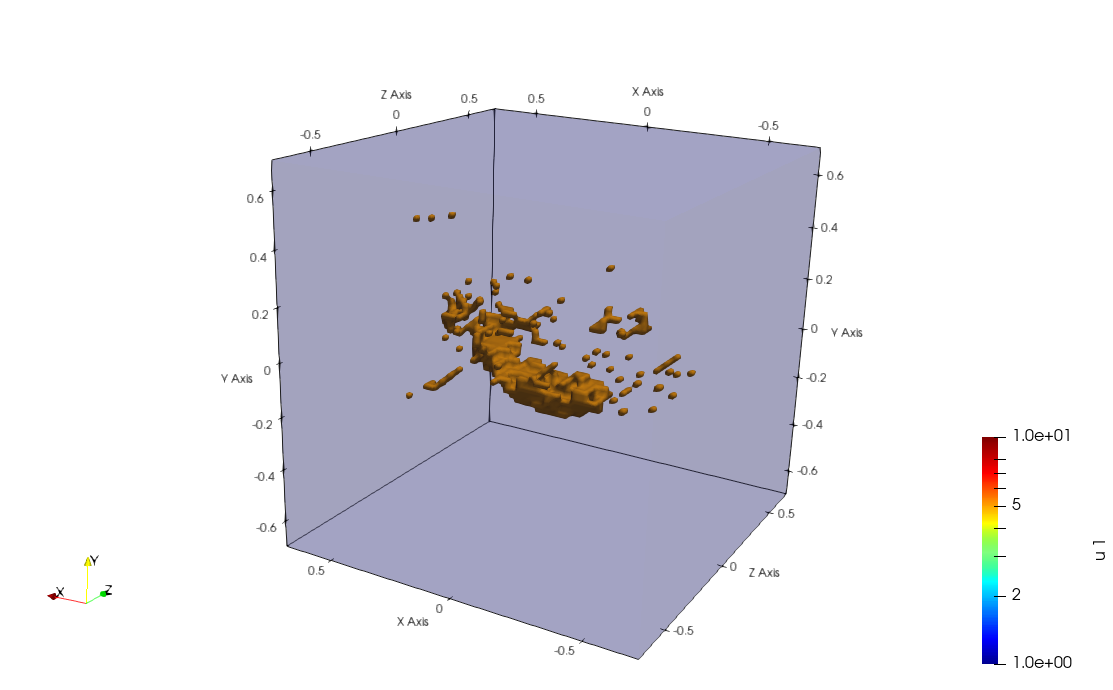}} &
      {\includegraphics[scale=0.3, trim = 8.0cm 0.0cm 8.0cm 0.0cm, clip=true,]{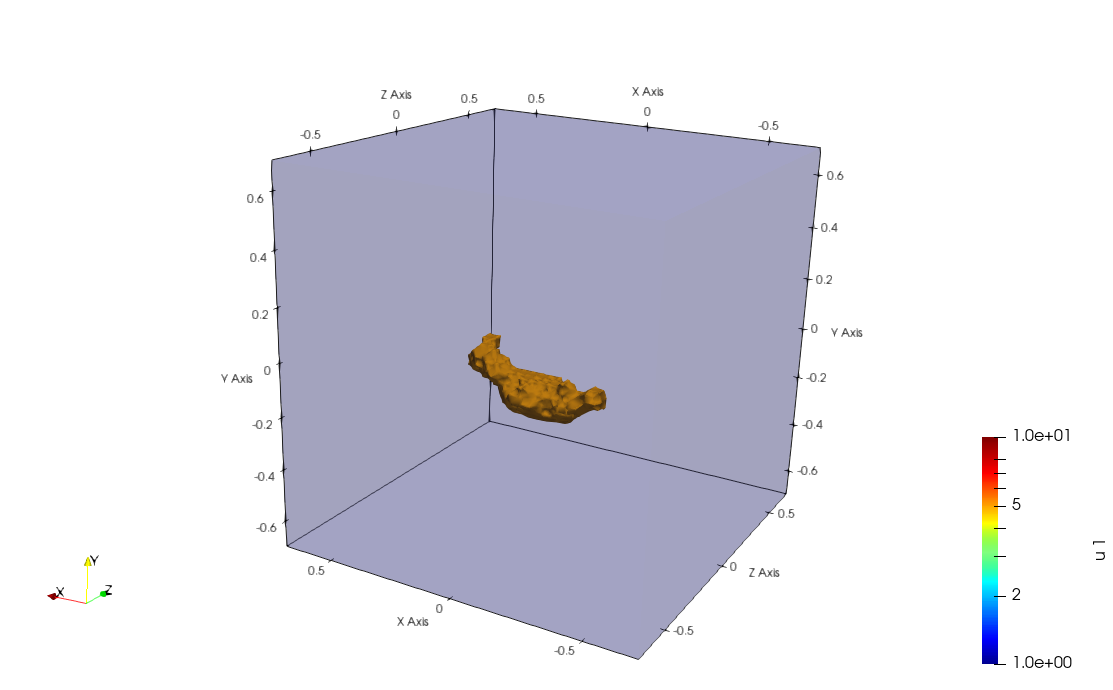}} \\
       a) & b)
    \end{tabular}
  \end{center}
  \caption{\small\emph{ Isosurface of weighted exact dielectric permittivity
    with value $\varepsilon_r \approx 5$    corresponding to   tissue type ``fibroconnective/glandular-1'': a) In Test 1  and b) In Test 2.  Table   \ref{tab:table1}  clarifies description of different tissue types. }}
  \label{fig:exacteps}
\end{figure}

\subsection{Computational set-up}

\label{sec:setup}

We have used the domain decomposition Algorithm  2 of section \ref{sec:hybalg}
to solve forward and adjoint problems  in the adaptive reconstruction Algorithm  4. To do this, we
  set the  dimensionless computational domain $\Omega$   as
  \begin{equation*}
    \Omega = \left\{ x= (x_1,x_2, x_3) \in
           (-0.8840,  0.8824)
    \times (-0.8630,  0.8648)
    \times (-0.8945,  0.8949)
    \right\},
 \end{equation*}
and  the  domain $\Omega_{\rm FEM}$ as
 \begin{equation*}
   \Omega_{\rm FEM}  = \left\{ x= (x_1,x_2,x_3) \in
          (-0.7, 0.6984)
   \times (-0.7, 0.7018)
   \times (-0.7, 0.7004)
   \right\}.
 \end{equation*}
 We choose the coarse mesh sizes $h_1 = 0.0368, h_2 = 0.0326, h_3 = 0.0389$  in
 $x_1, x_2, x_3$  directions, respectively, 
 in $\Omega = \Omega_{\rm FEM} \cup
 \Omega_{\rm FDM}$, as well as in the overlapping regions between $\Omega_{\rm FEM}$
 and  $\Omega_{\rm FDM}$.
  Corresponding physical domains in meters are $\tilde{\Omega} = 0.17664 \times 0.17278 \times 0.17894$ m  for $\Omega$
  and $\tilde{\Omega}_{ \rm FEM} = 0.13985 \times 0.14018 \times 0.14004$  m  for $\Omega_{\rm FEM}$.

The
boundary $\partial \Omega$ of the domain $\Omega$ is  decomposed into
 three different parts and is such that $\partial
\Omega =\partial _{1} \Omega \cup \partial _{2} \Omega \cup \partial
_{3} \Omega$ where $\partial _{1} \Omega$ and $\partial _{2} \Omega$
are, respectively,   front and   back sides of $\Omega$, and $\partial
_{3} \Omega$ is the union of left, right, top and bottom sides of this
domain. We will collect time-dependent observations  at $\Gamma_2 :=
\partial_2 \Omega \times (0,T)$, or at the transmitted side $\partial_2
\Omega$ of $\Omega$.  We also define $\Gamma_{1,1} := \partial_1
\Omega \times (0,t_1]$, $\Gamma_{1,2} := \partial_1 \Omega \times
(t_1,T)$, and $\Gamma_3
:= \partial_3 \Omega \times (0, T)$.

 The following
 model problem was used in  all computations:
\begin{equation}\label{model1}
\begin{split}
 \frac{1}{c^2}   \varepsilon_r \frac{\partial^2 E}{\partial t^2} +
 \nabla ( \nabla \cdot E) - \triangle E  -
 \varepsilon_0 \nabla  (\nabla \cdot ( \varepsilon_r  E))
  &= - \mu_0\sigma \frac{\partial E}{\partial t}~ \mbox{in}~~ \Omega_T, \\
  E(x,0) = 0, ~~~\frac{\partial E}{\partial t}(x,0) &= 0~ \mbox{in}~~ \Omega,     \\
\frac{\partial E}{\partial n}& = f(t) ~\mbox{on}~ \Gamma_{1,1},
\\
\frac{\partial E}{\partial n}& =-\frac{\partial E}{\partial t}  ~\mbox{on}~ \Gamma_{1,2} \cup \Gamma_2,\\
\frac{\partial E}{\partial n} & =0~\mbox{on}~ \Gamma_3.\\
\end{split}
\end{equation}
We initialize   a plane wave  $f(t) = (0,f_2,0)(t)$ for  one component $E_2$ of the electric field
$E=(E_1, E_2, E_3)$  at $\Gamma_{1,1}$
in (\ref{model1}).
The function $f_2(t)$  
 represents the single
direction of a plane wave which is initialized at $\partial_1 \Omega$
in time $t=[0,3.0]$ and is defined as
 \begin{equation}\label{f}
 \begin{split}
 f_2(t) =\left\{ 
 \begin{array}{ll}
 \sin \left( \omega t \right) ,\qquad &\text{ if }t\in \left( 0,\frac{2\pi }{\omega }
 \right) , \\ 
 0,&\text{ if } t>\frac{2\pi }{\omega }.
 \end{array}
 \right. 
 \end{split}
 \end{equation}

\begin{figure}[tbp]
  \begin{center}
    \begin{tabular}{ccc}
      $t= 1.2$ & $t = 1.8$ & $t= 2.4$ \\
      {\includegraphics[scale=0.15,  trim = 8.0cm 0.0cm 5.0cm 0.0cm, clip=true,]{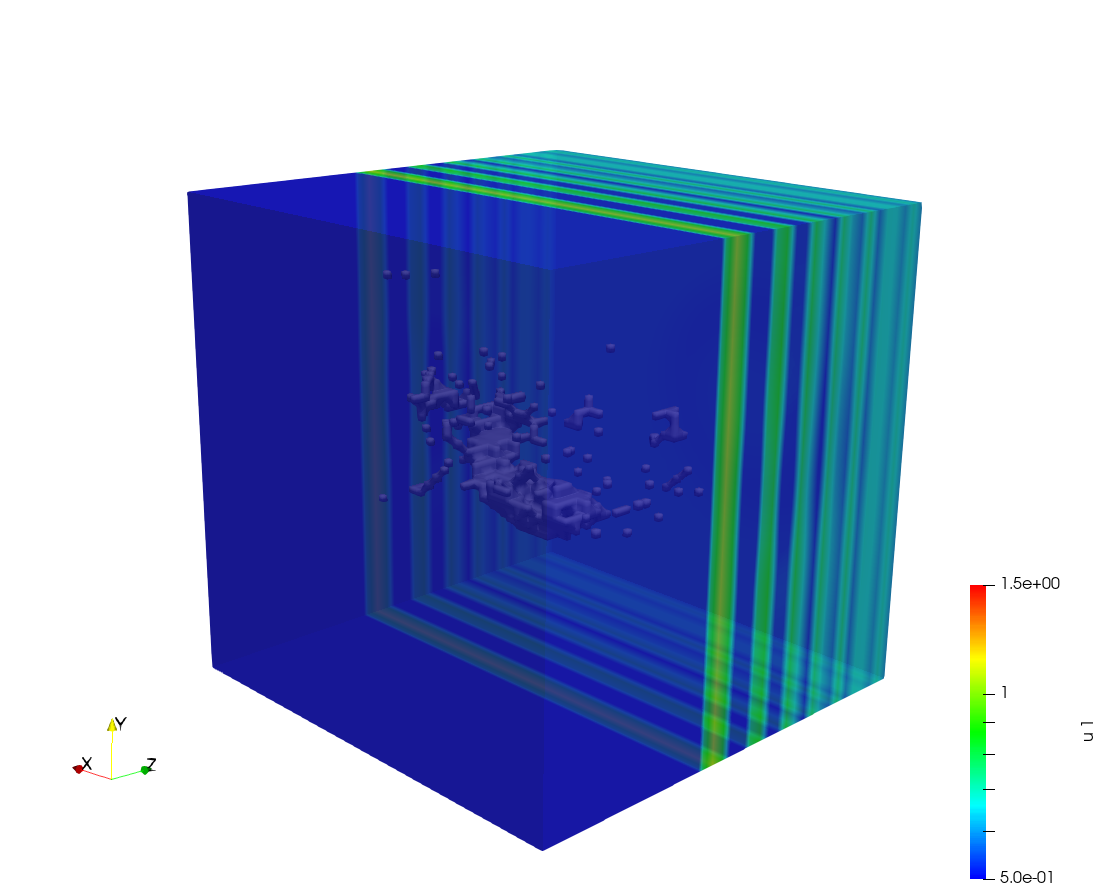}} &
      {\includegraphics[scale=0.15,  trim = 8.0cm 0.0cm 5.0cm 0.0cm, clip=true,]{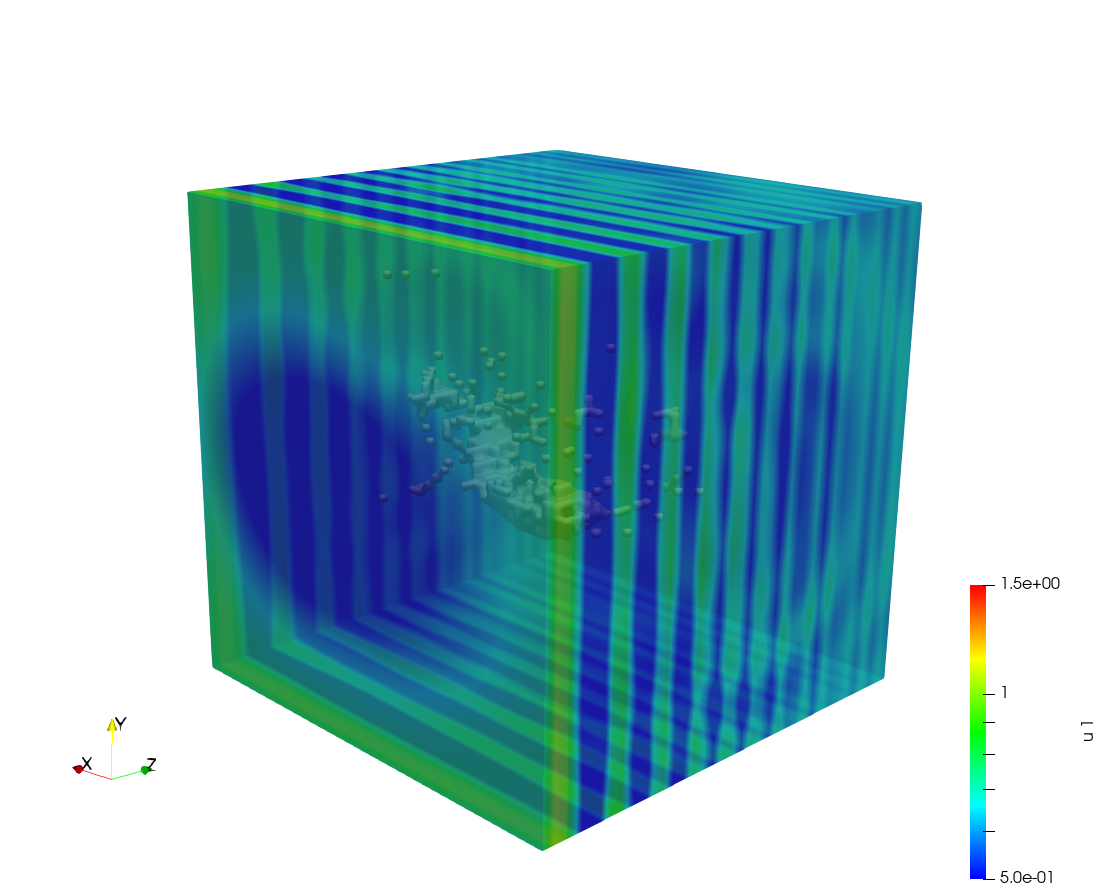}} &
       {\includegraphics[scale=0.15,  trim = 8.0cm 0.0cm 5.0cm 0.0cm, clip=true,]{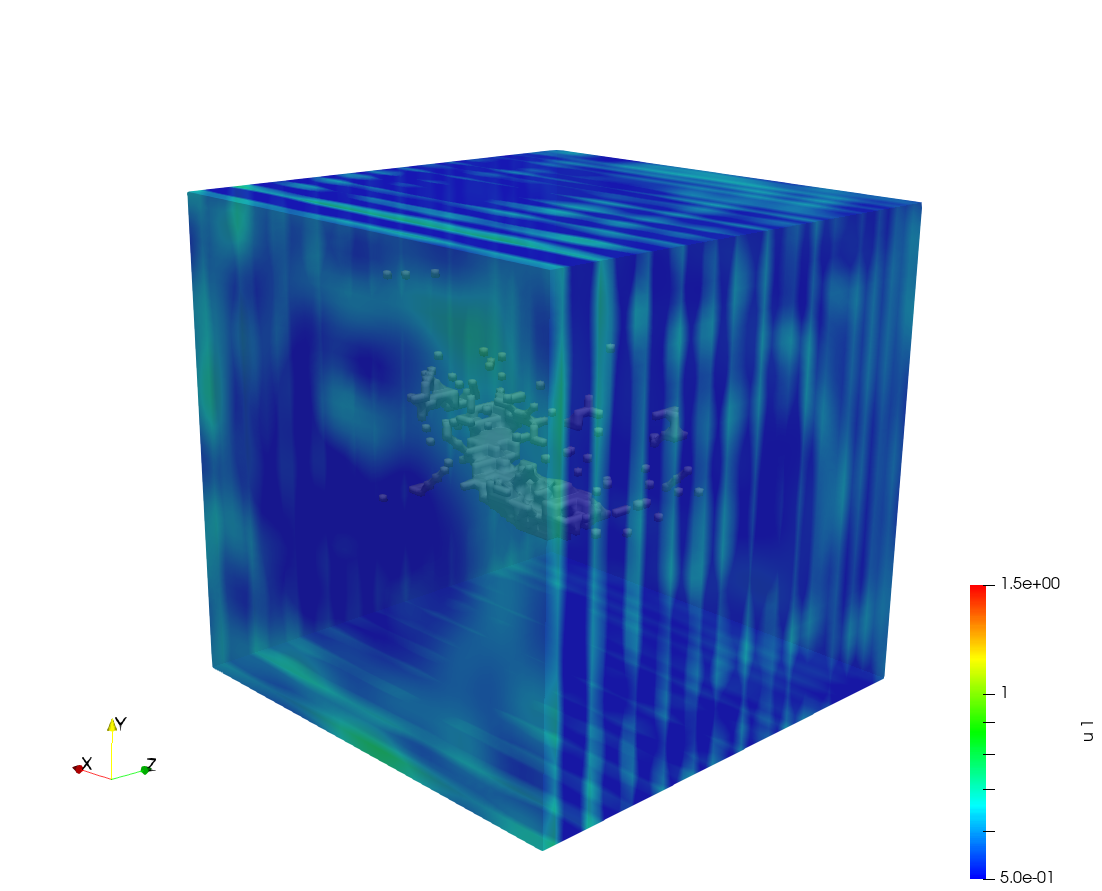}} \\
       a)  &  b)  &  c)
       \\
  \end{tabular}
     \begin{tabular}{ccc}
  {\includegraphics[scale=0.3, clip=true,]{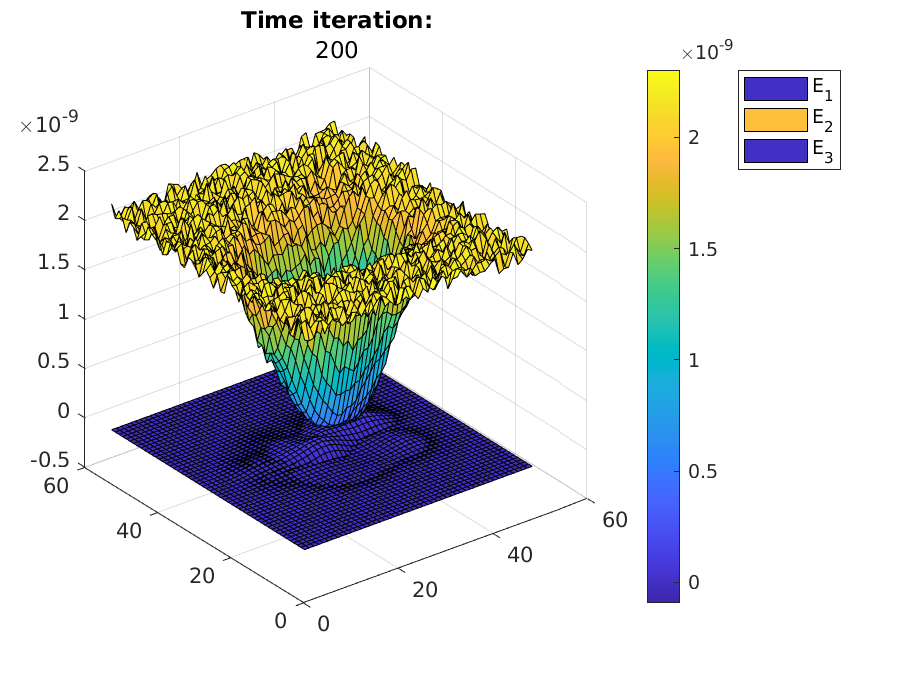}} &
 {\includegraphics[scale=0.3, clip=true,]{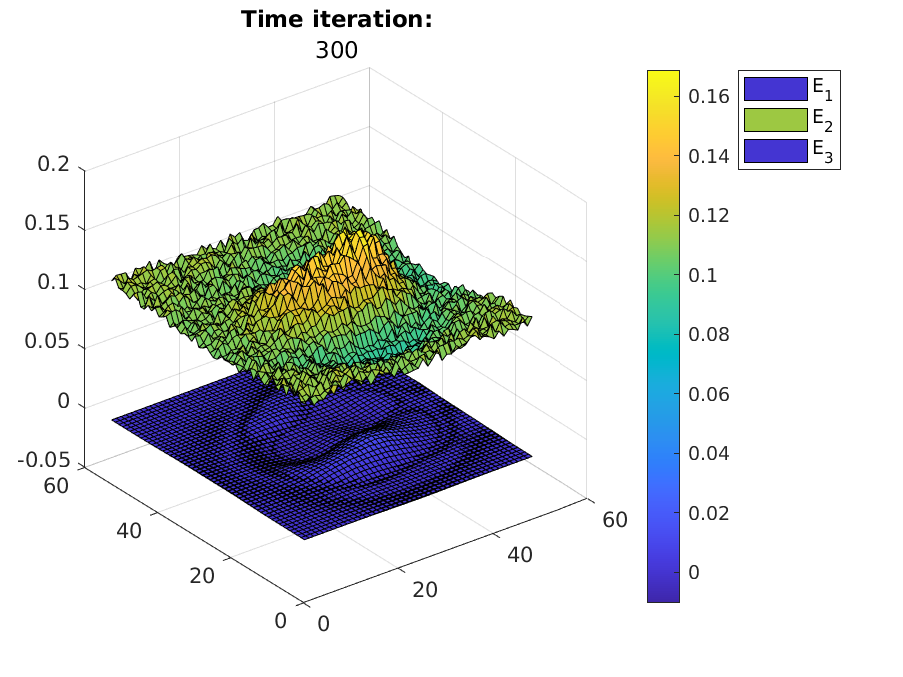}} &
  {\includegraphics[scale=0.3, clip=true,]{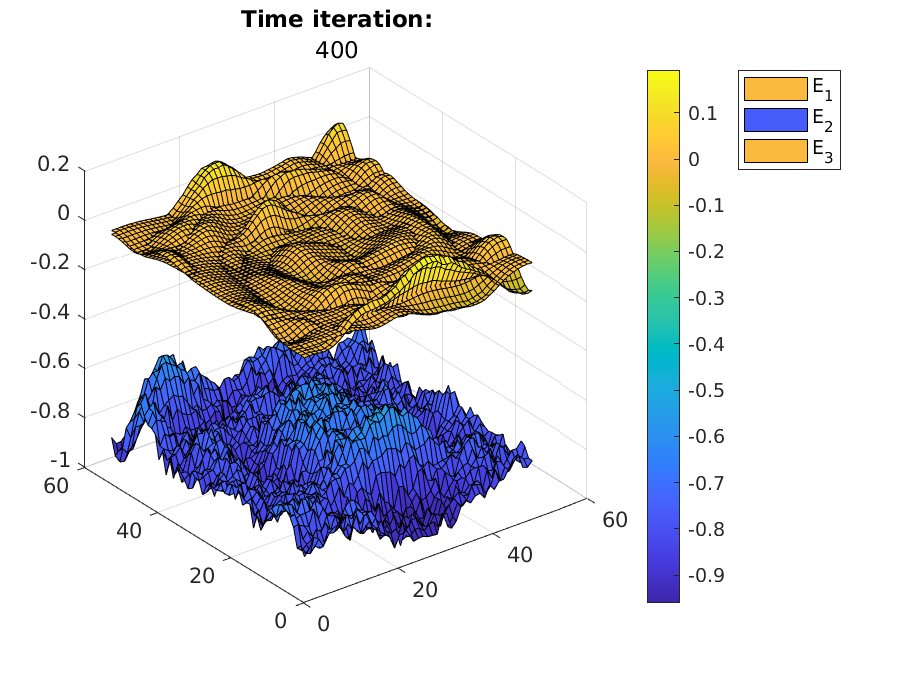}} \\
  d)  &  e)  &  f) \\
   \end{tabular}
\\
     \begin{tabular}{ccc}
      {\includegraphics[scale=0.3, clip=true,]{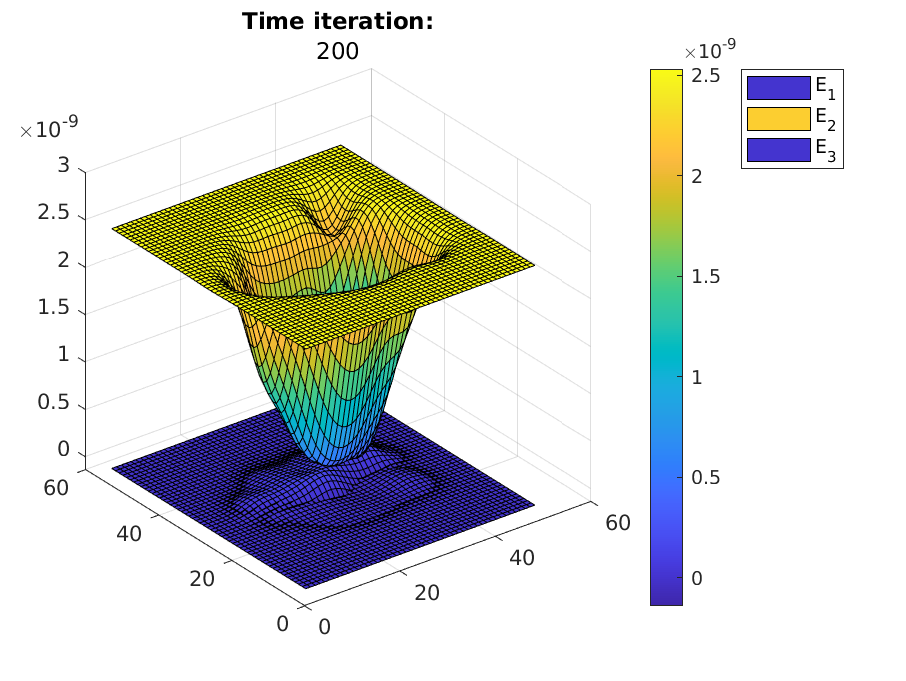}} &
      {\includegraphics[scale=0.3,clip=true,]{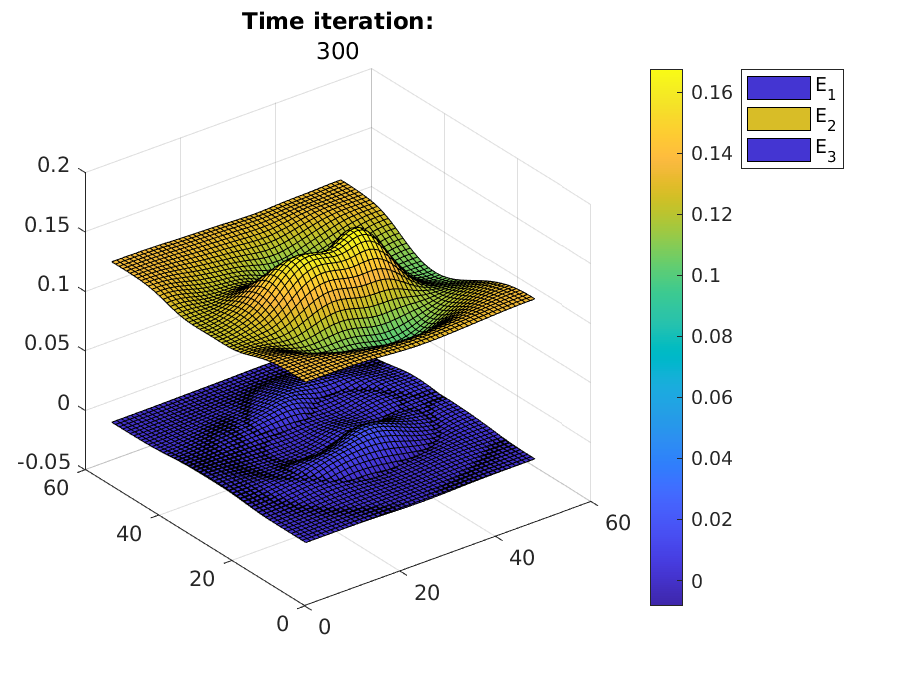}}  &
      {\includegraphics[scale=0.3, clip=true,]{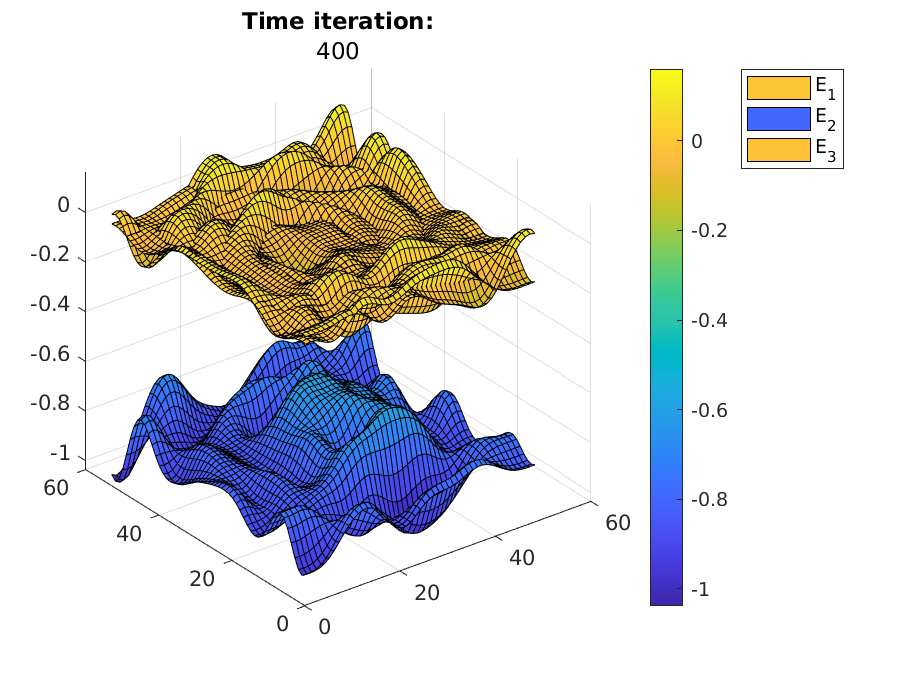}} \\
      g)  &  h)  &  i)
     \end{tabular}
  \end{center}
  \caption{\small\emph{ Test 1. a)-c): Solution $|E_h|$ of model problem \eqref{model1} at different times  for $\omega=40$ in \eqref{f}.
    d)-f): Transmitted noisy scattered data   $E_h = ({E_1}_h, {E_2}_h, {E_3}_h)$   of components of electric field
    $E=(E_1, E_2, E_3)$    at different times.   g)-i): Smoothed transmitted scattered data   $E_h = ({E_1}_h, {E_2}_h, {E_3}_h)$   of components of electric field
    $E=(E_1, E_2, E_3)$    at different times.  The noise level in data is $\delta=10\%$.  }}
  \label{fig:test1datanoise10}
\end{figure}

\begin{figure}[tbp]
  \begin{center}
    \begin{tabular}{ccc}
      $t= 1.2$ & $t = 1.8$ & $t= 2.4$ \\
         {\includegraphics[scale=0.15,  trim = 8.0cm 0.0cm 5.0cm 0.0cm, clip=true,]{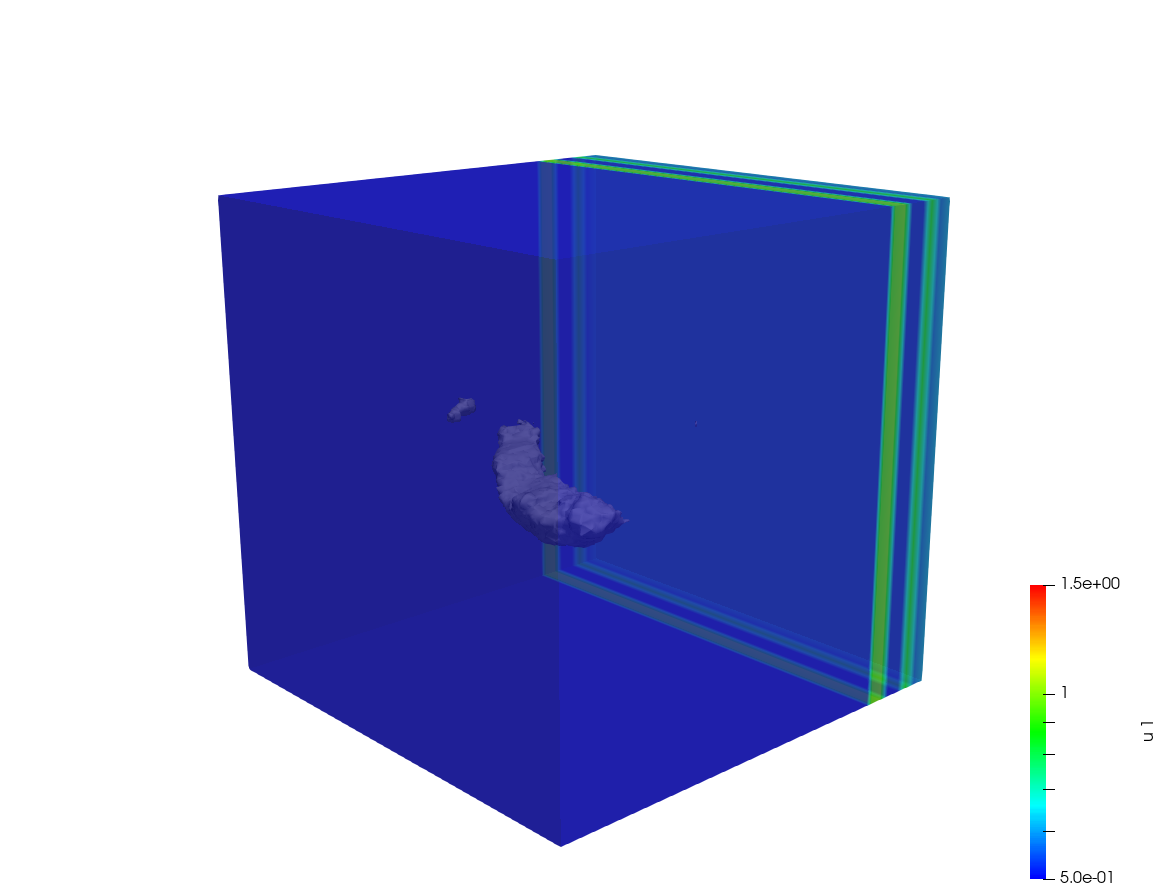}} &
      {\includegraphics[scale=0.15,  trim = 8.0cm 0.0cm 5.0cm 0.0cm, clip=true,]{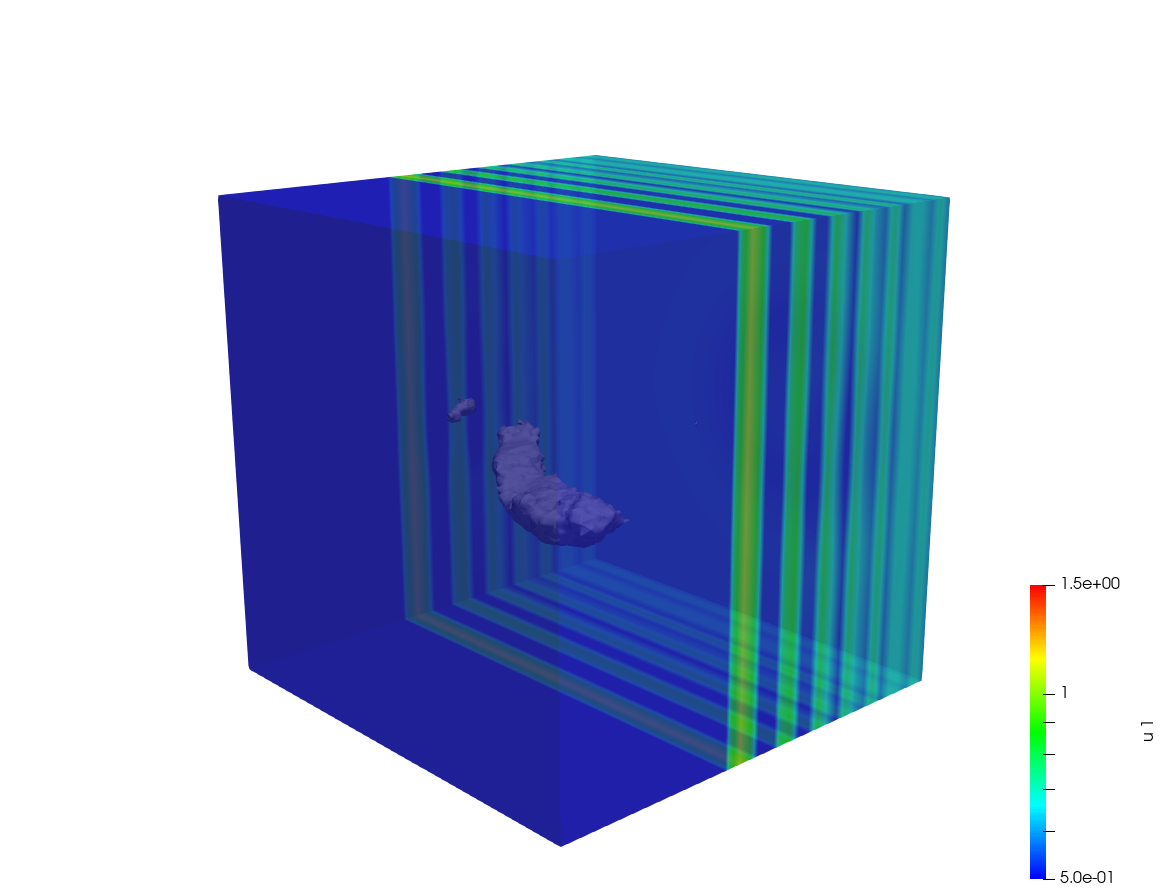}} &
       {\includegraphics[scale=0.15,  trim = 8.0cm 0.0cm 5.0cm 0.0cm, clip=true,]{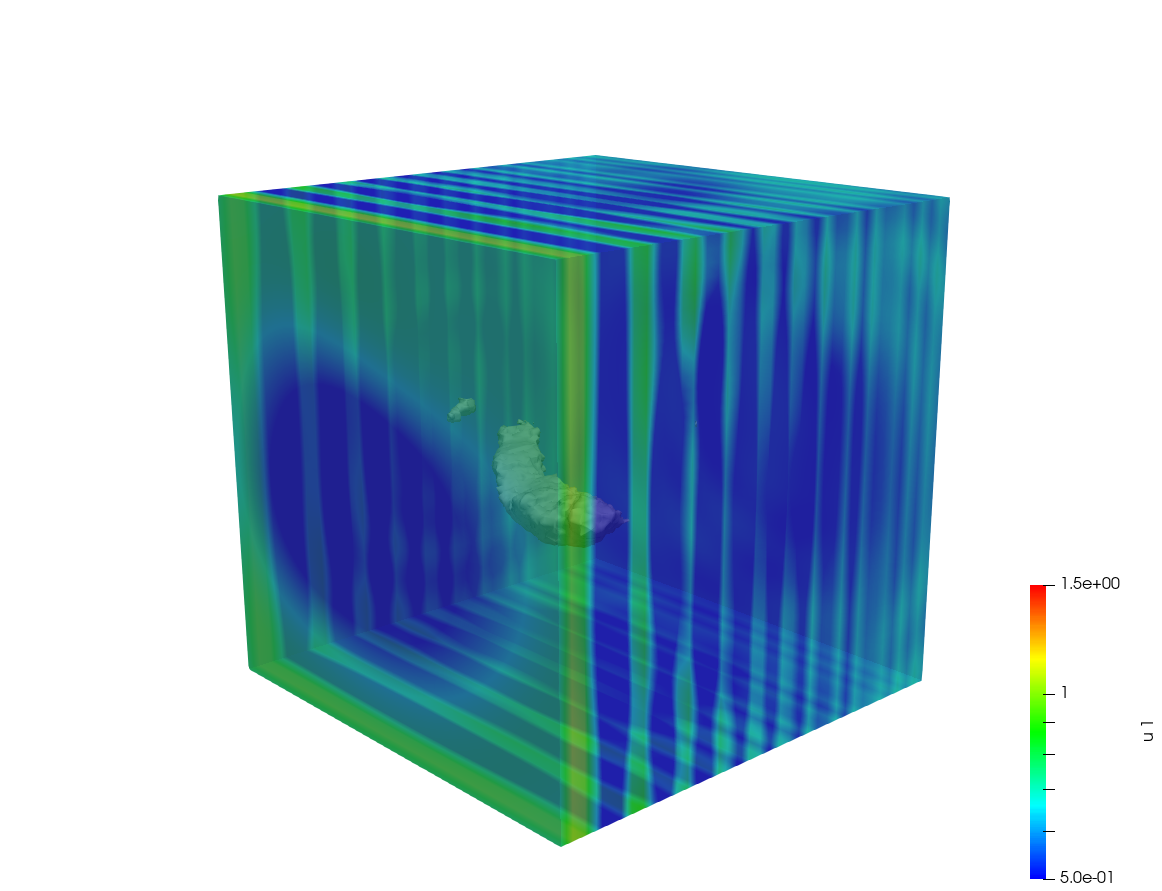}} \\
       a)  &  b)  &  c) \\
  {\includegraphics[scale=0.3, clip=true,]{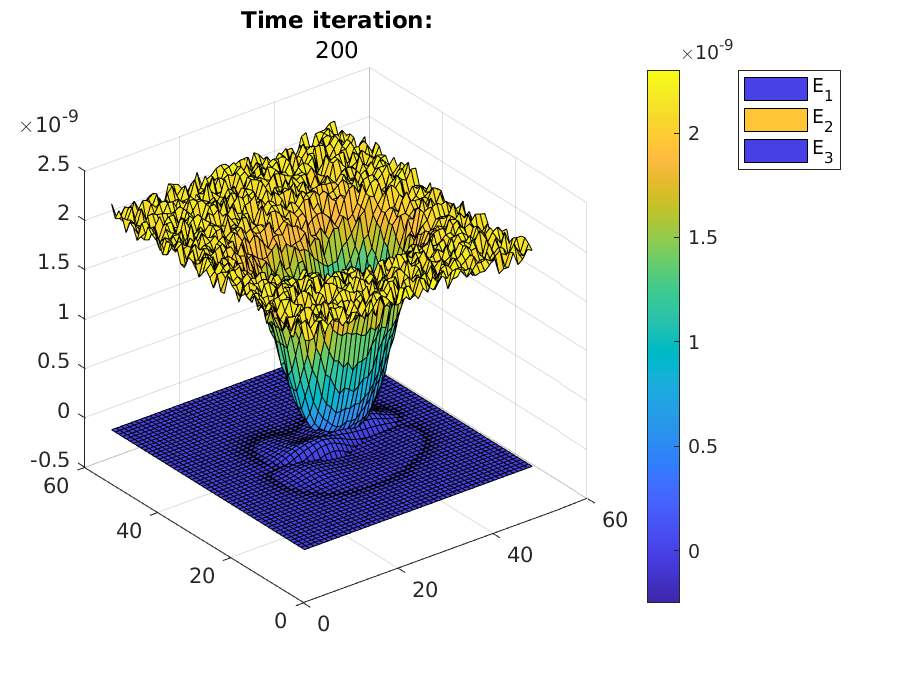}} &
 {\includegraphics[scale=0.3, clip=true,]{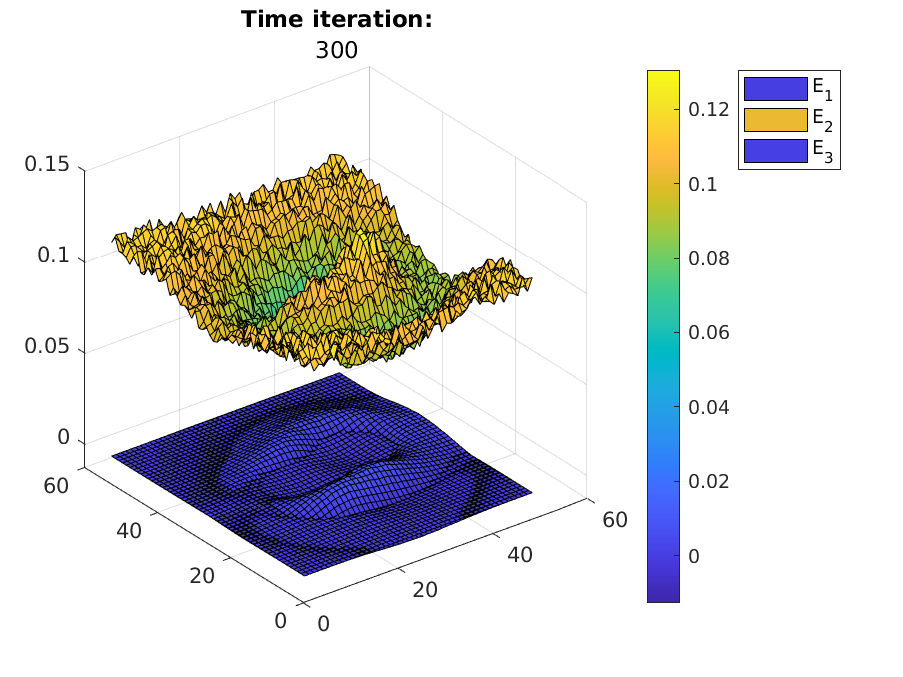}} &
  {\includegraphics[scale=0.3, clip=true,]{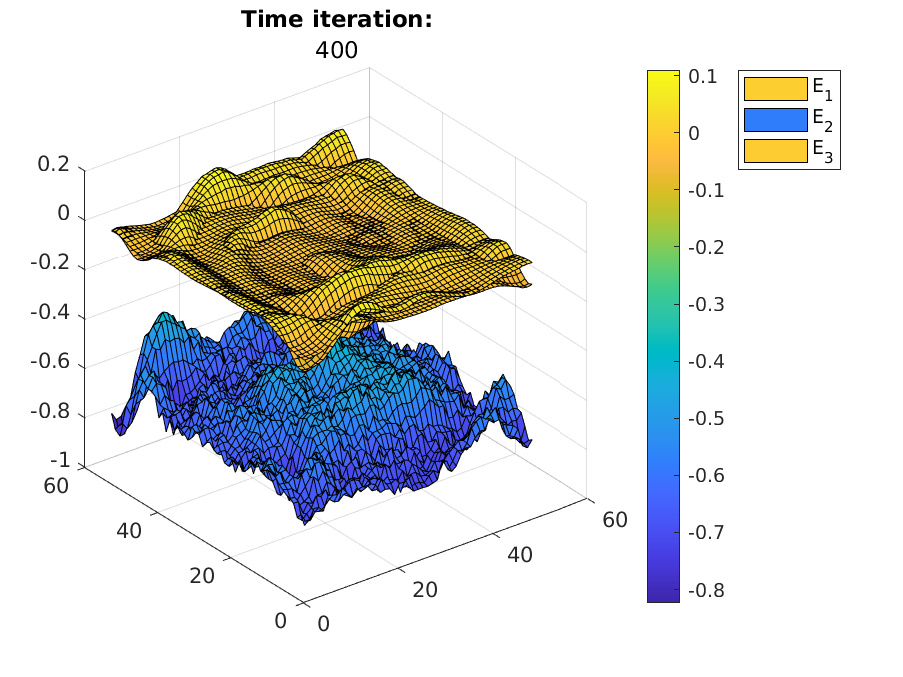}} \\
    d)  &  e)  &  f) \\
        {\includegraphics[scale=0.3, clip=true,]{omega40ref3noise10bary_EData_T200.png}} &
      {\includegraphics[scale=0.3,clip=true,]{omega40ref3noise10bary_EData_T300.png}}  &
      {\includegraphics[scale=0.3, clip=true,]{omega40ref3noise10bary_EData_T400.png}} \\
      g)  &  h)  &  i)
    \end{tabular}
  \end{center}
  \caption{\small\emph{ Test 2. a)-c): Solution $|E_h|$ of model problem \eqref{model1} at different times  for $\omega=40$ in \eqref{f}.
    d)-f): Transmitted noisy scattered data  $E_h = ({E_1}_h, {E_2}_h, {E_3}_h)$  of components of electric field
    $E=(E_1, E_2, E_3)$    at different times.   g)-i): Smoothed transmitted scattered data  $E_h = ({E_1}_h, {E_2}_h, {E_3}_h)$  of components of electric field
    $E=(E_1, E_2, E_3)$    at different times.  The noise level in data is $\delta=10\%$.  }}
  \label{fig:test2datanoise10}
\end{figure}

 The goal of our  numerical  tests Test 1, Test 2 was to reconstruct weighted dielectric permittivity function $\varepsilon_r$
 shown in Figures \ref{fig:exacteps}-a), b).  Figures \ref{fig:test1datanoise10}-a)-c), \ref{fig:test2datanoise10}-a)-c)  present simulated  solution $|E_h|$ in $\Omega_{FEM}$  of model problem \eqref{model1}  for Test 1 and Test 2, correspondingly.

\begin{figure}
     \centering
     \begin{tabular}{ccc}
        {\includegraphics[scale=0.07,  trim = 10cm 0.0cm 3.0cm 0.0cm, clip=true,]{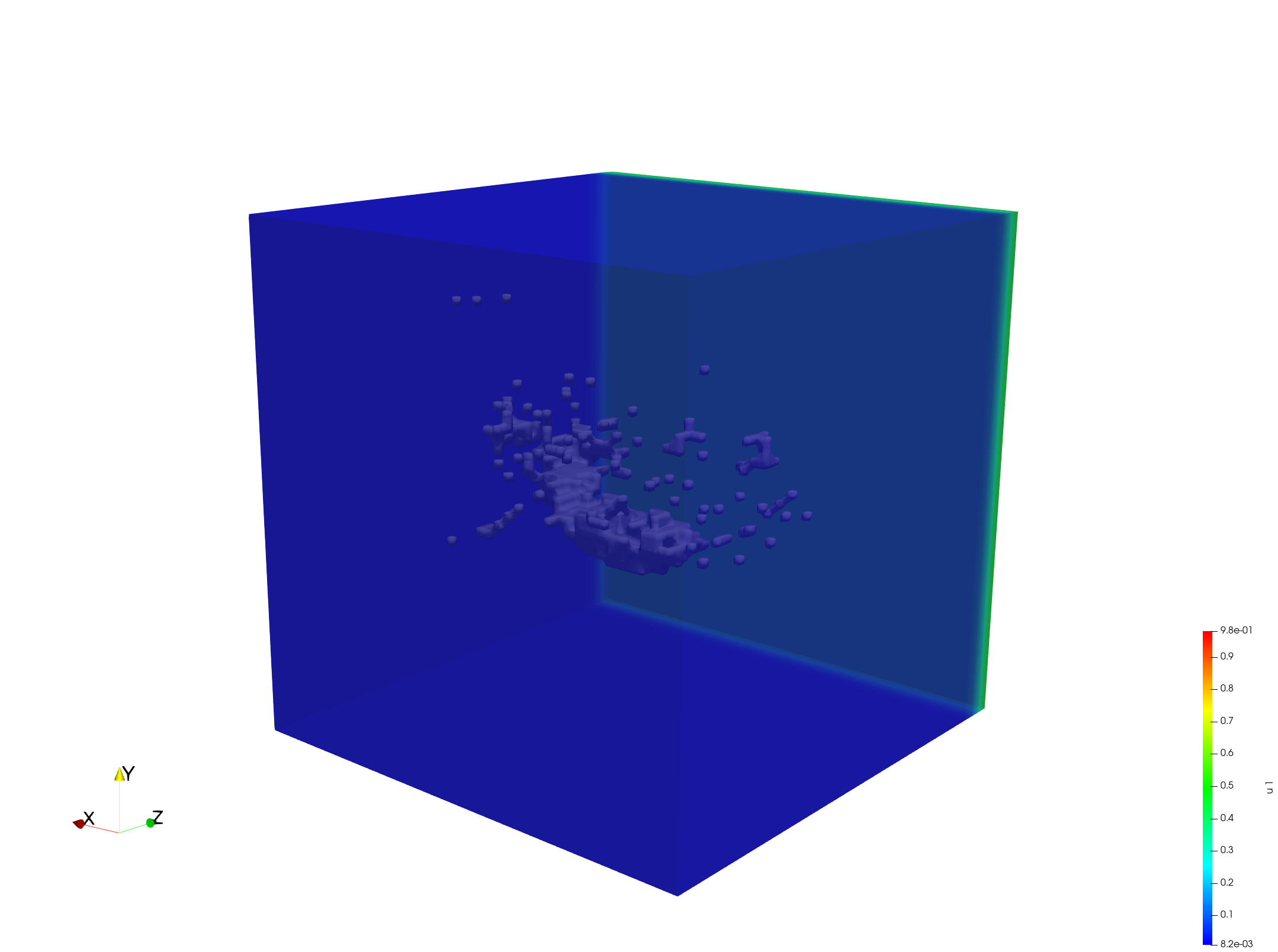}} &
        {\includegraphics[scale=0.07,  trim = 10cm 0.0cm 3.0cm 0.0cm, clip=true,]{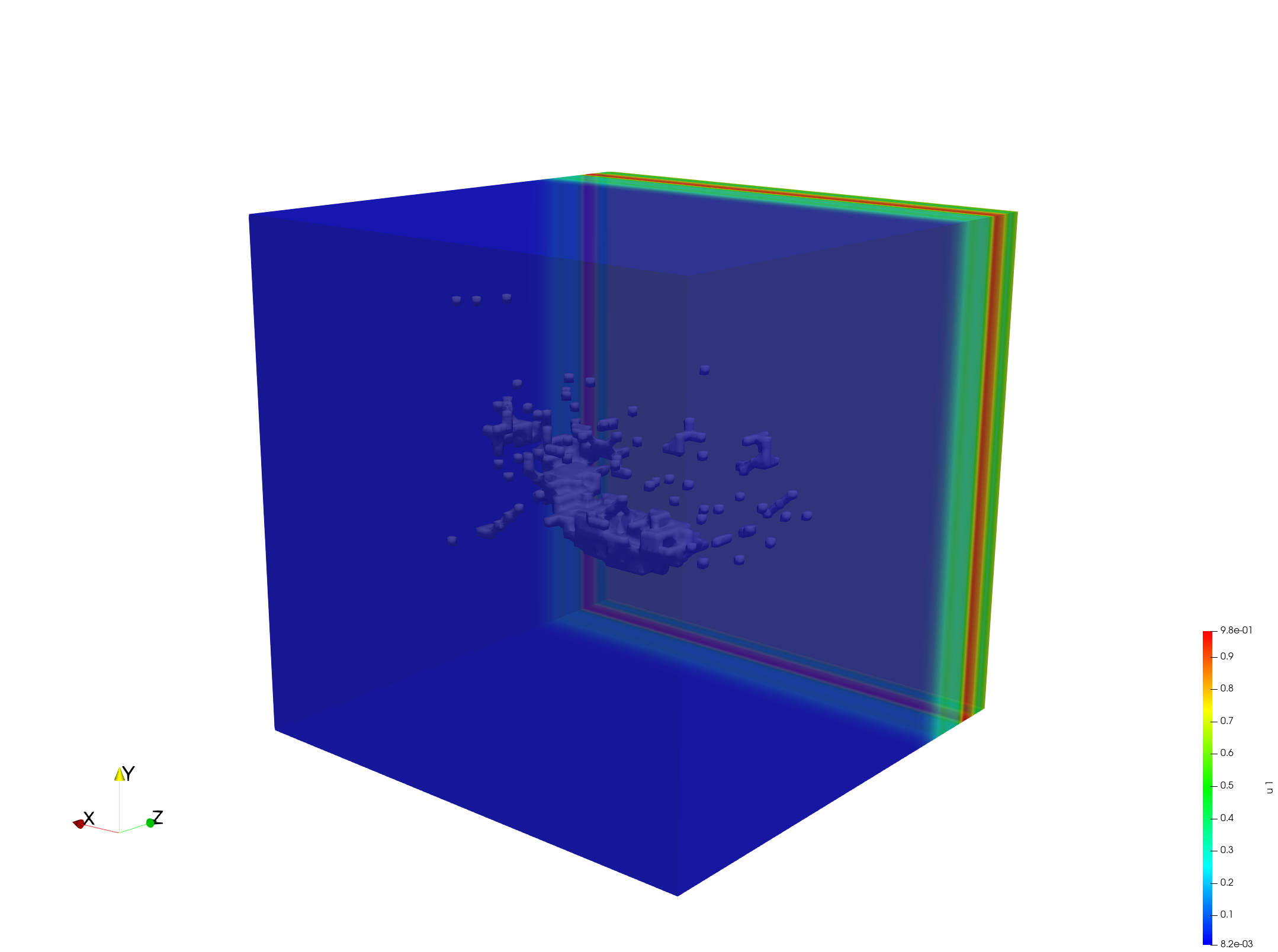}} &
        {\includegraphics[scale=0.07,  trim = 10cm 0.0cm 3.0cm 0.0cm, clip=true,]{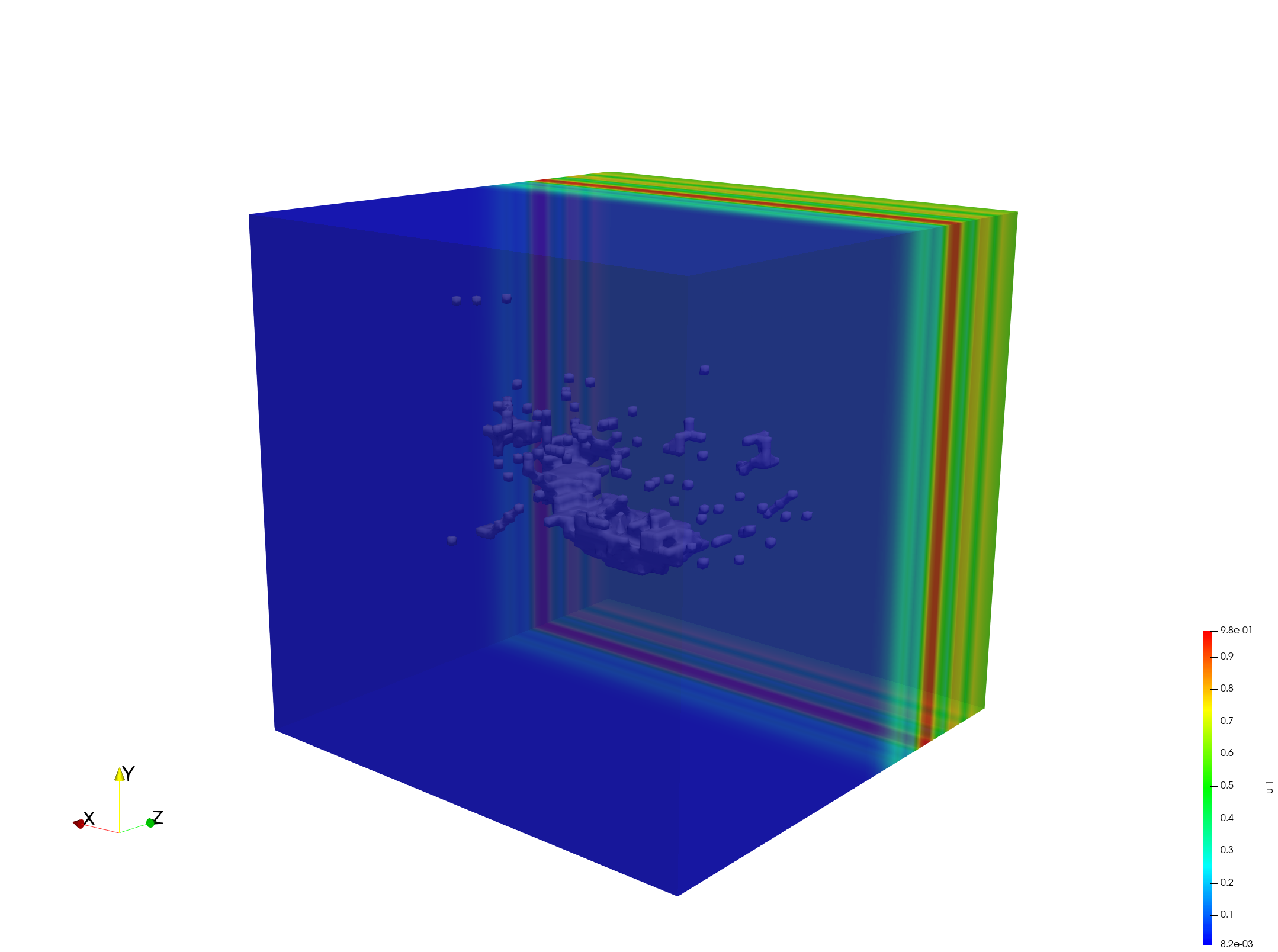}} \\
        a) $t = 0.24$ &  b) $t = 0.48$ &  c) $t = 0.72$ \\
        \hline \\
        {\includegraphics[scale=0.07,  trim = 10cm 0.0cm 3.0cm 0.0cm, clip=true,]{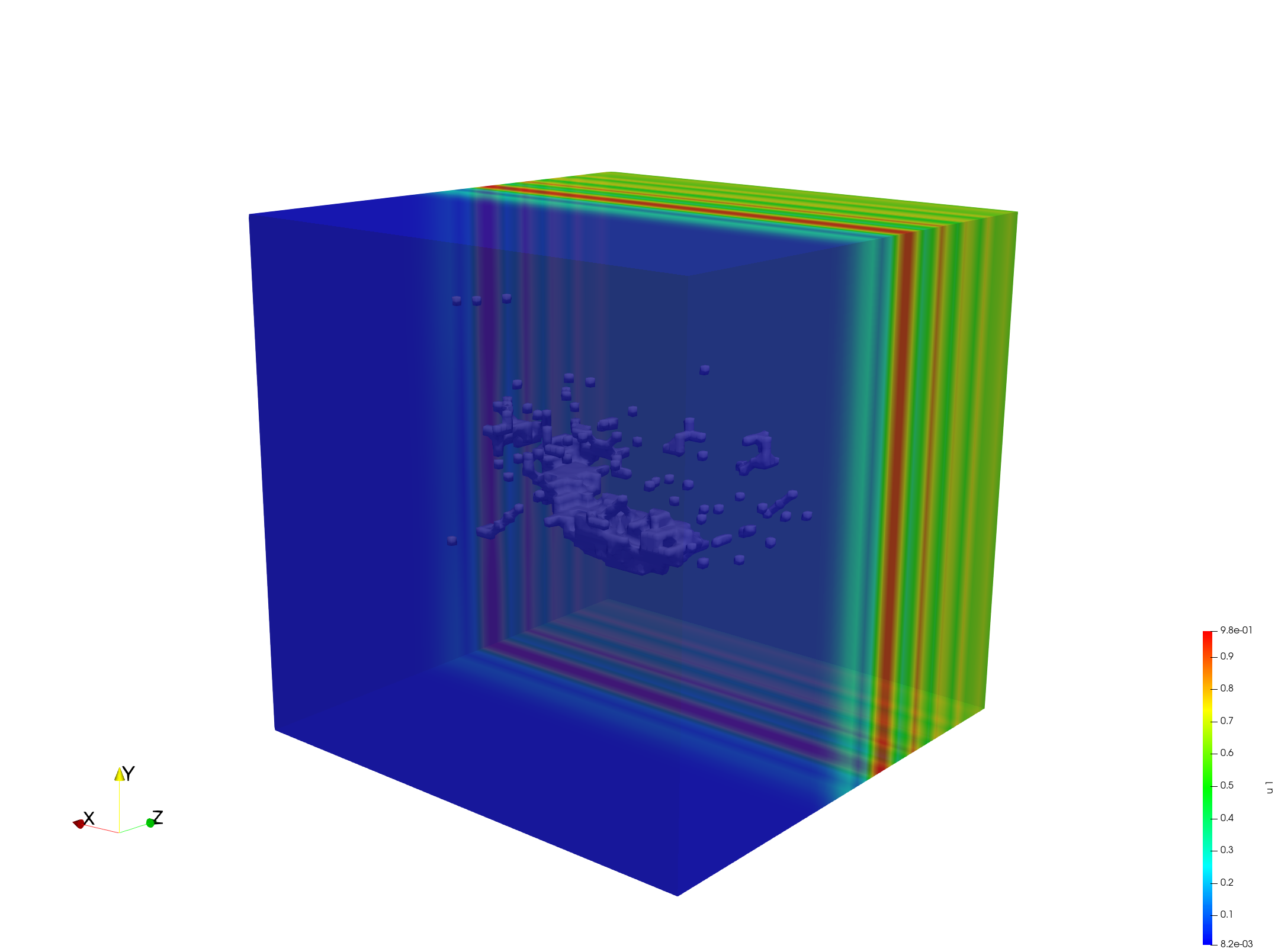}} &
        {\includegraphics[scale=0.07,  trim = 10cm 0.0cm 3.0cm 0.0cm, clip=true,]{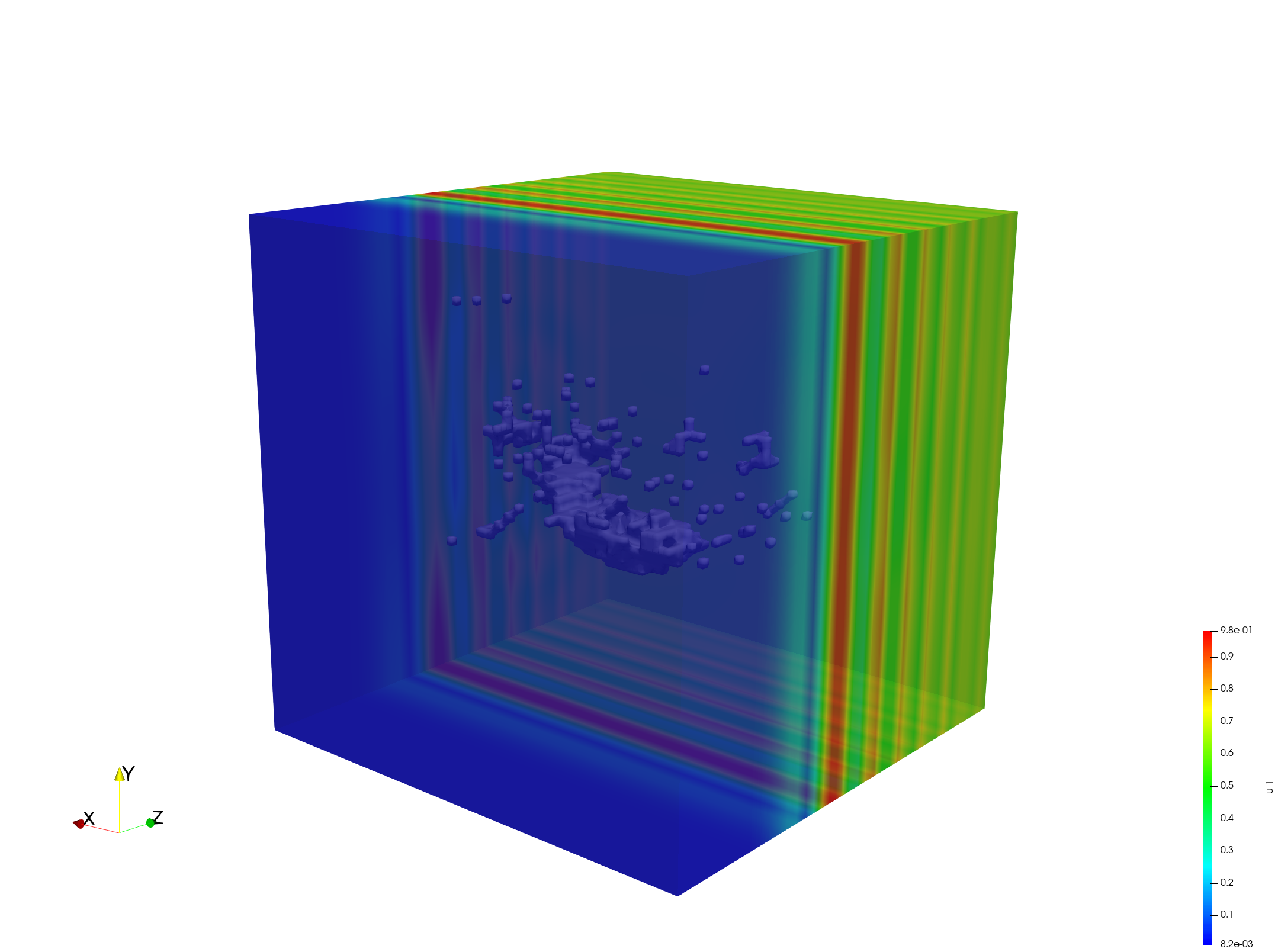}} &
        {\includegraphics[scale=0.07,  trim = 10cm 0.0cm 3.0cm 0.0cm, clip=true,]{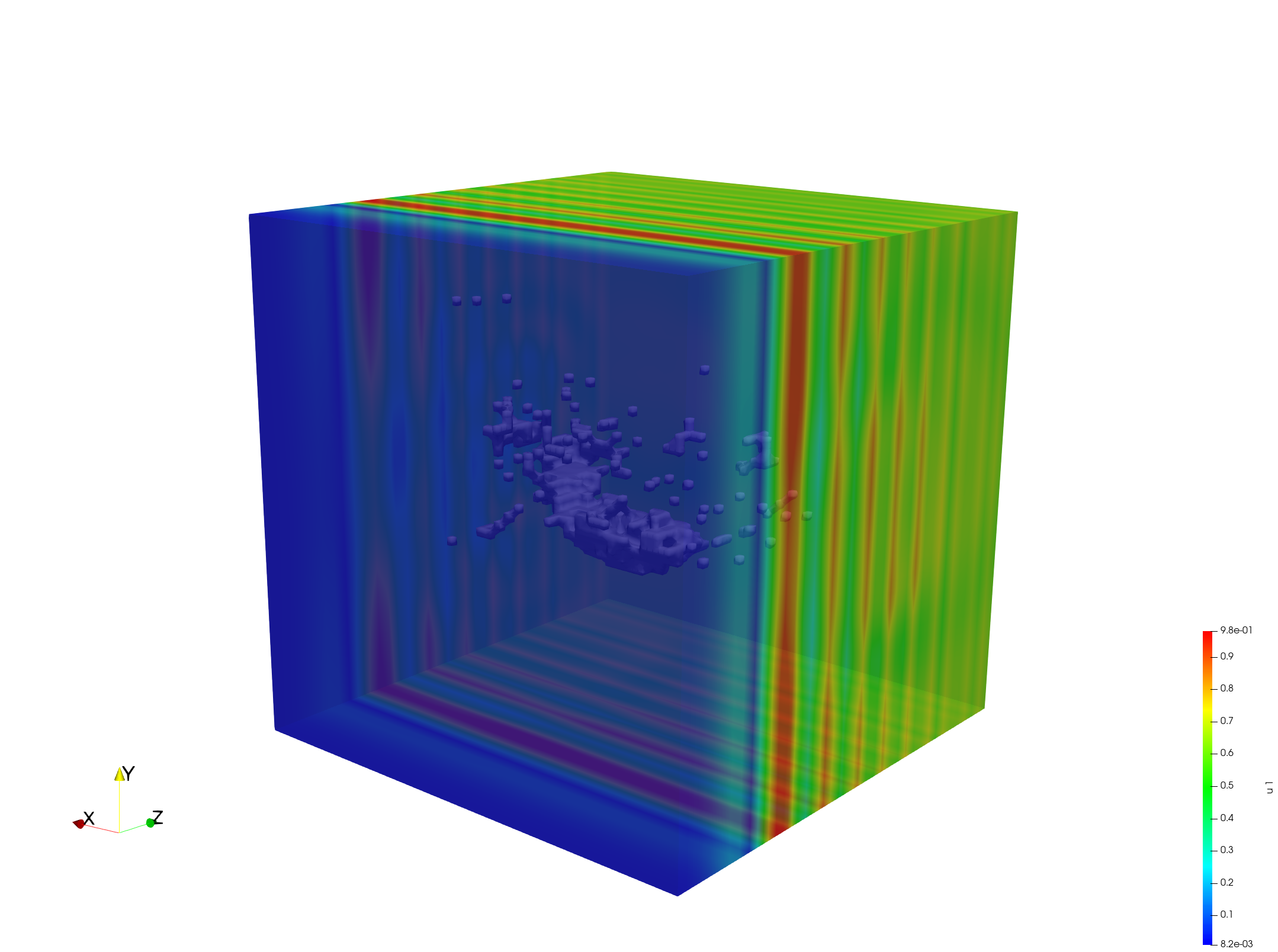}} \\
        d) $t = 0.96$ &  e) $t = 1.20$ &  f) $t = 1.44$\\
        \hline \\
        {\includegraphics[scale=0.07,  trim = 10cm 0.0cm 3.0cm 0.0cm, clip=true,]{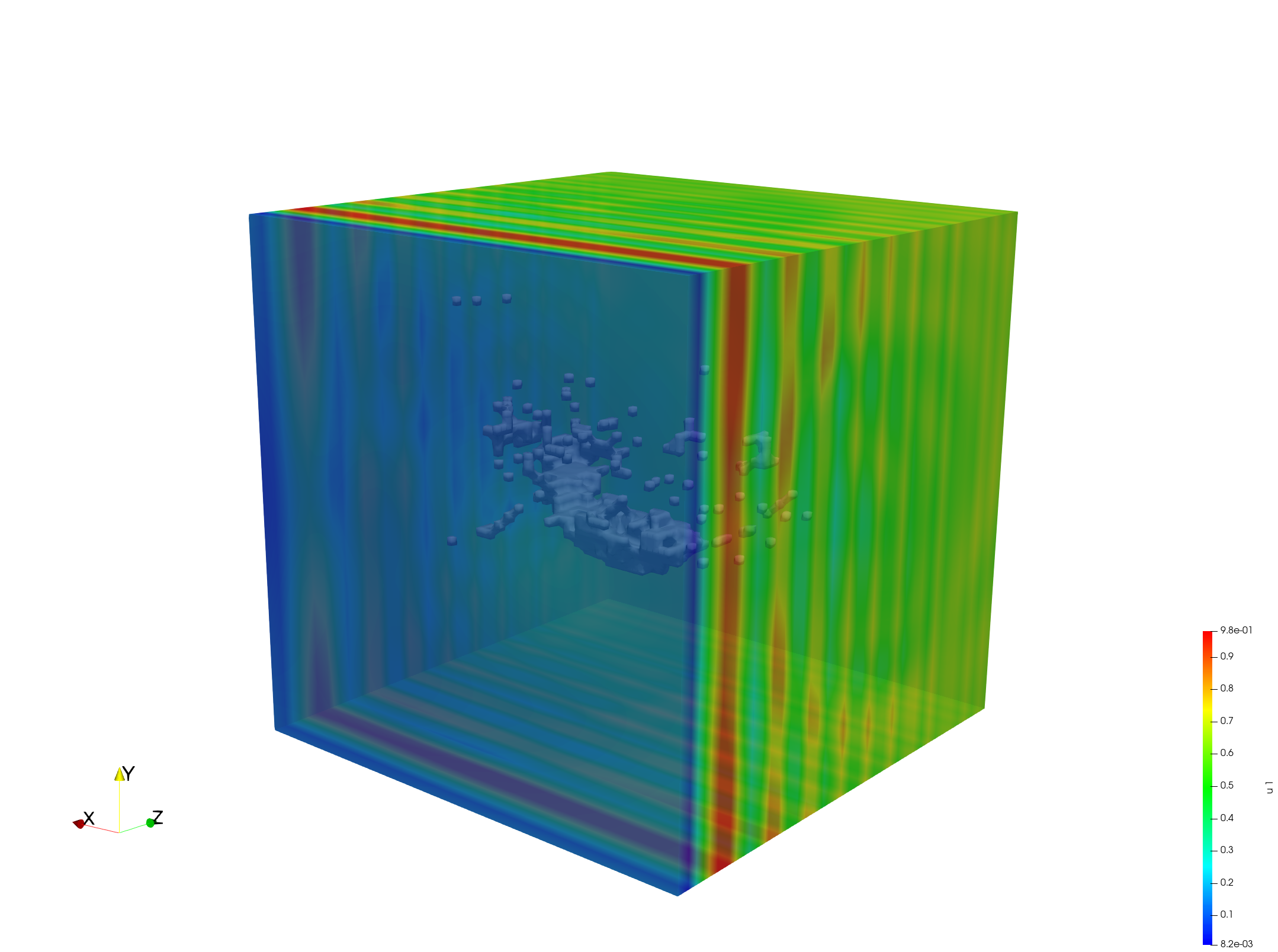}} &
        {\includegraphics[scale=0.07,  trim = 10cm 0.0cm 3.0cm 0.0cm, clip=true,]{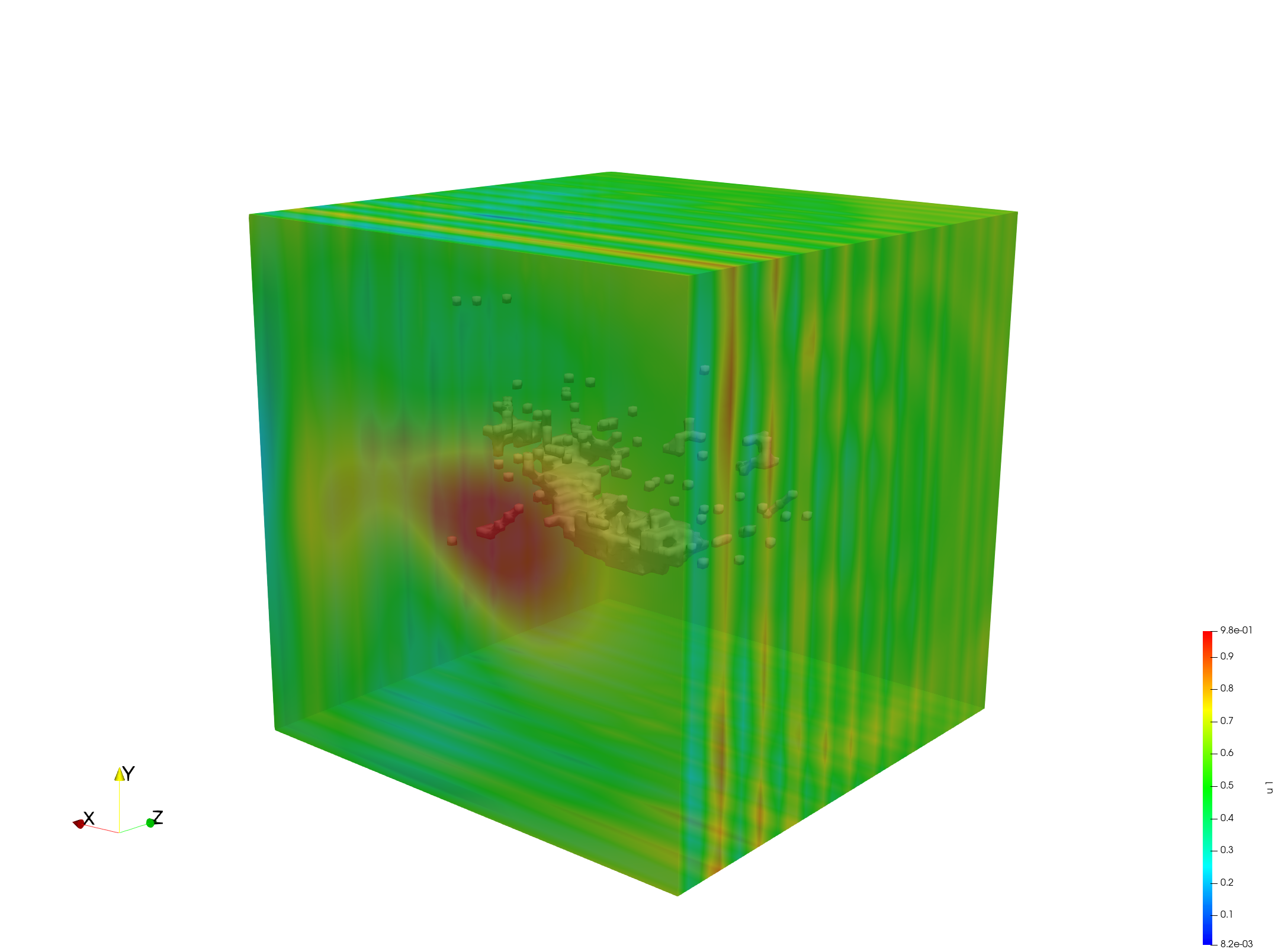}} &
        {\includegraphics[scale=0.07,  trim = 10cm 0.0cm 3.0cm 0.0cm, clip=true,]{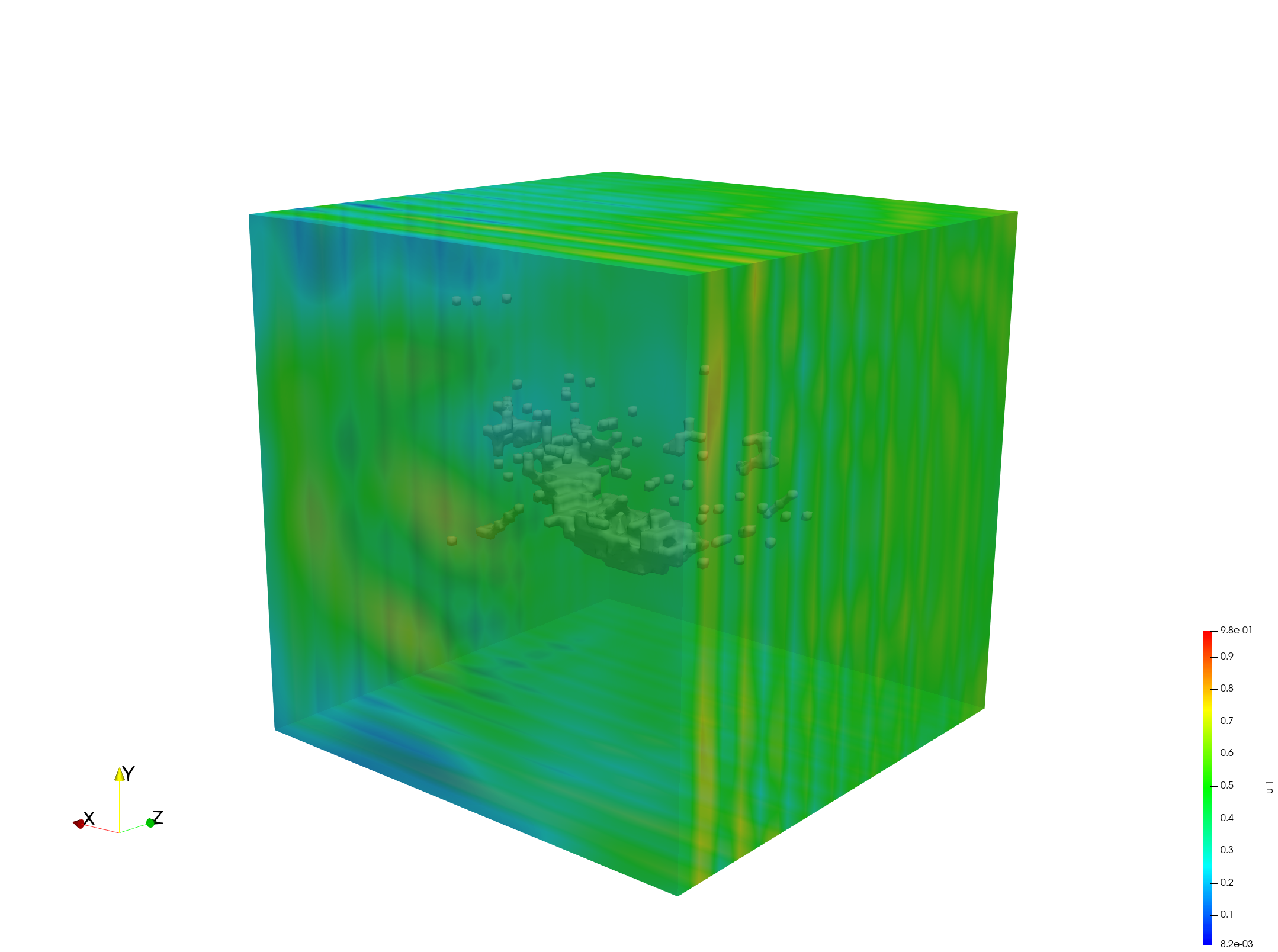}} \\
        g) $t = 1.68$ &  h) $t = 1.92$ &  i) $t = 2.16$\\
        \hline
     \end{tabular}
     \caption{\small \emph{
  The figures a)-i) illustrate how our planar wave $f$ used in implementations propagates through the medium by color plotting $|E_h|$ when $\omega = 40$, to clarify it's direction.
       }}
     \label{fig:wavepropex}
 \end{figure}

To perform computations for solution of inverse problem, we add
normally distributed Gaussian noise with mean $\mu = 0$ to simulated
electric field at the transmitted boundary $\partial_2 \Omega$. Then
we have smoothed out this data in order to get reasonable
reconstructions, see details of data-preprocessing in \cite{TBKF1, TBKF2}.
Computations of forward and inverse problems were done in time
$T=[0,3]$ with equidistant time step $\tau=0.006$ satisfying to CFL
condition.
Thus,  it took $500$ timesteps  at every iteration of
reconstruction Algorithm 4 to solve  forward or adjoint problem.
The time interval  $T=[0,3]$  was chosen  computationally such that the initialized plane wave could reach the transmitted boundary $\partial_2 \Omega$ in order to obtain meaningful reflections from the object inside the domain $\Omega_{\rm FEM}$. Figures  \ref{fig:wavepropex}-a)-i), \ref{fig:test1datanoise10}-a)-c), \ref{fig:test2datanoise10}-a)-c)   show these reflections in different tests.
Experimentally such signals can be produced by a Picosecond
Pulse Generator connected with a horn antenna, and scattered
time-dependent signals can be measured by a Tektronix real-time
oscilloscope, see \cite{TBKF1, TBKF2} for details of experimental
set-up for generation of a plane wave and collecting time-dependent
data. For  example, in our computational set-up, the  experimental
 time step between two   signals  can be $\tilde{\tau} = 6$
  picoseconds and every signal should be recorded   during $\tilde{T} = 3$ nanoseconds.

We have chosen  following set of admissible parameters for
reconstructed function $\varepsilon_r(x)$
 \begin{equation}\label{admpar}
   M_{\varepsilon_r} = \{ \rho\in C^2(\overline{\Omega })|1\leq \varepsilon_r(x)\leq 10\},
 \end{equation}
as well as tolerance $\theta =
10^{-5}$ at step 3 of the conjugate gradient Algorithm 3.
Parameters $\beta_k$ in the refined procedure of Algorithm 4 was chosen as
 the constant $\beta_k =
0.8$  for all refined meshes ${K_h}_k$.

Figures \ref{fig:test1datanoise10}-d)-i) -
\ref{fig:test2datanoise10}-d)-i) show simulated data of model problem
\eqref{model1} for all components $(E_1,E_2,E_3)(x,t)$ of electric
field $E(x,t)$ at different times at the transmitted boundary
$\partial_2 \Omega$.  Figures \ref{fig:test1datanoise10}-d)-f) -
\ref{fig:test2datanoise10}-d)-f) show randomly distributed noisy data
and Figures \ref{fig:test1datanoise10}-g)-i) -
\ref{fig:test2datanoise10}-g)-i) show smoothed noisy data used for
solution of inverse problem.

These figures show that largest by
amplitude reflections, or transmitted data, are obtained from the
second component $E_2$ of the electric field $E$.
 The same observation is obtained in previous
works \cite{BMaxwell, BondestaB} where was used a similar
computational set-up with a plane wave. However, comparison of all
three components was not presented in \cite{BondestaB}.  Domination of
reflections  at the transmitted boundary from the $E_2$ component can be explained by the fact that
we initialize only one component of the electric field $E=(E_1, E_2,
E_3)$ as a plane wave $f(t) = (0,f_2,0)(t)$ at $\Gamma_{1,1}$ in the
model problem \eqref{model1}, and thus, two other components $E_1,
E_3$ will be smaller  by  amplitude than the $E_2$  when we  use the explicit
scheme \eqref{femod1new} for computations. See also theoretical
justification of this fact in \cite{RK}.

Numerical tests of \cite{BondestaB} show that the best
reconstruction results of the space-dependent function $\varepsilon_r(x)$
for $\sigma=0$ in $\Omega$ are obtained for $\omega = 40$ in (\ref{f}).
Thus, we  performed simulations of the forward problem \eqref{model1} taking
 $\sigma=0$
for different $\omega = 40, 60, 80, 100$ in \eqref{f}. It turned out that
for chosen computational set-up with final time $T=3$ 
maximal   values of   scattered function $E_2$ are obtained for $\omega=40$.
Thus,
we take $\omega = 40$ in (\ref{f}) in all our tests.

 We assume that both functions $\varepsilon_r, \sigma$ satisfy
 conditions \eqref{2.3}: they are known inside $\Omega_{\rm out} \cup
 \Omega_{\rm FDM}$ and unknown inside $\Omega_{IN}$.
 The goal of our numerical tests is to
 reconstruct the function $\varepsilon_r$ of the domain $\Omega_{\rm FEM}$
 of Figure \ref{fig:exacteps} under conditions \eqref{2.3} and the additional condition that the function
 $\sigma(x)$ of this domain is known. See Table \ref{tab:table1} for distribution of $\varepsilon_r, \sigma$ in
 $\Omega_{\rm FEM}$.

 The computational set-up for solution of inverse problem is as
 follows.  We generate transmitted data by solving the model problem
 \eqref{model1} on three times adaptively refined mesh. In this way we
 avoid variational crimes when we solve the inverse problem. The
 transmitted data is collected at receivers located at every point of
 the transmitted boundary $\partial_2 \Omega$, and then normally distributed Gaussian noise
 $\delta = 3 \%, 10\%$ with mean $\mu=0$ is added to this data,
 see Figures \ref{fig:test1datanoise10}-d)-f) - \ref{fig:test2datanoise10}-d)-f).
The next step is data pre-processing: the noisy data is
 smoothed out,
 see Figures \ref{fig:test1datanoise10}-g)-i) - \ref{fig:test2datanoise10}-g)-i).
  Next, to reconstruct
  $\varepsilon_r$  we minimize the
  Tikhonov functional \eqref{functional}.  For solution of the
  minimization problem we introduce Lagrangian and search for a
  stationary point of it  using an adaptive Algorithm 4, see details in section \ref{sec:ref}.

  We take the initial approximation
  $\varepsilon_0 = 1$ at all points of the computational domain
  what corresponds to starting of our computations from the homogeneous
  domain. This is done because of  previous
computational works \cite{BondestaB} as well as
experimental works of \cite{meaney1, ieee12, ieee19} where was shown that a such
choice gives good results of reconstruction of dielectric permittivity function.

\begin{table}[tbp] 
\vspace{2mm}
\centerline{
  \begin{tabular}{|c|c|}
     \hline
     \multicolumn{2}{|c|}
                 {Test 1}
                 \\
 \hline
   $\delta=3\%$ &  $\delta = 10\%$ 
 \\
\begin{tabular}{l|l|l|l} \hline
Mesh & $\max_{\Omega_{\rm FEM}} {\varepsilon_h}_k$ & $\frac{\max_{\Omega_{\rm FEM}} |\varepsilon_r - {\varepsilon_h}_k | }{\max_{\Omega_{\rm FEM}} |\varepsilon_r|}$ & $M^k$    \\ \hline
${K_h}_0$ & 6.535  & 0.274 & 2   \\
${K_h}_1$ & 7.865 &  0.126 &   2 \\
${K_h}_2$ & 10.0  &  0.111 &  2 \\
\end{tabular}
 & 
\begin{tabular}{l|l|l|l} \hline
  Mesh & $\max_{\Omega_{\rm FEM}} {\varepsilon_h}_k$ &    
$\frac{\max_{\Omega_{\rm FEM}} |\varepsilon_r - {\varepsilon_h}_k | }{\max_{\Omega_{\rm FEM}} |\varepsilon_r|}$ & $M^k$  \\ \hline
 ${K_h}_0$  & 7.019  & 0.220   &  2    \\
${K_h}_1$   & 7.481  & 0.167  &   4 \\
${K_h}_2$  & 9.234  &  0.026 & 4   \\
\end{tabular} \\
\hline
\end{tabular}}
\caption{\textit{ Test 1.
Computational results of the
    reconstructions $\max_{\Omega_{\rm FEM}}  {\varepsilon_h}_k$ on a coarse and on adaptively refined meshes together with  relative 
    errors computed in the maximal contrast of $ \max_{\Omega_{\rm FEM}}  \varepsilon_r,\max_{\Omega_{\rm FEM}}  {\varepsilon_h}_k$. Here,
    $\max_{\Omega_{\rm FEM}}  {\varepsilon_h}_k$ denotes  maximum of the computed function $\varepsilon_h$
    on $k$ times refined mesh ${K_h}_k$ in the domain $\Omega_{\rm FEM}$, and
    $M^k$ denotes the final number of
    iterations in the conjugate gradient Algorithm 3 on $k$ times refined mesh ${K_h}_k$  for reconstructed
    function ${\varepsilon_h}_k, k=0,1,2$. 
}}
    \label{tab:test1}
\end{table}


\begin{table}[tbp] 
\vspace{2mm}
\centerline{
  \begin{tabular}{|c|c|}
     \hline
     \multicolumn{2}{|c|}
                 {Test 1, Computational Time}
                 \\
 \hline
   $\delta=3\%$ &  $\delta = 10\%$ 
 \\
\begin{tabular}{l|l|l|l|l} \hline
Mesh & nno &  Time (sec)& Rel. time  & $M^k$    \\ \hline
${K_h}_0$ & 63492  & 1183  & 3.73 $\cdot 10^{-5}$ & 2   \\
${K_h}_1$ &  64206  & 1199   & 3.74 $\cdot  10^{-5}$ & 2 \\
${K_h}_2$ &  65284 & 1212   & 3.71 $\cdot  10^{-5}$ &  2 \\
\end{tabular}
 & 
\begin{tabular}{l|l|l|l|l} \hline
  Mesh &  nno  &     Time (sec) & Rel. time  & $M^k$  \\ \hline
 ${K_h}_0$  &  63492  & 1180  &  $3.71\cdot 10^{-5}$    &  2    \\
${K_h}_1$   &  64766  & 2415  &  $7.43 \cdot 10^{-5}$     &  4 \\
${K_h}_2$   &  67965  & 2525  & $7.435 \cdot 10^{-5}$  &  4   \\
\end{tabular} \\
\hline
\end{tabular}}
\caption{
  \textit{ Test 1.  Performance of the reconstruction Algorithm 4
    (in seconds)
  on adaptively refined meshes. Here,
  $k$ is number of the refined mesh ${K_h}_k$ of the domain $\Omega_{\rm FEM}$,
 ${\rm nno}$ is number of the nodes in the computational   mesh  ${K_h}_k$,
  and
    $M^k$ denotes the final number of
  iterations in the conjugate gradient Algorithm 3.
}}
    \label{tab:test1time}
\end{table}


\begin{figure}[tbp]
  \begin{center}
    \begin{tabular}{ccc}
     {\includegraphics[scale=0.20, trim = 8.0cm 0.0cm 8.0cm 0.0cm, clip=true,]{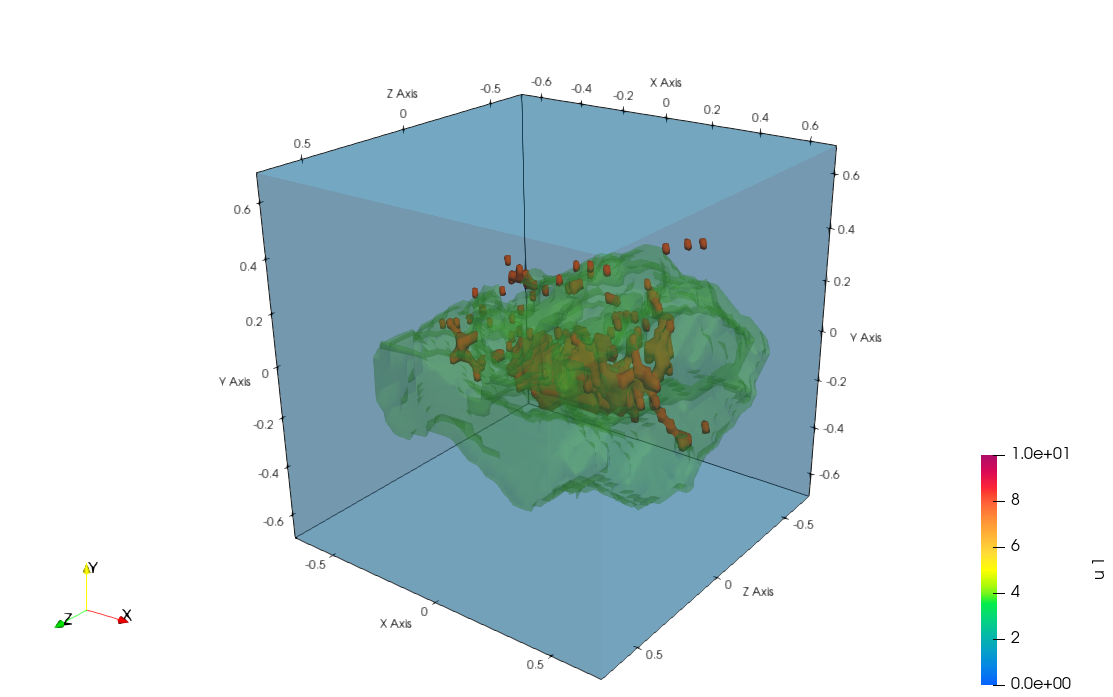}}
     &
      {\includegraphics[scale=0.20, trim = 8.0cm 0.0cm 8.0cm 0.0cm, clip=true,]{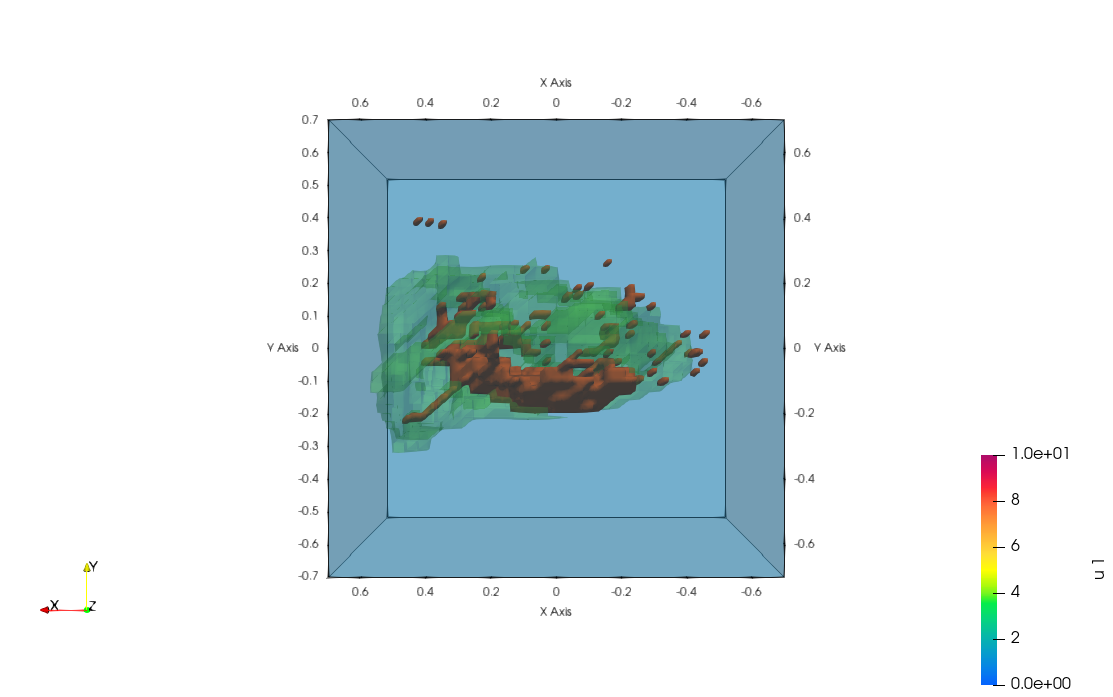}} &
      {\includegraphics[scale=0.20, trim = 8.0cm 0.0cm 8.0cm 0.0cm, clip=true,]{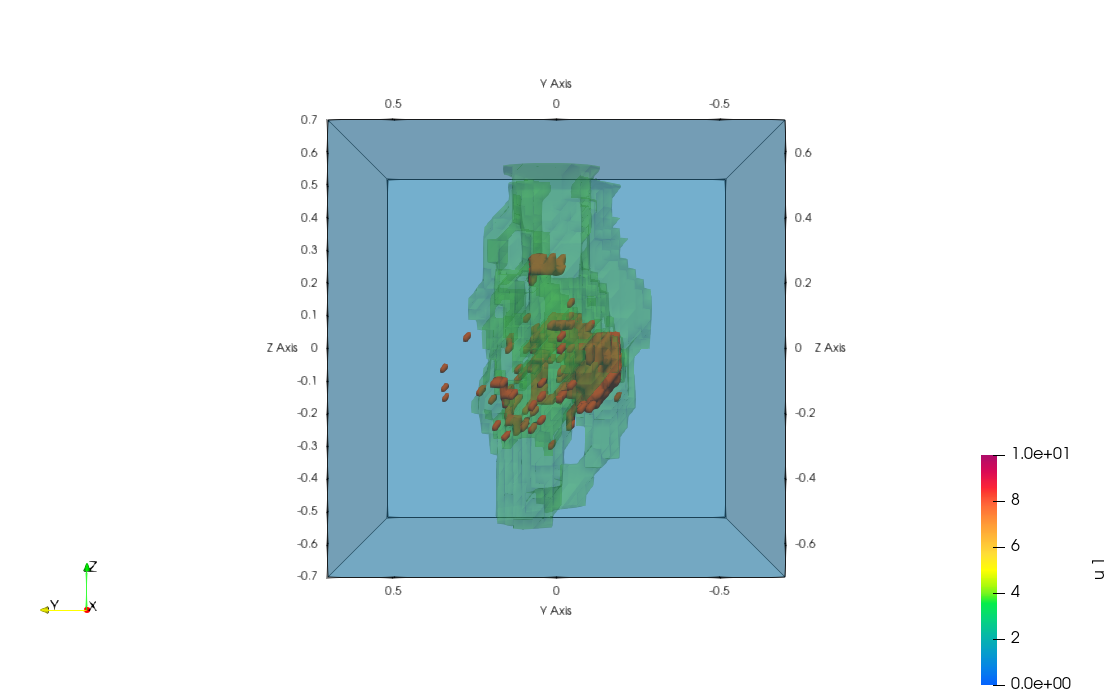}}\\
    a)  perspective view  & b) $x_1 x_2$ view & c) $x_2 x_3$ view \\
      \\
      \hline
          {\includegraphics[scale=0.20, trim = 8.0cm 0.0cm 8.0cm 0.0cm, clip=true,]{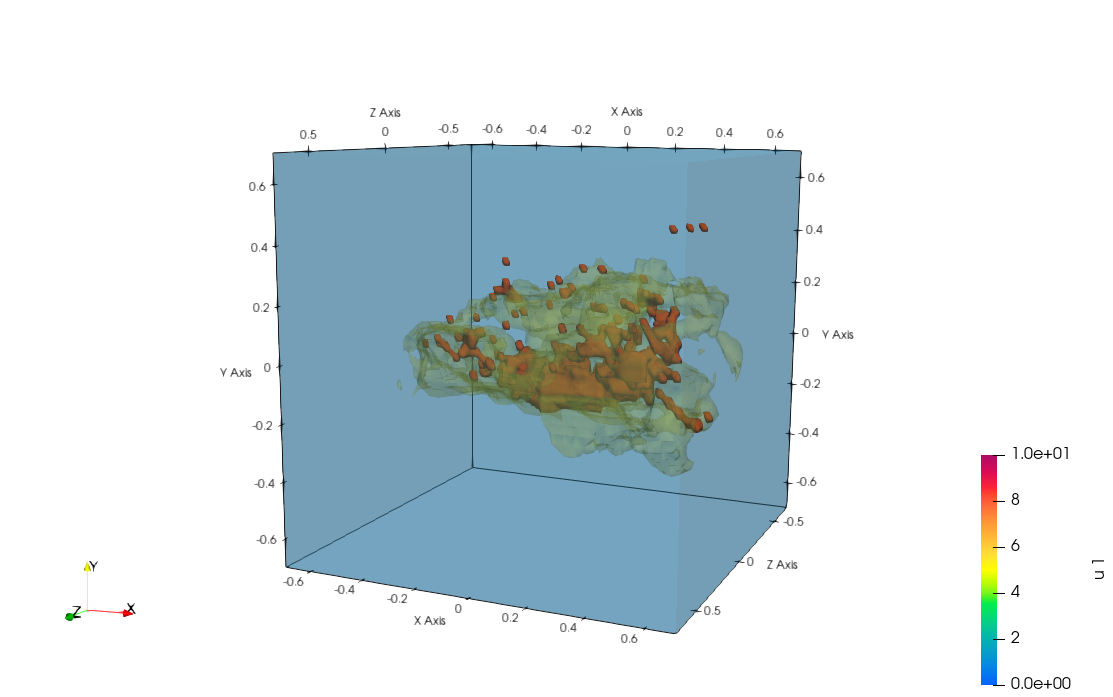}}
          &
      {\includegraphics[scale=0.20, trim = 8.0cm 0.0cm 8.0cm 0.0cm, clip=true,]{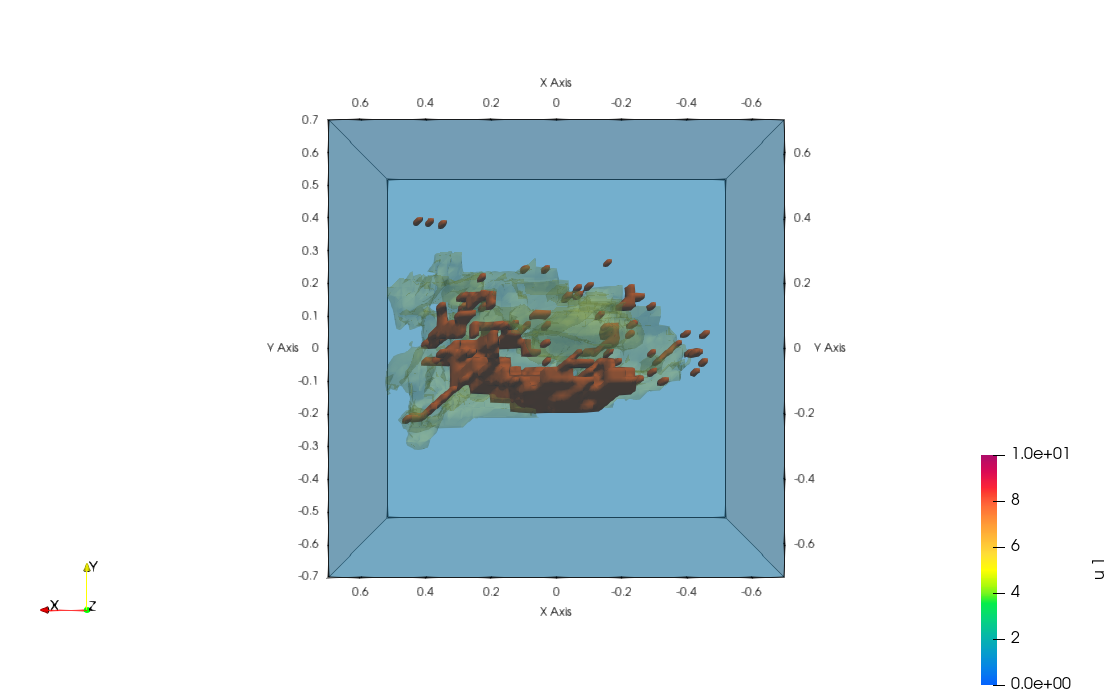}} &
      {\includegraphics[scale=0.20, trim = 8.0cm 0.0cm 8.0cm 0.0cm, clip=true,]{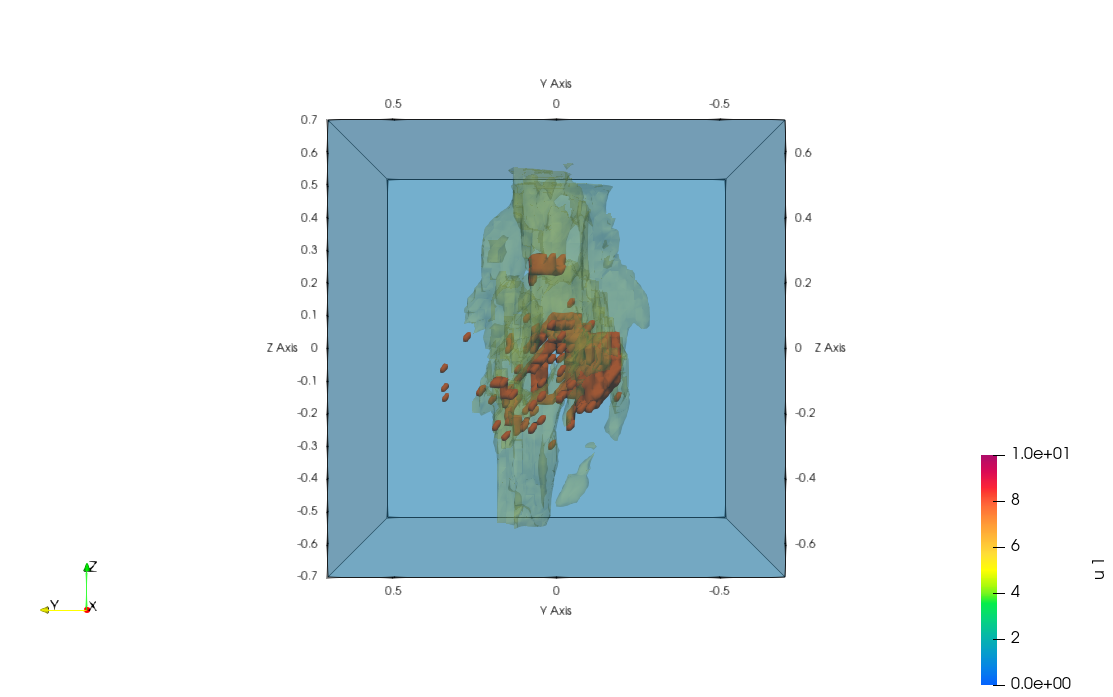}} \\
   d)   perspective view  & e) $x_1 x_2$ view & f) $x_2 x_3$ view \\
      \\
      \hline
      \\
      {\includegraphics[scale=0.15, trim = 5.0cm 0.0cm 6.0cm 0.0cm, clip=true,]{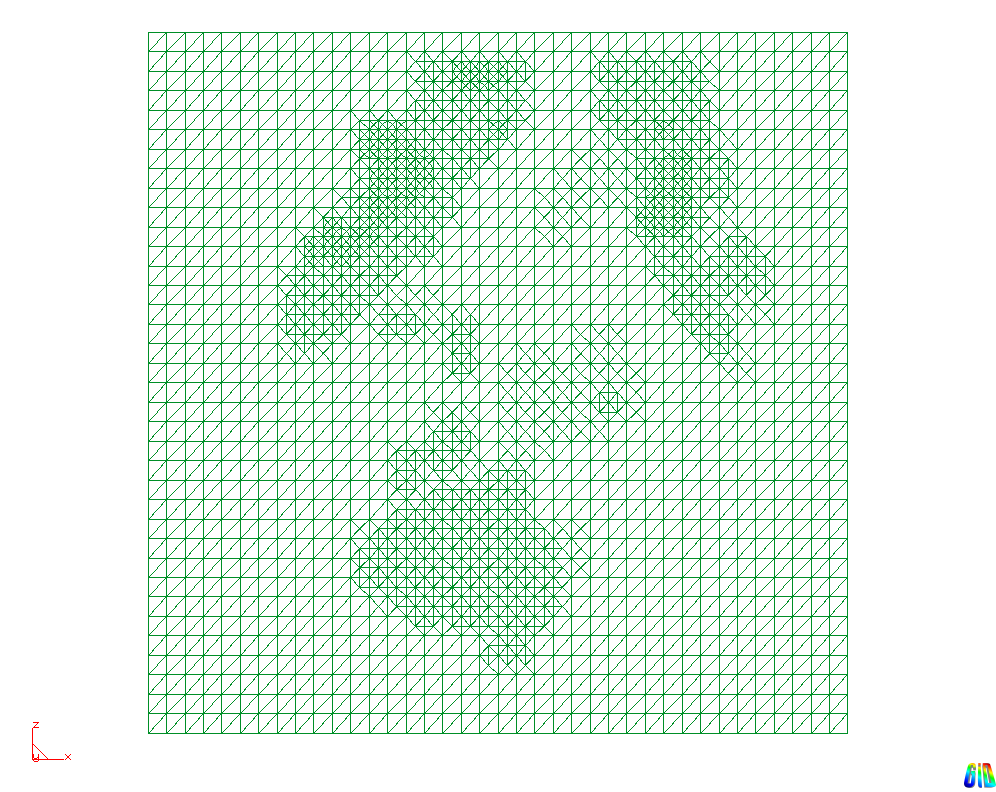}}
          &
      {\includegraphics[scale=0.15, trim = 5.0cm 0.0cm 6.0cm 0.0cm, clip=true,]{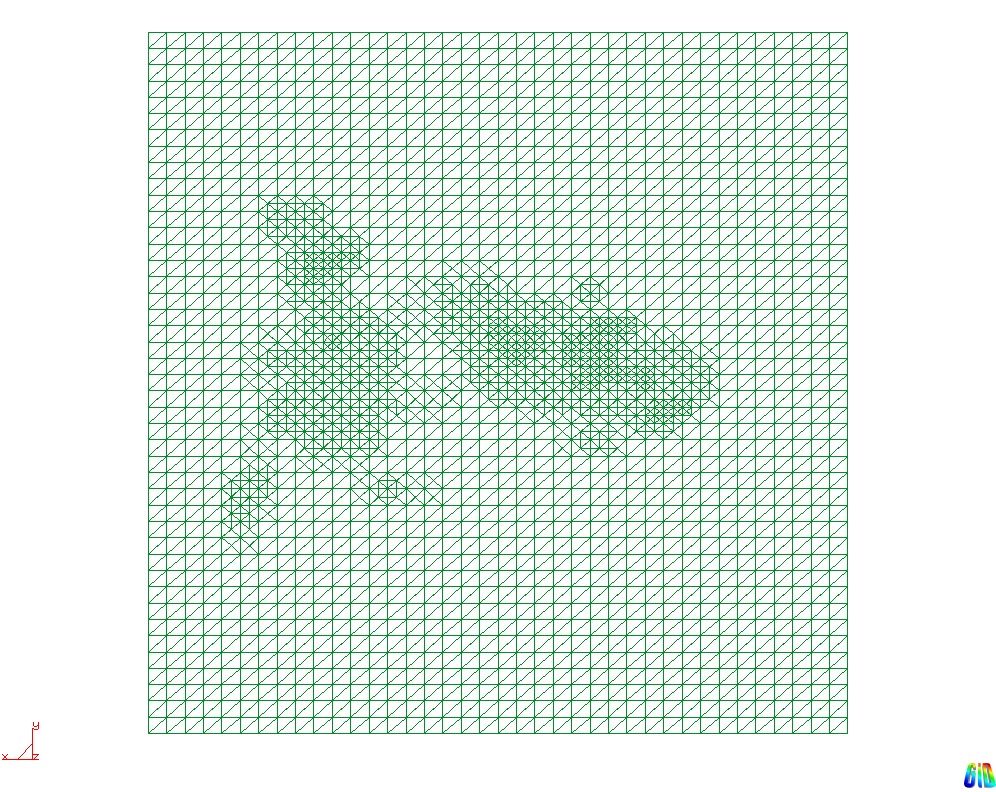}} &
      {\includegraphics[scale=0.15, trim = 5.0cm 0.0cm 6.0cm 0.0cm, clip=true,]{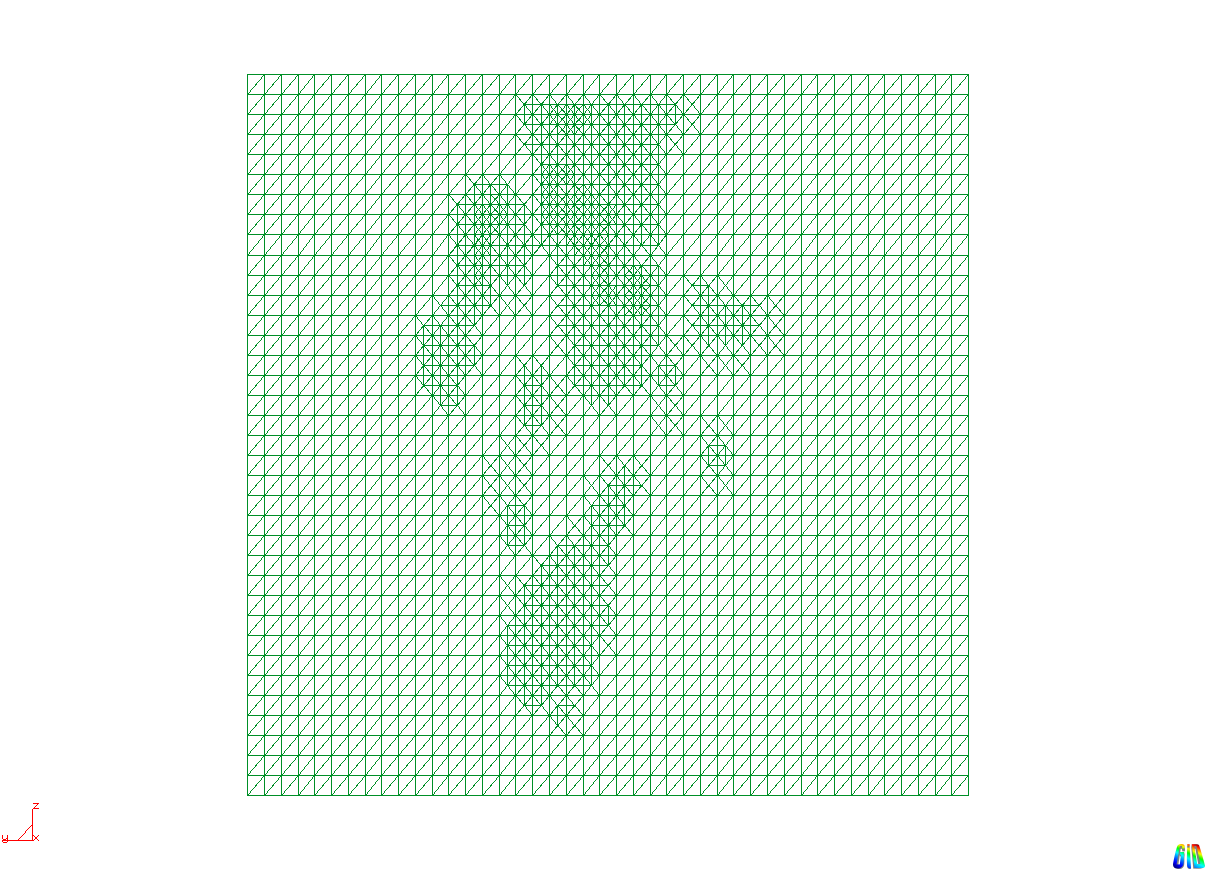}} 
      \\
    g)  $x_1 x_3$ view & h) $x_1 x_2$ view & i) $x_2 x_3$ view \\
       \hline
    \end{tabular}
  \end{center}
  \caption{
    \small\emph{Test 1. a)-c): Reconstructions
    $\varepsilon_{h_0} \approx 5$ (outlined in transparent green color) of
    $\varepsilon_r$ obtained on the coarse mesh.  d)-f): Reconstructions ${\varepsilon_h}_2 \approx 5$
    obtained on refined mesh ${K_h}_2$.  g) - i): refined mesh ${K_h}_2$.
    The noise level in the data
    is $\delta= 10\%$.
  See Table \ref{tab:test1}  for obtained contrasts  $\max_{\Omega_{\rm FEM}}  {\varepsilon_h}_k, k=0,1,2$.
    For comparison we also present exact isosurface  with value corresponding to reconstructed one
    and outlined by red color.
}}
  \label{fig:test1noise10v2}
\end{figure}

\begin{figure}[tbp]\vspace*{-3cm}
  \begin{center}
    \begin{tabular}{ccc}
      \\
      {\includegraphics[scale=0.16, trim = 0.0cm 0.0cm 0cm 0.0cm, clip=true,]{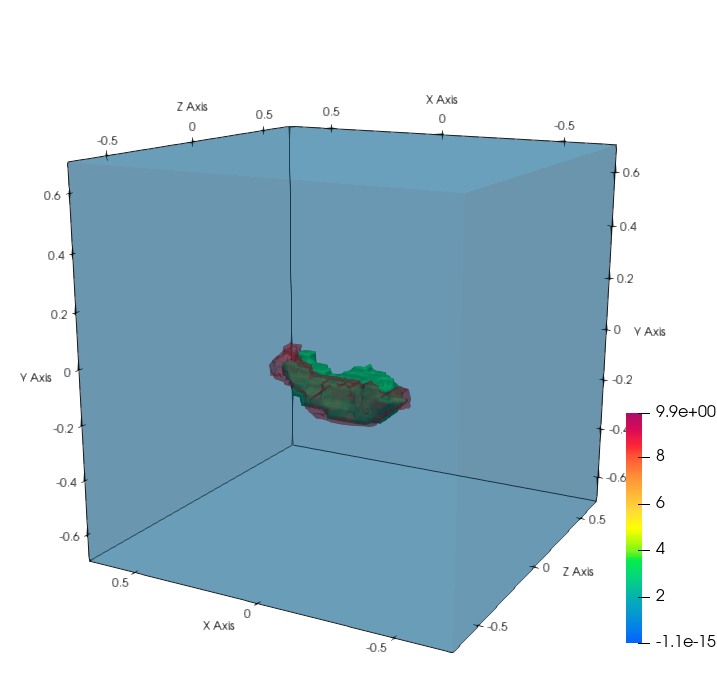}}   &
      {\includegraphics[scale=0.16, trim = 0.0cm 0.0cm 0cm 0.0cm, clip=true,]{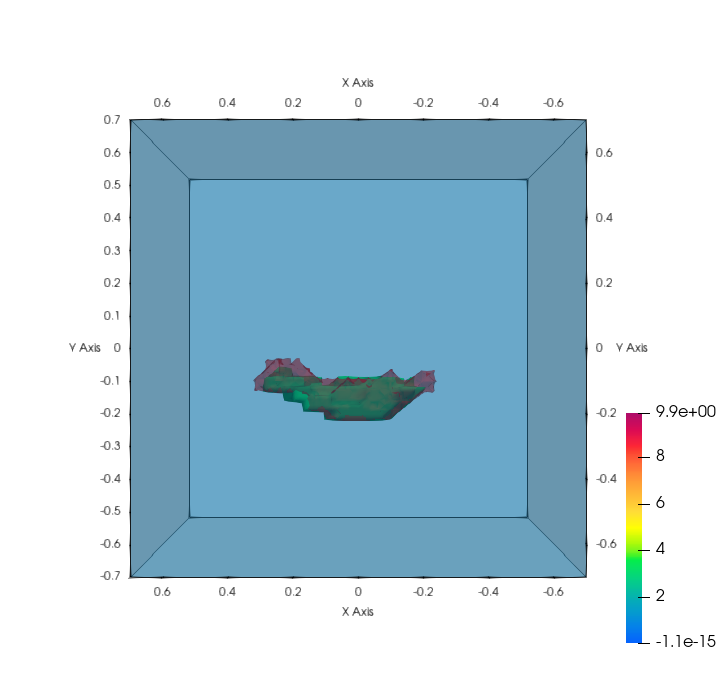}} &
      {\includegraphics[scale=0.16, trim = 0.0cm 0.0cm 0cm 0.0cm, clip=true,]{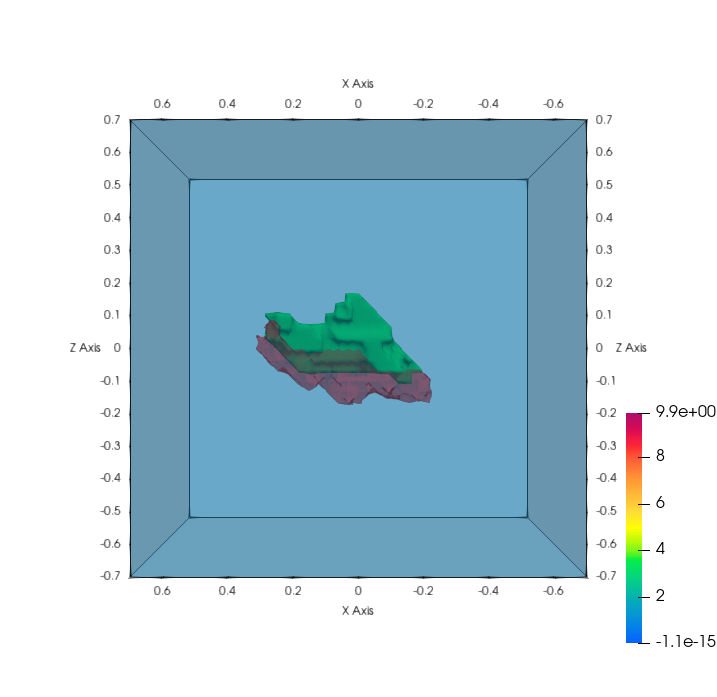}}\\
     a) prospect view & b) $x_1 x_2$ view & c) $x_1 x_3$ view\\
      \\
      \hline
       {\includegraphics[scale=0.16, trim = 0.0cm 0.0cm 0cm 0.0cm, clip=true,]{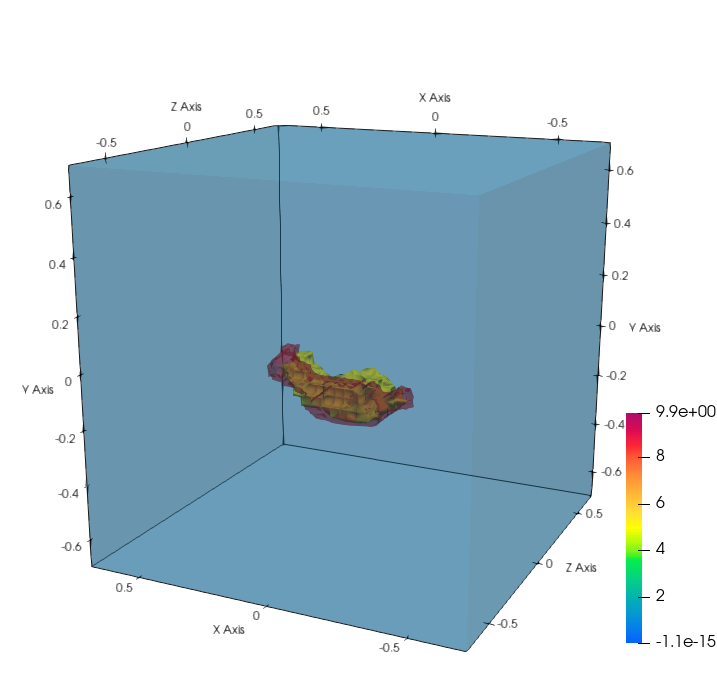}}   &
      {\includegraphics[scale=0.16, trim = 0.0cm 0.0cm 0cm 0.0cm, clip=true,]{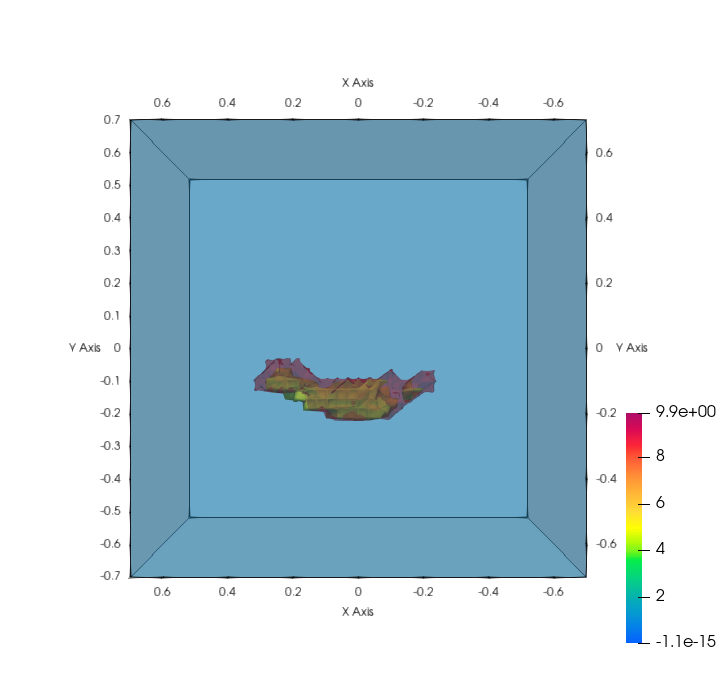}} &
      {\includegraphics[scale=0.16, trim = 0.0cm 0.0cm 0cm 0.0cm, clip=true,]{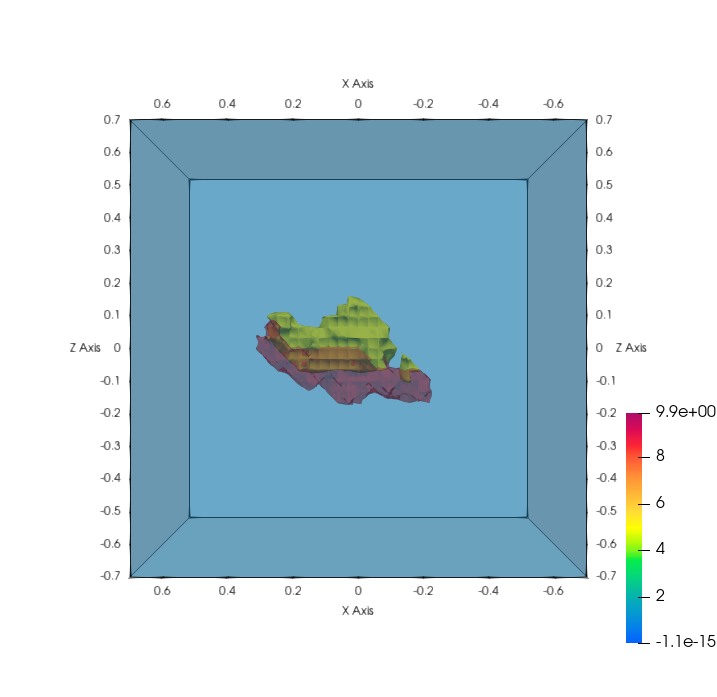}}\\
   d)   prospect view &  e) $x_1 x_2$ view & f) $x_1 x_3$ view\\
        \hline
       {\includegraphics[scale=0.16, trim = 0.0cm 0.0cm 0cm 0.0cm, clip=true,]{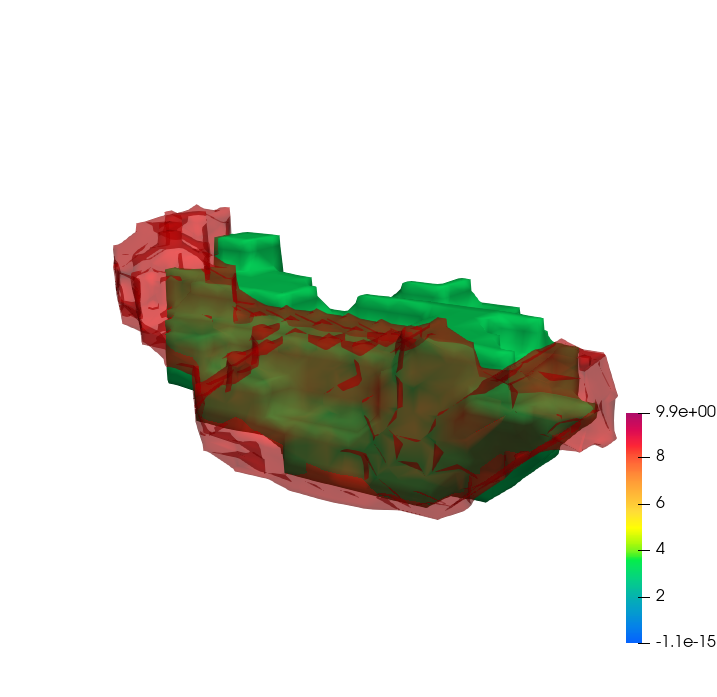}}   &
      {\includegraphics[scale=0.16, trim = 0.0cm 0.0cm 0cm 0.0cm, clip=true,]{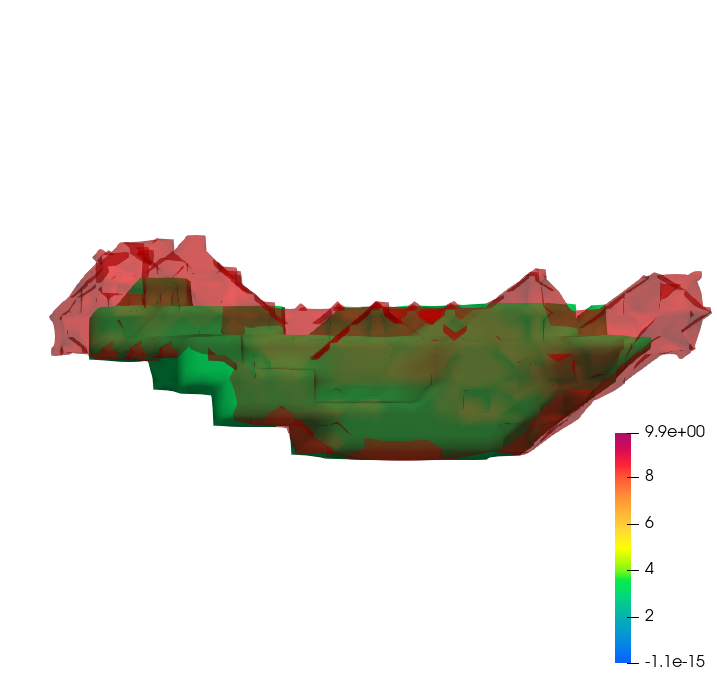}} &
      {\includegraphics[scale=0.16, trim = 0.0cm 0.0cm 0cm 0.0cm, clip=true,]{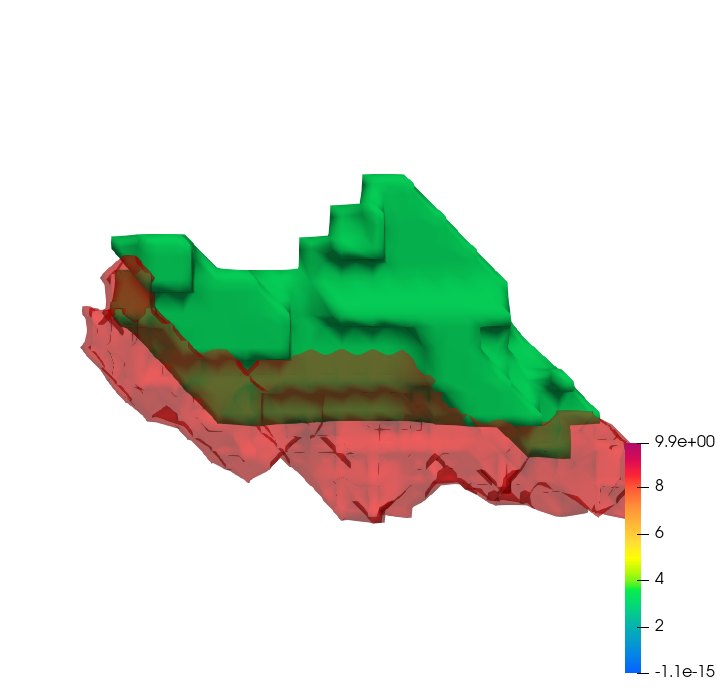}}\\
    g)  zoomed prospect view & h)  zoomed $x_1 x_2$ view & i) zoomed $x_1 x_3$ view\\
      \\
      \hline
      \\
      {\includegraphics[scale=0.13, trim = 5.0cm 0.0cm 6.0cm 0.0cm, clip=true,]{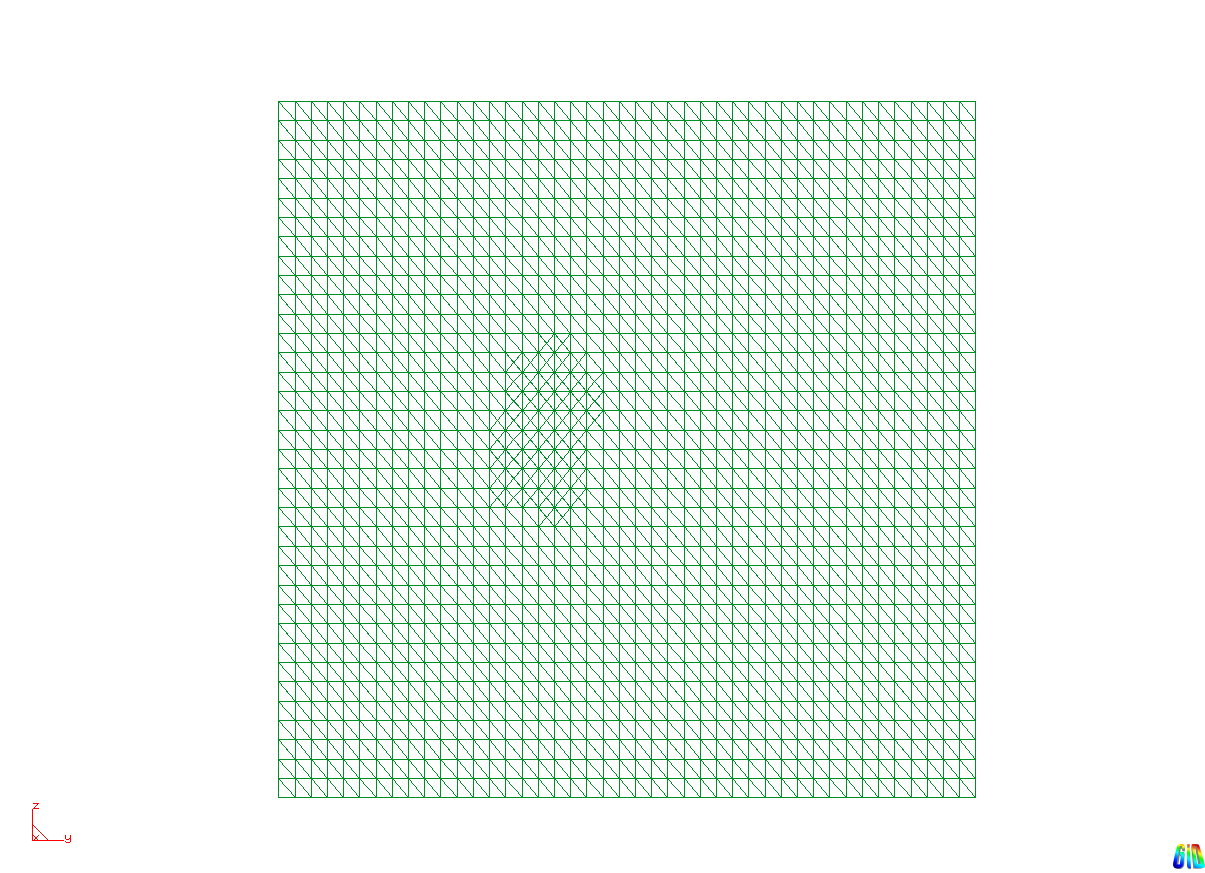}}
          &
      {\includegraphics[scale=0.13, trim = 5.0cm 0.0cm 6.0cm 0.0cm, clip=true,]{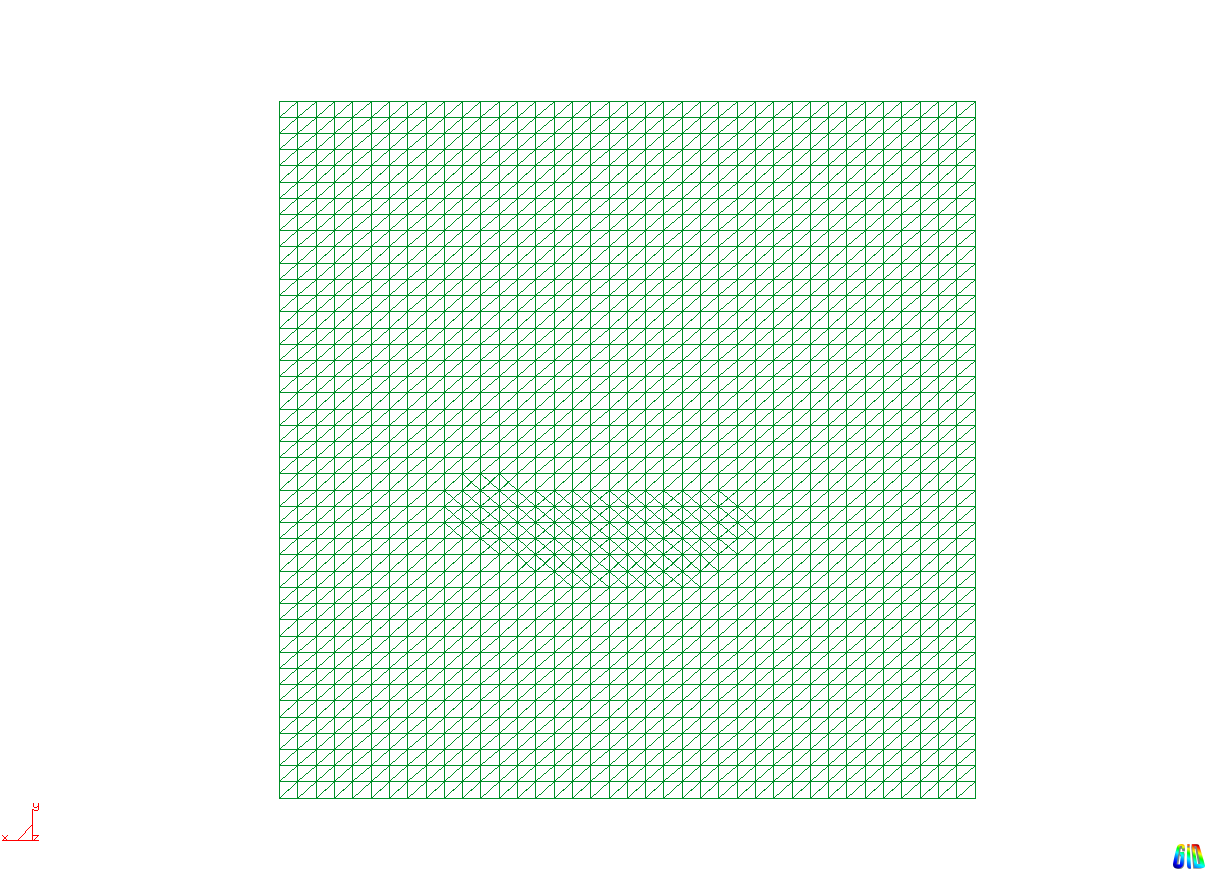}} &
      {\includegraphics[scale=0.13, trim = 5.0cm 0.0cm 6.0cm 0.0cm, clip=true,]{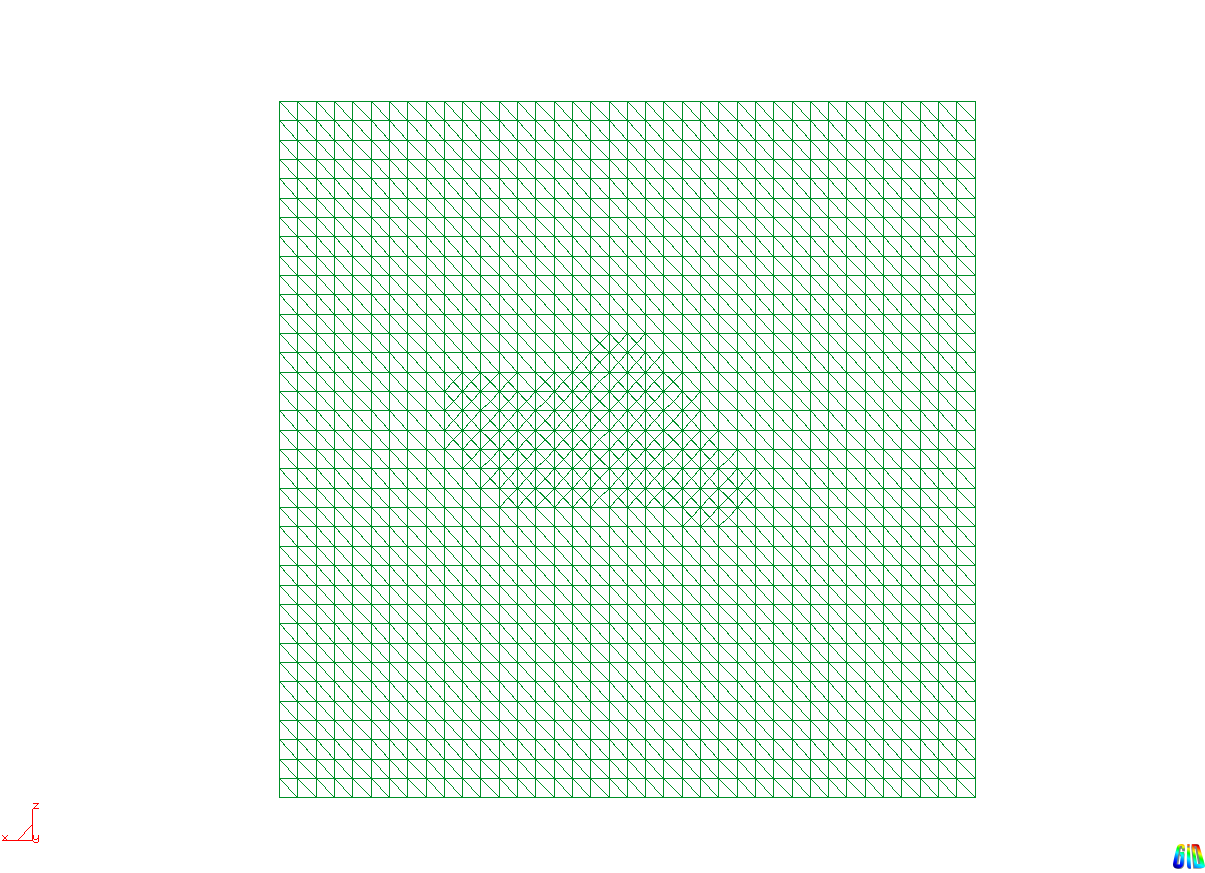}} \\
    j)   $x_2 x_3$ view & k) $x_1 x_2$ view & l) $x_1 x_3$ view \\
      \hline
      \\
    \end{tabular}
  \end{center}
  \caption{
    \small\emph{Test 2. a) -c): Isosurfaces of reconstructions
    $\varepsilon_{h_0} \approx 5$ (in green color) of
    $\varepsilon_r$ obtained on the coarse mesh ${K_h}_0$. d)-f): Isosurfaces of reconstructions ${\varepsilon_h}_1 \approx 5$
    obtained on refined mesh ${K_h}_1$ (in yellow color). g)-i) Zoomed reconstructions. 
  j) -l): Refined mesh ${K_h}_1$. The noise level in the data
    is $\delta= 10\%$. 
     See Table \ref{tab:test2}  for obtained contrasts  $\max_{\Omega_{\rm FEM}}  {\varepsilon_h}_k, k=0,1$.
   For comparison we also present exact isosurface  of $\varepsilon_r$  with value corresponding to reconstructed one
   and outlined by red color.
}}
  \label{fig:test2noise10}
\end{figure}

 \subsection{Test 1}

\label{sec:test1}

In this test we present numerical results of reconstruction of
  $\varepsilon_r$ when exact values of this function are  given in Table
\ref{tab:table1}, see Test 1.  Isosurface of the exact
function $\varepsilon_r$ to be reconstructed in this test is shown
in Figure ~\ref{fig:exacteps}-a).  We note
that the exact function  $\varepsilon_r$ has  complicated   structure.
 Using Figure ~\ref{fig:exacteps}-a) one can observe that isosurface presents a discontinuous function with a lot of big and
small inclusions in the domain  $\Omega_{\rm FEM}$.

  Figures  \ref{fig:test1noise10v2}-a)-i)
  show results of the reconstruction on
adaptively  locally refined meshes when
noise  level in the data was $\delta = 10\%$.
We start computations on a coarse mesh ${K_h}_0$.
  Figure  \ref{fig:test1noise10v2}-a)-c)
  shows that
the location of the reconstructed function
 ${\varepsilon_h}_0$ 
is imaged correctly and the reconstructed isosurface  covers the domain where the exact $\varepsilon_r$ is located.
 We refer to Table  \ref{tab:test1} for the
 reconstruction of the maximal contrast in ${\varepsilon_h}_0$.
  For improvement of the contrast  and shape obtained on a coarse mesh  ${K_h}_0$, we run
  computations on locally adaptivelly refined meshes.
  Figures
  \ref{fig:test1noise10v2}-d)-f)  show reconstruction
  obtained on the final two times refined mesh ${K_h}_2$.
   Table \ref{tab:test1} presents results of
  reconstructions for ${\varepsilon_h}_k$ obtained on the refined meshes  ${K_h}_k,  k=0, 1,2$. We observe that with
  mesh refinements we achieve better contrast for function
  $\varepsilon_r$.
  Also reconstructed isosurface of this function
  more precisely covers the domain where the exact $\varepsilon_r$ is located, compare Figure  \ref{fig:test1noise10v2}-a) with Figure  \ref{fig:test1noise10v2}-d). Figures  \ref{fig:test1noise10v2}-g)-i)  show locally adaptively refined mesh ${K_h}_2$.

\begin{table}[tbp] 
\vspace{2mm}
\centerline{
  \begin{tabular}{|c|c|}
     \hline
     \multicolumn{2}{|c|}
                 {Test 2}
                 \\
 \hline
   $\delta=3\%$ &  $\delta = 10\%$ 
 \\
\begin{tabular}{l|l|l|l} \hline
Mesh & $\max_{\Omega_{\rm FEM}} {\varepsilon_h}_k$ & $\frac{\max_{\Omega_{\rm FEM}} |\varepsilon_r - {\varepsilon_h}_k | }{\max_{\Omega_{\rm FEM}} |\varepsilon_r|}$ & $M^k$    \\ \hline
${K_h}_0$ & 6.874  & 0.236 & 2   \\
${K_h}_1$ & 7.558 &  0.160 & 5 \\
${K_h}_2$ & 10.0  &  0.111 & 2 \\
\end{tabular}
 & 
\begin{tabular}{l|l|l|l} \hline
  Mesh & $\max_{\Omega_{\rm FEM}} {\varepsilon_h}_k$ &    
$\frac{\max_{\Omega_{\rm FEM}} |\varepsilon_r - {\varepsilon_h}_k | }{\max_{\Omega_{\rm FEM}} |\varepsilon_r|}$ & $M^k$  \\ \hline
 ${K_h}_0$  & 5.350  &  0.406  &  2    \\
${K_h}_1$   & 9.450  &  0.05 &   4 \\
            &   & & \\
\end{tabular} \\
\hline
\end{tabular}}
\caption{
  \textit{ Test 2. Computational results of the
    reconstructions $\max_{\Omega_{\rm FEM}}  {\varepsilon_h}_k$ on a coarse and on adaptively refined meshes together with  relative 
    errors computed in the maximal contrast of $ \max_{\Omega_{\rm FEM}}  \varepsilon_r,\max_{\Omega_{\rm FEM}}  {\varepsilon_h}_k$. Here,
    $\max_{\Omega_{\rm FEM}}  {\varepsilon_h}_k$ denotes  maximum of the computed function $\varepsilon_h$
    on $k$ times refined mesh ${K_h}_k$ in the domain $\Omega_{\rm FEM}$, and
    $M^k$ denotes the final number of
    iterations in the conjugate gradient Algorithm 3 on $k$ times refined mesh ${K_h}_k$  for reconstructed
    function ${\varepsilon_h}_k, k=0,1,2$.
}}
    \label{tab:test2}
\end{table}


\begin{table}[tbp] 
\vspace{2mm}
\centerline{
  \begin{tabular}{|c|c|}
     \hline
     \multicolumn{2}{|c|}
                 {Test 2, Computational Time}
                 \\
 \hline
   $\delta=3\%$ &  $\delta = 10\%$ 
 \\
\begin{tabular}{l|l|l|l|l} \hline
Mesh & nno &  Time (sec)& Rel. time  & $M^k$    \\ \hline
${K_h}_0$ & 63492  & 1186  & 3.72 $\cdot 10^{-5}$ & 2   \\
${K_h}_1$ &  64096  & 3588 & 9.34 $\cdot  10^{-5}$ & 5 \\
${K_h}_2$ &  66112 &  1228   & 3.72 $\cdot  10^{-5}$ &  2 \\
\end{tabular}
 & 
\begin{tabular}{l|l|l|l|l} \hline
  Mesh &  nno  &     Time (sec) & Rel. time  & $M^k$  \\ \hline
 ${K_h}_0$  &  63492  & 1214  &  $3.82 \cdot 10^{-5}$    &  2    \\
${K_h}_1$   &  63968  & 2384  &  $7.44 \cdot 10^{-5}$     &  4 \\
  & & &  & \\
\end{tabular} \\
\hline
\end{tabular}}
\caption{\textit{
    Test 2.  Performance of the reconstruction Algorithm 4
    (in seconds)
  on adaptively refined meshes. Here,
  $k$ is number of the refined mesh ${K_h}_k$ of the domain $\Omega_{\rm FEM}$,
 ${\rm nno}$ is number of the nodes in the computational   mesh  ${K_h}_k$,
  and
    $M^k$ denotes the final number of
  iterations in the conjugate gradient Algorithm 3.
}}
    \label{tab:test2time}
\end{table}


 \subsection{Test 2}

\label{sec:test2}

Since it is quite demanding  reconstruct very complicated structure
of $\varepsilon_r$  taken in Test 1, in this test we will
reconstruct $\varepsilon_r$ with exact isosurface as it is presented in the
Figure \ref{fig:exacteps}-b).
Exact values of this function are
taken as in fibroconnective/glandular-1 media (see Table \ref{tab:table1})
inside isosurface of Figure
\ref{fig:exacteps}-b),  and  outside of this isosurface all values of
$\varepsilon_r =1$.

Figures \ref{fig:test2noise10}-a)-i)
show results of the reconstruction on
adaptively refined meshes when
noise  level in the data was $\delta = 10\%$. We refer to the Table 
 \ref{tab:test2}  for 
  reconstruction of the contrast in $\varepsilon_r$.
  Using the Table  \ref{tab:test2}
  we now observe that with mesh refinements we achieve
slightly higher  maximal contrast $9.45$  in reconstruction ${\varepsilon_h}_1$ compared to the exact one $9$.
Moreover, on the mesh ${K_h}_1$  for $\sigma = 10\%$ we get more than 8 times smaller relative error in the reconstruction compared to the error obtained on  the coarse mesh  ${K_h}_0$.
  Figures \ref{fig:test2noise10}-d)-i)  show   good
  matching
   of the reconstructed  ${\varepsilon_h}_1$ compared with the exact one.
Figures  \ref{fig:test1noise10v2}-j)-l)  show locally adaptively refined mesh ${K_h}_2$.

\begin{table}[tbp] 
\vspace{2mm}
\centerline{
   \begin{tabular}{|c|c|}
     \hline
      \multicolumn{2}{|c|}
                 { Computational time}
                 \\
 \hline
   $\delta=3\%$ &  $\delta = 10\%$ 
 \\
\begin{tabular}{l|l|l|l} \hline
         &  Time (sec) &  Relative time            & $n$     \\ \hline
  Test 1 & 110.59      & $1.779\cdot 10^{-6}$     & 71360 \\
  Test 2 & 106.58      &  $1.714 \cdot 10^{-6}$ & 69699  \\
\end{tabular}
 & 
\begin{tabular}{l|l|l|l} \hline
 &     Time (sec)       &  Relative time & $n$  \\ \hline
 Test 1    &    116.22      &  $1.869\cdot 10^{-6}$      &  75052    \\
 Test 2    &    111.53      &  $1.793\cdot 10^{-6}$      &   65359 \\
\end{tabular} \\
\hline
\end{tabular}}
\caption{\textit{
    Performance of solution of forward problem  \eqref{model1}  in Tests 1 and 2  of section \ref{sec:numex} on the mesh ${K_h}_0$   in terms of computational time (in seconds) and relative computational time computed by \eqref{relt1}.
    Here, $n$ is number of the nodes  on three times adaptively refined original coarse mesh
 (consisting of 63492 nodes)
    which we used for generation of transmitted data.
}}
    \label{tab:time}
\end{table}

\subsection{Performance comparison }

All computations were performed on a linux
workstation Intel Core i7-9700 CPU with one processor using software
package WavES \cite{waves} efficiently implemented in C++/PETSc \cite{petsc}.

We have estimated
the relative computational time $T_r$ of the forward problem using the following formula
  \begin{equation}\label{relt1}
T_r = \frac{t}{  n_t \cdot n}.
  \end{equation}
  Here, $t$ is the total computational time  of the forward problem on
  the mesh ${K_h}_l$ where $l =0,1,2,...$ is number of the refined mesh, $n$ is
  the total number of nodes in the mesh ${K_h}_l$, $n_t$ is number of
  timesteps.
  We take $n_t = 500$  in all computational tests, see clarification in section \ref{sec:setup}.
  Computational times (in seconds) for solution of forward
  problem are presented in Table \ref{tab:time}.  Using this table we
  observe that the relative time is approximately the same for all
  tests and we can take it as $T_r \approx 1.8 \cdot
  10^{-6}$. Next, using this relative time we can estimate approximate
  computational time for solution of forward problem for any mesh
  consisting of $n$ nodes. For example, if we will take original mesh
  consisting of $n= 34 036 992$ nodes, then computational time will be
  already $t = T_r \cdot n_t \cdot n = 1.8 \cdot 10^{-6} \cdot 500
  \cdot 34 036 992 = 30633 $ seconds, and this time is not
  computationally efficient. Clearly, computing of the solution of
  inverse problem on the sampled mesh allows significantly reduce
  computational times.

We have estimated  also
the relative computational time $T_r^{ip}$ of the   solution of inverse problem
using the formula
  \begin{equation}\label{relt2}
T_r^{ip} = \frac{t^{ip}}{  n_t \cdot nno}.
  \end{equation}
  Here, $t^{ip}$ is the total computational time to run inverse Algorithm 4
  on the mesh ${K_h}_l$ where $l =0,1,2,...$ is number of the refined
  mesh, $nno$ is the total number of nodes in the mesh ${K_h}_l$,
  $n_t$ is number of timesteps.  Computational times (in seconds) for
  solution of inverse problem for Test 1 and Test 2 are presented in
  Tables \ref{tab:test1time},\ref{tab:test2time}, respectivelly.
  Using these tables we observe that computational times are depend on
  the number of iterations
 $M^k$
  in the conjugate gradient method  (CGM)  and
  number of the nodes $nno$ in the meshes ${K_h}_l$. We took $n_t=500$ for
  all tests and thus, computational times presented in these tables
  are not depend on number of times steps for different refined
  meshes. We note, that the number of time steps $n_t$ can be chosen
  adaptively as well. However, we are performing only adaptive mesh
  refinement in space and not in time.  The full space-time adaptive
  algorithm can be considered as a topic for future research.

  Using Table \ref{tab:test1time} we observe that computational time
  in Test 1 is
  around 20 minutes for both noise levels  $\sigma = 3\%$ and $ \sigma =10\%$.
  On every mesh ${K_h}_l, l=0,1,2$, was performed
  two iterations CGM , or $M^K=2$.  Thus, the total computational time
  to obtain final reconstruction in Test 1 is 60 min.
  
   Table \ref{tab:test2time}  shows that 
   computational time in Test 2 with noise in data $\delta = 3\%$ is
   around 20 minutes for non-refined mesh ${K_h}_0$, 60 min for one
   time refined mesh ${K_h}_1$, and 20 minutes for twice refined mesh
   ${K_h}_2$.  Thus, the total computational time to obtain final
   reconstruction in Test 2 is 100 minutes. Computational time in this test
   is larger than in the previous Test 1 since  CGM converged only at 5-th iteration
    on the one time refined mesh ${K_h}_1$.
   However, the total
   computational time with noise in data $\delta = 10 \%$ is around 60
   minutes. This is because the  solution was obtained already  on the one time refined mesh
   ${K_h}_1$.
   Tables
 \ref{tab:test1time},  \ref{tab:test2time}
   also demonstrate that it takes  around 10 minutes to compute solution of inverse problem on the one iteration of the conjugate gradient algorithm.

  We note that
  PETSc supports parallel implementation and thus, current version of
  code can be extended to the version with parallel implementation
  such that times reported in Tables \ref{tab:time},  \ref{tab:test1time}
   and Table 
 \ref{tab:test2time} 
  can be significantly reduced.
 
\section{Conclusions}

\label{sec:concl}

This work describes reconstruction methods for determination of the
relative dielectric permittivity function in conductive media using
scattered data of the time-dependent electric field at number of
detectors placed at the boundary of the investigated domain.

Reconstruction methods use optimization approach where a functional
 is minimized via a domain decomposition finite element/finite
 difference method.  In an adaptive reconstruction method the space
 mesh is refined only in the domain where a finite element method is
 used with a feedback from a posteriori error indicators.  Developed
 adaptive algorithms allow
 us to
 obtain correct values and shapes of the
 dielectric permittivity function to be determined.
Convergence and stability analysis of the developed methods is ongoing work and will be presented in forthcoming publication.
 The algorithms of the current work are designed from previous
 adaptive algorithms  developed in  \cite{BTKM2, BondestaB}
  which reconstruct  the wave speed or the
 dielectric permittivity function.
 However, all previous algorithms
 are developed for non-conductive medium.
 
Our computational tests show qualitative and quantitative
reconstruction of dielectric permittivity function using anatomically
realistic breast phantom which capture the heterogeneity of normal
breast tissue at frequency 6 GHz taken from online repository
\cite{wisconsin}.
In all tests we used assumption that the  conductivity function is known. Currently we are working on algorithms when both dielectric permittivity and conductivity functions   can be reconstructed. Results of this work
will be presented in  our future research.

All computations are performed in real time presented in Tables
\ref{tab:test1time}, \ref{tab:test2time} and \ref{tab:time}.  Some
data (Matlab code to read data of database \cite{wisconsin}, visualize
 and produce discretized values of $\varepsilon_r, \sigma$,
etc.)  used in computations of this work is available for download and
testing, see \cite{projectwaves}.  Additional data (computational
FE/FD meshes, transmitted data, C++/PETSc code) can be provided upon
request.

In summary, the main features of algorithms of this work are as
follows:

\begin{itemize}

\item Ability to reconstruct shapes, locations and maximal values of
  dielectric permittivity function of targets in conductive media
    under the condition that the conductivity of this media is a known function.

  \item More exact reconstruction of   shapes and
     maximal values of
  dielectric permittivity function
  of inclusions because of local adaptive
  mesh refinement.

\item Computational greater efficiency  
  because of usage  software package WavES \cite{waves} implemented in C++/PETSc \cite{petsc}.

\end{itemize}

\section*{Acknowledgment}
The research of authors
is supported by the Swedish Research Council grant VR 2018-03661.

\end{document}